\documentclass[1p,times]{elsarticle}

\usepackage[utf8]{inputenc}
\usepackage{amsmath}
\usepackage{dutchcal}
\usepackage{amsfonts}
\usepackage{amssymb}
\usepackage{amsthm}
\usepackage[dvipsnames,table]{xcolor}
\usepackage{graphicx}
\usepackage{todonotes}
\usepackage{tikz}
\usetikzlibrary{spy}
\usetikzlibrary{shapes.geometric}
\usetikzlibrary{pgfplots.groupplots}
\pgfdeclarelayer{bg}    
\pgfsetlayers{bg,main}  

\usepackage{pgfplots}
\usepackage{pgfplotstable} 
\usepackage{subcaption}
\usepackage{algorithm}
\usepackage{algorithmicx}
\usepackage{algpseudocode}

\usepackage{booktabs}
\usepackage{multirow}
\usepackage{bm}
\usepackage{mathtools}
\usepackage{enumerate}
\usepackage{hyperref}
\hypersetup{
    colorlinks = true,
}
\usepackage[capitalize]{cleveref}

\usepackage{physics}

\newcommand{\HV}[1]{\textcolor{Green}{#1}}

\usepackage{comment}

\usepackage{lineno}

\biboptions{sort&compress}

\journal{}

\definecolor{tudelft-cyan}{cmyk}{1,0,0,0}
\definecolor{tudelft-black}{cmyk}{0,0,0,1}
\definecolor{tudelft-white}{cmyk}{0,0,0,0}
\definecolor{tudelft-dark-blue}{cmyk}{1,0.8,0.08,0.7}
\definecolor{tudelft-turquoise}{cmyk}{0.72,0,0.24,0}
\definecolor{tudelft-blue}{cmyk}{0.98,0.4,0.0,0}
\definecolor{tudelft-purple}{cmyk}{0.65,1,0,0.12}
\definecolor{tudelft-pink}{cmyk}{0,0.7,0,0}
\definecolor{tudelft-raspberry}{cmyk}{0.05,1,0.48,0.3}
\definecolor{tudelft-red}{cmyk}{0,0.85,0.75,0}
\definecolor{tudelft-orange}{cmyk}{0,0.7,0.75,0}
\definecolor{tudelft-yellow}{cmyk}{0,0.31,0.98,0}
\definecolor{tudelft-light-green}{cmyk}{0.63,0,0.84,0}
\definecolor{tudelft-dark-green}{cmyk}{1,0,0.68,0.04}

\colorlet{TUcol1}{tudelft-cyan}
\colorlet{TUcol2}{tudelft-black}
\colorlet{TUcol3}{tudelft-white}
\colorlet{TUcol4}{tudelft-dark-blue}
\colorlet{TUcol5}{tudelft-turquoise}
\colorlet{TUcol6}{tudelft-blue}
\colorlet{TUcol7}{tudelft-purple}
\colorlet{TUcol8}{tudelft-pink}
\colorlet{TUcol9}{tudelft-raspberry}
\colorlet{TUcol10}{tudelft-red}
\colorlet{TUcol11}{tudelft-orange}
\colorlet{TUcol12}{tudelft-yellow}
\colorlet{TUcol13}{tudelft-light-green}
\colorlet{TUcol14}{tudelft-dark-green}

\colorlet{gray1}{black!20}
\colorlet{gray2}{black!40}
\colorlet{gray3}{black!60}
\colorlet{gray4}{black!80}

\colorlet{col0}{TUcol4}
\colorlet{col1}{TUcol9}
\colorlet{col2}{TUcol12}
\colorlet{col3}{TUcol13}
\colorlet{col4}{TUcol6}
\colorlet{col5}{TUcol7}

\colorlet{plotcol1}{TUcol10}
\colorlet{plotcol2}{TUcol12}
\colorlet{plotcol3}{TUcol13}
\colorlet{plotcol4}{TUcol14}
\colorlet{plotcol5}{TUcol6}
\colorlet{plotcol6}{TUcol7}

\pgfplotscreateplotcyclelist{mark list mod}{
every mark/.append style={solid,fill=\pgfplotsmarklistfill,line width = 0.1},mark size=1.0,mark=*,mark options=solid,fill opacity=0.5\\
every mark/.append style={solid,fill=\pgfplotsmarklistfill,line width = 0.1},mark size=1.0,mark=square*,mark options=solid,fill opacity=0.5\\
every mark/.append style={solid,fill=\pgfplotsmarklistfill,line width = 0.1},mark size=1.0,mark=triangle*,mark options=solid,fill opacity=0.5\\
every mark/.append style={solid,fill=\pgfplotsmarklistfill,line width = 0.1},mark size=1.0,mark=halfsquare*,mark options=solid,fill opacity=0.5\\
every mark/.append style={solid,fill=\pgfplotsmarklistfill,line width = 0.1},mark size=1.0,mark=pentagon*,mark options=solid,fill opacity=0.5\\
every mark/.append style={solid,fill=\pgfplotsmarklistfill,line width = 0.1},mark size=1.0,mark=halfcircle*,mark options=solid,fill opacity=0.5\\
every mark/.append style={solid,fill=\pgfplotsmarklistfill,rotate=180,line width = 0.1},mark size=1.0,mark=halfdiamond*,mark options=solid,fill opacity=0.5\\
every mark/.append style={solid,fill=\pgfplotsmarklistfill!40,line width = 0.1},mark size=1.0,mark=otimes*,mark options=solid,fill opacity=0.5\\
every mark/.append style={solid,fill=\pgfplotsmarklistfill,line width = 0.1},mark size=1.0,mark=diamond*,mark options=solid,fill opacity=0.5\\
every mark/.append style={solid,fill=\pgfplotsmarklistfill,line width = 0.1},mark size=1.0,mark=halfsquare right*,mark options=solid,fill opacity=0.5\\
every mark/.append style={solid,fill=\pgfplotsmarklistfill,line width = 0.1},mark size=1.0,mark=halfsquare left*,mark options=solid,fill opacity=0.5\\
}

\pgfplotscreateplotcyclelist{linestyles mod}{
thin,solid\\
thin,densely dashed\\
thin,densely dotted\\
thin,densely dashdotted\\
thin,densely dashdotdotted\\
thin,dashed\\
thin,dotted\\
thin,dashdotted\\
thin,dashdotdotted\\
}

\pgfplotscreateplotcyclelist{colorlist mod}{
plotcol1\\
plotcol2\\
plotcol3\\
plotcol4\\
plotcol5\\
plotcol6\\
}

\pgfplotsset{
	cycle multiindex* list={
							mark list mod\nextlist
							colorlist mod\nextlist
							linestyles mod\nextlist
							},
    width =\linewidth,
    height=\normalplotheight,
    grid=major,
    ylabel near ticks,
    xlabel near ticks,
    enlargelimits,
    grid=major,
	legend style={fill=white, fill opacity=0.6, draw opacity=1,draw=black,text opacity=1, font=\small},
    legend cell align={left},
    legend image with text/.style={
    legend image code/.code={%
        \node[anchor=center] at (0.3cm,0cm) {#1};
    }
	},
}

\pgfplotsset{select coords between index/.style 2 args={
		x filter/.code={
			\ifnum\coordindex<#1\fi
			\ifnum\coordindex>#2\fi
		}
}}

\tikzset{%
style0/.style = {plotcol1,thin,solid,every mark/.append style={solid,fill=\pgfplotsmarklistfill,line width = 0.1},mark size=1.0,mark=*,mark options=solid,fill opacity=0.75,},
style1/.style = {plotcol2,thin,densely dashed,every mark/.append style={solid,fill=\pgfplotsmarklistfill,line width = 0.1},mark size=1.0,mark=square*,mark options=solid,fill opacity=0.75},
style2/.style = {plotcol3,thin,densely dotted,every mark/.append style={solid,fill=\pgfplotsmarklistfill,line width = 0.1},mark size=1.0,mark=triangle*,mark options=solid,fill opacity=0.75},
style3/.style = {plotcol4,thin,densely dashdotted,every mark/.append style={solid,fill=\pgfplotsmarklistfill,line width = 0.1},mark size=1.0,mark=halfsquare*,mark options=solid,fill opacity=0.75},
style4/.style = {plotcol5,thin,densely dashdotdotted,every mark/.append style={solid,fill=\pgfplotsmarklistfill,line width = 0.1},mark size=1.0,mark=pentagon*,mark options=solid,fill opacity=0.75},
}

\tikzset{%
singlepatch/.style 	= {style=style0,thick,black,mark=0},
dpatch/.style 		= {style=style0,thick},
dpatch2/.style 		= {style=style0,thick,plotcol1!50!black},
approxC1/.style 	= {style=style1,thick},
exactC1/.style 		= {style=style2,thick},
almostC1/.style 	= {style=style3,thick},
weak1/.style 		= {style=style4,thick,plotcol5!75!black},
weak2/.style 		= {style=style4,thick,plotcol5!50!black},
weak3/.style 		= {style=style4,thick,plotcol5!25!black},
weak4/.style 		= {style=style4,thick,plotcol5!20!black}
}

\newcommand{\logLogSlopeTriangle}[5]
{
    \small

    \pgfplotsextra
    {
        \pgfkeysgetvalue{/pgfplots/xmin}{\xmin}
        \pgfkeysgetvalue{/pgfplots/xmax}{\xmax}
        \pgfkeysgetvalue{/pgfplots/ymin}{\ymin}
        \pgfkeysgetvalue{/pgfplots/ymax}{\ymax}

        \pgfmathsetmacro{\xArel}{#1}
        \pgfmathsetmacro{\yArel}{#3}
        \pgfmathsetmacro{\xBrel}{#1-#2}
        \pgfmathsetmacro{\yBrel}{\yArel}
        \pgfmathsetmacro{\xCrel}{\xArel}

        \pgfmathsetmacro{\lnxB}{\xmin*(1-(#1-#2))+\xmax*(#1-#2)} 
        \pgfmathsetmacro{\lnxA}{\xmin*(1-#1)+\xmax*#1} 
        \pgfmathsetmacro{\lnyA}{\ymin*(1-#3)+\ymax*#3} 
        \pgfmathsetmacro{\lnyC}{\lnyA+#4*(\lnxA-\lnxB)}
        \pgfmathsetmacro{\yCrel}{\lnyC-\ymin)/(\ymax-\ymin)} 

        \coordinate (A) at (rel axis cs:\xArel,\yArel);
        \coordinate (B) at (rel axis cs:\xBrel,\yBrel);
        \coordinate (C) at (rel axis cs:\xCrel,\yCrel);

        \draw[#5]   (A)-- node[pos=0.5,anchor=north] {1}
                    (B)--
                    (C)-- node[pos=0.5,anchor=west] {#4}
                    cycle;
    }
}

\newcommand{\logLogSlopeTriangleRev}[5]
{
\small

    \pgfplotsextra
    {
        \pgfkeysgetvalue{/pgfplots/xmin}{\xmin}
        \pgfkeysgetvalue{/pgfplots/xmax}{\xmax}
        \pgfkeysgetvalue{/pgfplots/ymin}{\ymin}
        \pgfkeysgetvalue{/pgfplots/ymax}{\ymax}

        \pgfmathsetmacro{\xArel}{#1-#2}
        \pgfmathsetmacro{\yArel}{#3}
        \pgfmathsetmacro{\xBrel}{#1}
        \pgfmathsetmacro{\yBrel}{\yArel}
        \pgfmathsetmacro{\xCrel}{\xArel}

        \pgfmathsetmacro{\lnxB}{\xmin*(1-(#1-#2))+\xmax*(#1-#2)} 
        \pgfmathsetmacro{\lnxA}{\xmin*(1-#1)+\xmax*#1} 
        \pgfmathsetmacro{\lnyA}{\ymin*(1-#3)+\ymax*#3} 
        \pgfmathsetmacro{\lnyC}{\lnyA-#4*(\lnxA-\lnxB)}
        \pgfmathsetmacro{\yCrel}{\lnyC-\ymin)/(\ymax-\ymin)} 

        \coordinate (A) at (rel axis cs:\xArel,\yArel);
        \coordinate (B) at (rel axis cs:\xBrel,\yBrel);
        \coordinate (C) at (rel axis cs:\xCrel,\yCrel);

        \draw[#5]   (A)-- node[pos=0.5,anchor=south] {1}
                    (B)--
                    (C)-- node[pos=0.5,anchor=east] {#4}
                    cycle;
    }
}

\pgfkeys{
	/pgfplots/table/print last/.style={
		row predicate/.code={
			\pgfmathtruncatemacro\firstprintedrownumber{\pgfplotstablerows-#1}
			\ifnum##1<\firstprintedrownumber\relax
			\pgfplotstableuserowfalse
			\fi
		}
	},
	/pgfplots/table/print last/.default=1
}

\theoremstyle{definition}

\crefname{rem}{Remark}{Remarks}

\crefname{ex}{Example}{Examples}

\begin{document}
\begin{frontmatter}



\title{A comparison of smooth basis constructions for isogeometric analysis}

\author[label1,label2]{H.M. Verhelst}
\ead{h.m.verhelst@tudelft.nl}
\author[label3]{P. Weinmüller}
\ead{pascal.weinmueller@mtu.de}
\author[label4]{A. Mantzaflaris}
\ead{angelos.mantzaflaris@inria.fr}
\author[label5]{T. Takacs}
\ead{thomas.takacs@ricam.oeaw.ac.at}
\author[label1]{D. Toshniwal}
\ead{d.toshniwal@tudelft.nl}
\corref{cor1}
\cortext[cor1]{Corresponding Author. }
\address[label1]{Delft University of Technology, Department of Applied Mathematics, Mekelweg 4, Delft 2628 CD, The Netherlands}
\address[label2]{Delft University of Technology, Department of Maritime and Transport Technology, Mekelweg 2, Delft 2628 CD, The Netherlands}
\address[label3]{MTU Aero Engines AG, Dachauer Straße 665, 80995, Munich, Germany}
\address[label4]{Inria centre at Université Côte d’Azur, 2004 route des Lucioles - BP 93, 06902 Sophia Antipolis, France}
\address[label5]{Johann Radon Institute for Computational and Applied Mathematics, Austrian Academy of Sciences, Altenberger Str. 69, 4040 Linz, Austria}

\begin{abstract}
In order to perform isogeometric analysis with increased smoothness on complex domains, trimming, variational coupling or unstructured spline methods can be used. The latter two classes of methods require a multi-patch segmentation of the domain, and provide continuous bases along patch interfaces. In the context of shell modeling, variational methods are widely used, whereas the application of unstructured spline methods on shell problems is rather scarce. In this paper, we therefore provide a qualitative and a quantitative comparison of a selection of unstructured spline constructions, in particular the D-Patch, Almost-$C^1$, Analysis-Suitable $G^1$ and the Approximate $C^1$ constructions. Using this comparison, we aim to provide insight into the selection of methods for practical problems, as well as directions for future research. In the qualitative comparison, the properties of each method are evaluated and compared. In the quantitative comparison, a selection of numerical examples is used to highlight different advantages and disadvantages of each method. In the latter, comparison with weak coupling methods such as Nitsche's method or penalty methods is made as well. In brief, it is concluded that the Approximate $C^1$ and Analysis-Suitable $G^1$ converge optimally in the analysis of a bi-harmonic problem, without the need of special refinement procedures. Furthermore, these methods provide accurate stress fields. On the other hand, the Almost-$C^1$ and D-Patch provide relatively easy construction on complex geometries. The Almost-$C^1$ method does not have limitations on the valence of boundary vertices, unlike the D-Patch, but is only applicable to biquadratic local bases. Following from these conclusions, future research directions are proposed, for example towards making the Approximate $C^1$ and Analysis-Suitable $G^1$ applicable to more complex geometries.
\end{abstract}

\begin{keyword}
Isogeometric Analysis \sep Unstructured Splines \sep Kirchhoff--Love shell \sep Biharmonic Equation
\end{keyword}

\end{frontmatter}



\section{Introduction}
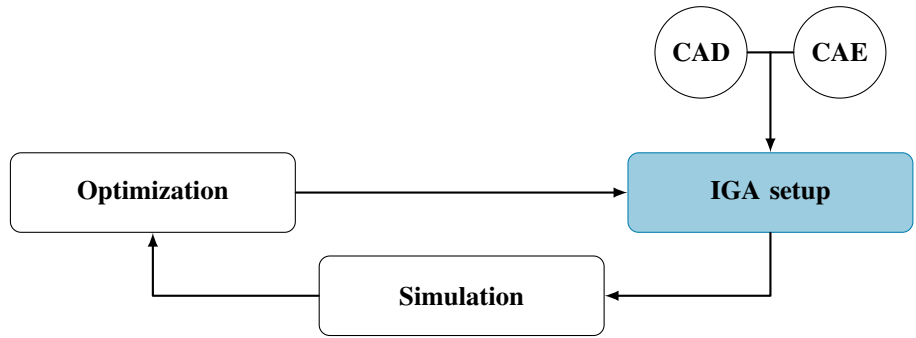
\begin{figure}
    \centering
    \tikzstyle{arrow} = [draw, -latex,thick,line cap=round,line join=round]
\tikzstyle{line} = [draw,thick,line cap=round,line join=round]
\tikzstyle{block} = [rectangle, draw, fill=white,
    text width=10em, text centered, rounded corners, minimum height=3em]
\tikzstyle{circleblock} = [circle, draw, fill=white,
    text width=2.5em, text centered, rounded corners, minimum height=1.5em]
\tikzstyle{diamondblock} = [diamond, draw, fill=white,rounded corners,
    text width=5em, text centered, minimum height=10em, minimum width=10em]
\tikzstyle{dummyblock} = [rectangle, fill=none,
    text width=10em, text centered, rounded corners, minimum height=3em]
\tikzstyle{noblock} = [rectangle, draw=none, fill=none,
    text width=10em, rounded corners, minimum height=0em,node distance = 0.05\linewidth, text = blue]
\tikzstyle{textblock} = [text width=7em,text centered]
\tikzstyle{textblock2} = [text width=14em,text centered]
\tikzstyle{textblock4} = [text width=28em,text centered]
\tikzstyle{octagon} = [shape=regular polygon, regular polygon sides=8, draw, fill=red!20,draw=red,
    text width=12em, text centered, rounded corners, minimum height=3em,minimum width=3em,inner sep=-2em]

\begin{tikzpicture}[node distance = 0.3cm]
    \node[block,fill=col1!40,draw=col1] (IGA)        {\textbf{IGA setup}};

    \node[above = of IGA,yshift=0.9cm] (CADCAE)        {};
	\node[circleblock,left = of CADCAE.center                     ] (CAD)        {\textbf{CAD}};
	\node[circleblock,right = of CADCAE.center                     ] (CAE)        {\textbf{CAE}};
	\path [line] (CAD) -- (CADCAE.center);
	\path [line] (CAE) -- (CADCAE.center);
	\path [arrow] (CADCAE.center) --  (IGA);

    \node[dummyblock,  below= of IGA    ] (empty1)       {};
    \node[block, left = of empty1                      ] (simulation)        {\textbf{Simulation}};

    \node[dummyblock,  left= of IGA    ] (empty2)       {};
    \node[block, left = of empty2                      ] (optimization)        {\textbf{Optimization}};

	\path [arrow] (IGA) |-  (simulation);
	\path [arrow] (simulation) -|  (optimization);
	\path [arrow] (optimization) --  (IGA);
\end{tikzpicture}
    \caption{General workflow for solving a physics problem and optimizing a geometry or topology coming from CAD and CAE processes. Starting from \emph{CAD} and \emph{CAE}, the \emph{IGA Setup} is performed. In this block, a computational basis is extracted from the geometry, to be used for simulation. Then, the \emph{Simulation} block involves assembly of the operators of the physics problem on the computational basis coming from the \emph{IGA Setup}. In case of shape or topology optimization problems, the simulation results are evaluated and the shape/topology is modified. From this changed shape/topology, a new computational basis can be obtained and the process can be repeated. The \emph{IGA Setup} block is marked to be elaborated further on in \cref{fig:IGA_Setup}.}
    \label{fig:workflow}
\end{figure}

Present day engineering disciplines depend on Computer-aided design (CAD) and numerical simulation models for physics for design and analysis. Typically, geometries designed in CAD are converted to meshes for an analysis with numerical techniques like Finite Element Methods (FEMs). Since the geometry description in CAD is based on splines whereas meshes for simulation are based on linear geometry approximations, geometric data is lost during this conversion. Isogeometric analysis \cite{Hughes2005} is the bridge between CAD and Computer-Aided Engineering (CAE), since it is employing splines as a basis for geometric design and numerical analysis. In practice, an isogeometric analysis and optimization workflow can be seen as depicted in \cref{fig:workflow}. Starting with a geometry from CAD as well as from material parameters, boundary conditions et cetera from CAE, isogeometric simulations and eventually geometry or topology optimization can be performed. The step connecting the inputs from CAD and CAE is referred to as \emph{IGA Setup} in \cref{fig:workflow}. This step takes care of the preparation for the simulation step, including the pre-processing of the geometry, if needed, and the construction of the isogeometric discretization space.\\

Due to the arbitrary smoothness of spline basis functions, isogeometric analysis has several advantages over conventional finite element methods. For example: (i) the introduction of $k$-refinements, which are proven to provide high accuracy per degree of freedom \cite{Bressan2019,Sande2020}; (ii) high accuracy in eigenvalue problems, e.g. for structural vibrations \cite{Cottrell2006,Hughes2008a,Hughes2014}; or (iii) geometric exactness in parametric design and interface problems, e.g. applied to the parametric design of prosthetic heart valves \cite{Xu2018}. Furthermore, the $C^1$-smooth discretization spaces allow to solve equations such as the biharmonic equation, the Cahn--Hilliard equations or the Kirchhoff--Love shell equations without introducing auxiliary variables. However, due to the tensor-product structure of the spline basis, higher-order smoothness can be enforced easily only on domains that allow simple patch partitions (e.g. an L-shape or an annulus), whereas on geometrically and topologically more complicated domains alternative approaches are required to solve equations that require basis functions of higher-order continuity.\\

\begin{figure}
    \centering
%
%
%
%

\tikzstyle{arrow} = [draw, -latex,thick,line cap=round,line join=round]
\tikzstyle{line} = [draw,thick,line cap=round,line join=round]
\tikzstyle{block} = [rectangle, draw, fill=white,
    text width=10em, text centered, rounded corners, minimum height=3em]
\tikzstyle{circleblock} = [circle, draw, fill=white,
    text width=2.5em, text centered, rounded corners, minimum height=1.5em]
\tikzstyle{diamondblock} = [diamond, draw, fill=white,rounded corners,
    text width=5em, text centered, minimum height=10em, minimum width=10em]
\tikzstyle{dummyblock} = [rectangle, fill=none,
    text width=10em, text centered, rounded corners, minimum height=3em]
\tikzstyle{noblock} = [rectangle, draw=none, fill=none,
    text width=10em, rounded corners, minimum height=0em,node distance = 0.05\linewidth, text = blue]
\tikzstyle{textblock} = [text width=7em,text centered]
\tikzstyle{textblock2} = [text width=14em,text centered]
\tikzstyle{textblock4} = [text width=28em,text centered]
\tikzstyle{octagon} = [shape=regular polygon, regular polygon sides=8, draw, fill=red!20,draw=red,
    text width=12em, text centered, rounded corners, minimum height=3em,minimum width=3em,inner sep=-2em]
\begin{tikzpicture}[node distance = 0.3cm,scale=0.5]
	\coordinate (CADCAE) at (0,0);
	\node[above] at (CADCAE)  {From CAD/CAE};
	\node[dummyblock,below=of CADCAE] (start){};
	\path[arrow] (CADCAE.center) -| (start.center);
	\node[dummyblock,below=of start] (mid) {};

	\node[block,left=of mid] (trimming) {Get trimming curves, etc.};
	\node[block,right=of mid] (mesh) {Construct/get quadrilateral mesh};
	\path[arrow] (start.center) -| (trimming);
	\path[arrow] (start.center) -| (mesh);

	\node[block,below=of trimming] (quadrule) {Setup problem using specialized quadrature, etc.};
	\path[arrow] (trimming) -- (quadrule);

	\node[block,below=of mesh] (weak) {Setup weak-coupling problem using penalty, etc.};
	\path[arrow] (mesh) -- (weak);
	\node[dummyblock,below=of mid] (bmid){};
	\node[octagon,below=of mid] (USconditions)
	{
		Unstructured Spline Constraints
		\begin{enumerate}
		\itemsep0em
		\item Manifold geometry
		\item Conforming mesh
		\end{enumerate}};
	\path[arrow] (mesh) -| (USconditions);
	\node[block,below=of USconditions,fill=col2!40,draw=col2] (USpre) {Unstructured spline pre-processing};
	\path[arrow] (USconditions) -- (USpre);
	\node[block,below=of USpre] (US) {Unstructured spline construction};
	\path[arrow] (USpre) -- (US);

	\node[dummyblock,below=of US] (bUS){};

	\path[line] (quadrule) |- (bUS.center);
	\path[line] (US) -- (bUS.center);
	\path[line] (weak) |- (bUS.center);

	\coordinate[below=of bUS] (bbUS);
	\path[arrow] (bUS.center) -- (bbUS.center);
	\node[below=of bbUS] {To Simulation};

	\path (trimming.north west)--(quadrule.south west) node[midway,above,rotate=90,text=col1]{Trimmed domain approaches};
	\path (mid.south west)--(US.south west) node[textblock,midway,above,rotate=90,text=col1]{Unstructured Splines};
	\path (mesh.north east)--(weak.south east) node[textblock,midway,below,rotate=90,text=col1!50]{Variational coupling};

\begin{pgfonlayer}{bg}    
	\draw[draw=col1!20,fill=col1!20][rounded corners]
([xshift=-2,yshift=2]trimming.north west) --
([xshift=2,yshift=2]trimming.north east) --
([xshift=2,yshift=-2]quadrule.south east) --
([xshift=-2,yshift=-2]quadrule.south west)
--cycle;
	\draw[draw=col1!20,fill=col1!20,opacity=0.5][rounded corners]
([xshift=-2,yshift=2]mid.north west) --
([xshift=2,yshift=2]mesh.north east) --
([xshift=2,yshift=-2]mesh.south east) --
([xshift=2,yshift=-2]mid.south east) --
([xshift=2,yshift=-2]US.south east) --
([xshift=-2,yshift=-2]US.south west)
--cycle;
	\draw[draw=col1!50,fill=col1!50,opacity=0.5][rounded corners]
	([xshift=-2,yshift=2]mesh.north west) --
([xshift=2,yshift=2]mesh.north east) --
([xshift=2,yshift=-2]weak.south east) --
([xshift=-2,yshift=-2]weak.south west);

--cycle;
\end{pgfonlayer}

\end{tikzpicture}
    \caption{Inside the \emph{IGA Setup} block from \cref{fig:workflow}, three methods are distinguished. Firstly, \emph{trimmed domain approaches} use trimming curves or surfaces to identify parts of a tensor-product domain as the actual domain. However, since elements can be trimmed poorly, \emph{specialized quadrature} rules and solver preconditioners are typically needed. Alternatives to trimming are \emph{weak coupling} or \emph{unstructured spline} methods. For both classes of methods, a geometry with a given topology needs to be decomposed into multiple sub-domains (i.e. patches) via \emph{quadrilateral meshing}. Given a quadrilateral mesh, \emph{weak methods} assemble extra penalty terms into the equation to be solved, or add extra equations to be solved to satisfy continuity constraints. Lastly, \emph{unstructured spline} constructions can be used to couple multiple domains by constructing a continuous basis. These methods, however can only be used on \emph{manifold geometries} and \emph{conforming meshes}. When these requirements are satisfied, \emph{unstructured spline pre-processing} is required before the \emph{unstructured spline construction} can take place. The \emph{pre-processing} is highlighted and will be elaborated on more in \cref{fig:USconstraints} in \cref{sec:qualitative}.}
    \label{fig:IGA_Setup}
\end{figure}
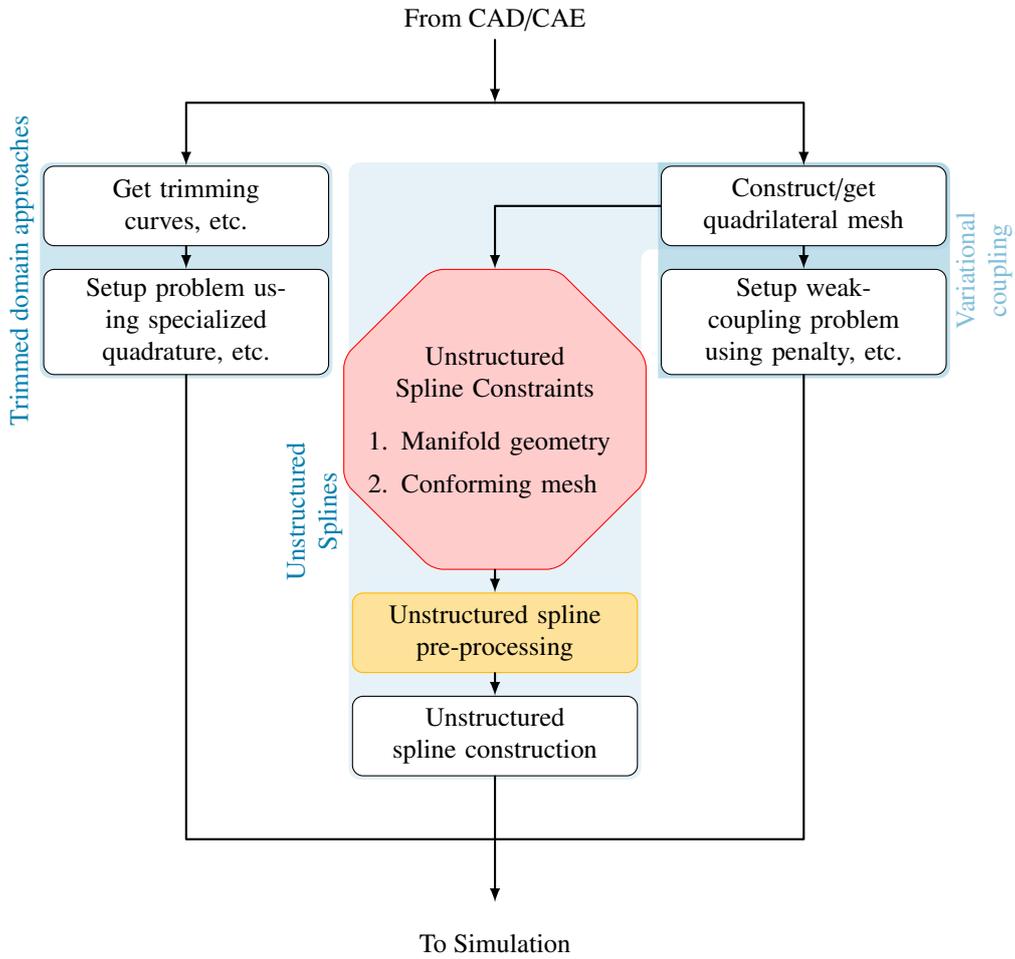

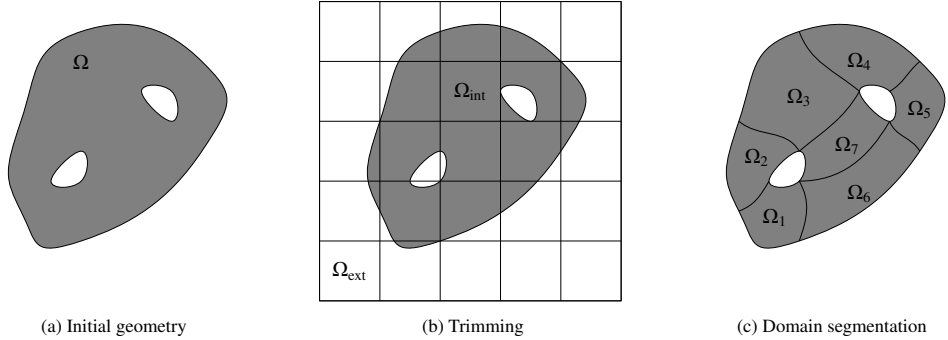
\begin{figure}
 \centering
 \begin{subfigure}{0.3\linewidth}
  \resizebox{\linewidth}{!}{\begin{tikzpicture}
\draw[white] (0,0) rectangle (5,5);
\filldraw [draw=black,fill=gray] plot [smooth cycle,tension=1] coordinates {(2,1) (1,1.5) (1,3) (2,4.5) (4,4) (4,2.5)};
\filldraw [draw=black,fill=white] plot [smooth cycle,tension=1] coordinates {(1.5,2) (2,2.5) (2,2)};
\filldraw [draw=black,fill=white] plot [smooth cycle,tension=1] coordinates {(3,3.5) (3.5,3) (3.5,3.5)};

\node[anchor=center] at (2.0,4.0) { $\Omega$};
\end{tikzpicture}}
  \caption{Initial geometry}
  \label{fig:initial_geometry}
 \end{subfigure}
\hfill
 \begin{subfigure}{0.3\linewidth}
  \resizebox{\linewidth}{!}{\begin{tikzpicture}
\draw (0,0) rectangle (5,5);
\filldraw [draw=black,fill=gray] plot [smooth cycle,tension=1] coordinates {(2,1) (1,1.5) (1,3) (2,4.5) (4,4) (4,2.5)};
\filldraw [draw=black,fill=white] plot [smooth cycle,tension=1] coordinates {(1.5,2) (2,2.5) (2,2)};
\filldraw [draw=black,fill=white] plot [smooth cycle,tension=1] coordinates {(3,3.5) (3.5,3) (3.5,3.5)};
\foreach \x in {0,1,...,5} 
{ 	
\draw (\x,0) -- (\x,5); 
\draw (0,\x) -- (5,\x); 
}

\node[anchor=center] at (0.5,0.5) { $\Omega_{\text{ext}}$};
\node[anchor=center] at (2.5,3.5) { $\Omega_{\text{int}}$};
\end{tikzpicture}}
  \caption{Trimming}
  \label{fig:trimming}
 \end{subfigure}
\hfill
 \begin{subfigure}{0.3\linewidth}
  \resizebox{\linewidth}{!}{\begin{tikzpicture}
\draw[white] (0,0) rectangle (5,5);
\filldraw [draw=black,fill=gray] plot [smooth cycle,tension=1] coordinates {(2,1) (1,1.5) (1,3) (2,4.5) (4,4) (4,2.5)};
\filldraw [draw=black,fill=white] plot [smooth cycle,tension=1] coordinates {(1.5,2) (2,2.5) (2,2)};
\filldraw [draw=black,fill=white] plot [smooth cycle,tension=1] coordinates {(3,3.5) (3.5,3) (3.5,3.5)};

\draw (2,1) to[in=-60,out=70] (2,2);
\draw (2,2) to[in=-120,out=00] (3.5,3);
\draw (3.5,3) to[in=120,out=-70] (4,2.5);
\draw (2,2.5) to[in=-120,out=40] (3,3.5);
\draw (2,2.5) to[in=-40,out=120] (1,3);
\draw (2,4.5) to[in=140,out=-70] (3,3.5);
\draw (1,1.5) to[in=-120,out=20] (1.5,2);
\draw (3.5,3.5) to[in=-140,out=30] (4,4);

\node[anchor=center] at (1.6,1.4) { $\Omega_1$};
\node[anchor=center] at (1.3,2.4) { $\Omega_2$};
\node[anchor=center] at (2.0,3.4) { $\Omega_3$};
\node[anchor=center] at (3.0,4.0) { $\Omega_4$};
\node[anchor=center] at (4.0,3.2) { $\Omega_5$};
\node[anchor=center] at (3.0,1.8) { $\Omega_6$};
\node[anchor=center] at (2.8,2.6) { $\Omega_7$};
\end{tikzpicture}}
  \caption{Domain segmentation}
  \label{fig:segmentation}
 \end{subfigure}
 \caption{Given an initial geometry $\Omega$ (a), trimming (b) uses the curves of the boundary of the original geometry to define the interior domain $\Omega_{\text{int}}$ and the exterior domain $\Omega_{\text{ext}}$. An alternative approach for modelling the domain is to use domain segmentation (c). Here, the domain is decomposed into several patches $\Omega_i$ which together define the full domain $\Omega$.}
 \label{fig:potatoes}
\end{figure}

For more complicated domains, the \emph{IGA setup} block in \cref{fig:workflow} involves a pre-processing step of either the geometry, the system of equations or the solution space to solve the original system of equations. In \cref{fig:IGA_Setup}, this pre-processing step is subdivided into three options: trimmed domain approaches, unstructured splines and variational coupling methods. Given an initial geometry (cf. \cref{fig:initial_geometry}), the trimmed domain approaches alter the tensor-product domain by defining parts of the domain that are physical or non-physical (cf. \cref{fig:trimming}). In case of unstructured splines or variational coupling methods, the geometry is decomposed into multiple different patches (cf. \cref{fig:segmentation}) on which continuity conditions are enforced by constructing a smooth basis (unstructured splines) or by adding extra terms to the system of equations (variational coupling approaches). In \cref{sec:literature} of this paper, a review of trimmed domain approaches, unstructured splines and variational coupling methods is provided. Examples include immersed methods, degenerate patches and Nitsche's method, respectively. In case of simple geometries (and given the right inputs) the methods are identical. \\

As shown in \cref{fig:IGA_Setup}, each class of methods has its own characteristics and previous work has provided several comparisons of methods among each other, which are elaborated more in \cref{sec:literature}. In the context of the workflow sketched in \cref{fig:workflow}, unstructured splines provide a valuable alternative to the other methods, since they are constructed for a fixed topology and hence the computational costs of their construction are not related to changing shapes or moving domains. However, recent developments mainly focused on different unstructured spline methods separately, rather than providing a valuable comparison. In this paper, we therefore provide a qualitative and a quantitative comparison of a selection of unstructured spline constructions. We consider finite, piece-wise polynomial spline constructions, hence we do not include rational constructions or infinite representations, such as subdivision surfaces. More precisely, we compare examples of (globally) $G^1$-smooth multi-patch constructions (the Analysis-Suitable $G^1$ construction of \cite{Farahat2023} and the Approximate $C^1$ construction of \cite{Weinmuller2021}), the D-Patch method of \cite{Toshniwal2017} and the Almost-$C^1$ construction of \cite{Takacs2023}, motivated in \cref{subsec:US}. The selected methods are qualitatively compared based on their properties, and quantitatively based on several different examples with biharmonic and Kirchhoff--Love shell equations. The aim of this paper is to provide a fair comparison\footnote{We believe that a comparison like the one presented in this paper is never fully unbiased, since the authors have contributed to different methods in previous publications and do not represent the entire research community.} of these methods, providing a good overview of the strengths and weaknesses of each method in different cases.\\

The paper is outlined as follows: in \cref{sec:literature}, a detailed overview of the methods appearing in \cref{fig:IGA_Setup} is provided. In \cref{sec:qualitative} we provide a qualitative analysis of the four constructions that are discussed in this paper, while in \cref{sec:quantitative} we provide a quantitative analysis of all methods. There, we present five benchmark problems solving either a biharmonic or a Kirchhoff--Love equation. These benchmark problems serve different purposes and we compare which method, in which setting, performs best. In \cref{sec:conclusions}, we conclude this paper based on the findings from the previous sections and we provide directions for future research.\\

\section{Multi-patch isogeometric analysis: literature review}\label{sec:literature}
As discussed in the introduction of this paper, in particular in \cref{fig:IGA_Setup}, three classes of methods for the modelling of complicated domains can be characterised: trimmed domain approaches, variational coupling methods and unstructured splines. The goal of all methods is to achieve a certain level of continuity across the whole analysis domain such that multi-patch isogeometric analysis can be performed for example for the Kirchhoff--Love shell model \cite{Kiendl2009}, the biharmonic equation or the Cahn--Hilliard equation \cite{Gomez2008}. As shown in \cref{fig:potatoes}, trimmed domain approaches use the fact that parts of tensor-product geometries are trimmed away, using trimming curves to separate regions of interest and regions that should be omitted, see \cref{fig:trimming}. Variational coupling approaches and unstructured splines are defined on multi-patch domains, typically following from a segmentation of the original domain, see \cref{fig:segmentation}. In case of variational coupling methods, the system of equations is enriched with terms that will enforce continuity (typically in a weak sense) between the patches. In case of unstructured spline constructions, a basis is constructed on the multi-patch object, where certain smoothness is enforced strongly. When starting from a trimmed geometry, the step of creating a multi-patch domain decomposition (i.e. \emph{untrimming}) from an arbitrary geometry with an arbitrary topology is a very important step in the application of weak coupling methods and unstructured spline constructions, as can be seen in the flowchart in \cref{fig:IGA_Setup}. In this paper, however, the topic of untrimming will not be discussed as it is out of scope of our study. Hence, the reader is referred to \cite{Massarwi2018b,Hiemstra2020} for an overview of these methods.\\

In this section, an overview of the trimmed domain approaches (\cref{subsec:trimming}), variational coupling methods (\cref{subsec:variational}) and unstructured splines (\cref{subsec:US}) is provided. A fourth method, which will not be discussed in this section, is to introduce auxiliary variables for derivatives of the solution, so that $C^1$ continuity requirements are reduced to $C^0$ and standard interface coupling can be used. These so-called mixed formulations are common in conventional FEM, although recent advances have also been made for Kirchhoff--Love plates and shells and the biharmonic eigenvalue problem \cite{rafetseder2018decomposition,Rafetseder2019,kosin2023new}.

\subsection{Trimming approaches}\label{subsec:trimming}
Trimming is a technique where so-called trimming curves or surfaces separate parts of tensor-product spline domains to define a geometry. Trimming is a common technique to represent complex geometries in CAD, and typically geometries consist of multiple trimmed patches with boundary and interface curves trimming the actual patches. We refer to the work  \cite{Marussig2018} for an overview of trimming methods in isogeometric analysis. Generalizing the idea of trimming to techniques where curves or surfaces are used to define the domain of interest as \emph{trimmed domain approaches}, several approaches have been proposed to perform simulations on complex geometries, including the finite cell method \cite{Parvizian2007,Duster2008,Schillinger2015a}, Cut-FEM \cite{Burman2015} or immersed methods \cite{Hsu2016,Kamensky2015}. The advantage of these methods is that the trimmed CAD geometries could directly be used for analysis. However, when only small parts of the physical domain are cut, leading to small cut elements, numerical difficulties can occur in the conditioning of the system, leading to solver instabilities or accuracy problems \cite{dePrenter2022}. Therefore, the analysis of complex trimmed geometries via methods like the FCM typically require special quadrature schemes to take into account small cut cells \cite{Divi2020a} or preconditioners to stabilize the numerical analysis \cite{dePrenter2020}. In the context of Kirchhoff--Love shell modelling, isogeometric analysis on trimmed geometries has been performed in several studies \cite{Guo2015a,Guo2018,Coradello2020} including some with focus on multi-patch coupling \cite{Coradello2020,Coradello2020d,Coradello2021,Guo2017,Leidinger2019}.\\

\subsection{Variational coupling methods}\label{subsec:variational}
We define variational coupling methods as methods that modify the system of equations to enforce certain continuity across patch interfaces. Examples of these methods are penalty methods, Nitsche's methods, mortar methods or Lagrangian penalized methods. In the context of Kirchhoff--Love shell analysis, these weak coupling methods have received a lot of attention in previous studies and an overview is provided by \cite{Paul2020}. Firstly, an in-plane coupling was proposed in \cite{Kiendl2011} together with a method for coupling non-manifold patches using the so-called bending strip method \cite{Kiendl2010}. Later, weak coupling approaches have been developed for multi-patch domains. Here, coupling terms can be added inside the existing variational formulation (referred to as Nitsche's or penalty methods) or imposed by Lagrange multipliers (referred to as mortar methods).\\

Several works on Nitsche techniques (cf. \cite{Nitsche1971}) for isogeometric analysis have been published starting from the imposition of boundary conditions \cite{Ruess2013}, towards multi-patch coupling and the coupling of patches \cite{Ruess2014}, later using a non-symmetric parameter-free Nitsche's method \cite{Schillinger2016}. Nitsche's methods have been applied to Kirchhoff plates \cite{Hu2018}, Kirchhoff--Love shells \cite{Guo2015,Guo2015a,Benzaken2021}, hyperelastic 2D elasticity \cite{Du2020} and the biharmonic equation \cite{Weinmuller2022,Weinmuller2021} and for modelling local subdomains \cite{Bouclier2018} for elasticity simulations. The advantages of Nitsche's methods are that the formulation is variationally consistent and requires only mild stabilization, which can be performed automatically, by estimating the stability parameter. However, the involved integral terms are complicated expressions that impose high implementation and assembly efforts. Therefore, coupling approaches using only penalization have been developed \cite{Duong2017,Breitenberger2015a,Herrema2019,Leonetti2020,Pasch2021a}. Although several improvements have been made in these works, the main disadvantage of penalty methods is that a suitable penalty parameter has to be chosen. Using the super penalty approach \cite{Coradello2021a,Coradello2021}, the computation of the penalty parameter can be automated. However, this method has not yet been tested for non-linear shell problems or on `dirty' geometries. Both Nitsche's and penalty methods can be used to couple geometries that are non-manifold, i.e. geometries that have out-of-plane connections like stiffened structures by penalizing changes in the angle of patches on an interface. Furthermore, the methods can handle interfaces with non-matching parameterizations.\\

Instead of adding coupling terms in the variational form, as is done in Nitsche’s and penalty methods, mortar methods \cite{Bernardi1993} add extra degrees of freedom by introducing Lagrange multipliers which are required to resolve additional coupling conditions. The use of mortar methods to couple non-conforming isogeometric sub-domains was first done by \cite{Hesch2012}. In \cite{Dornisch2015} the FEA-based approach of \cite{Bernardi1993} was extended for NURBS-based IGA, but the aim was to develop a method for $C^0$-coupling for Reissner-Mindlin shells, hence insufficient for isogeometric Kirchhoff--Love shells. A mortar method aiming to establish $C^1$ coupling is given in \cite{Bouclier2017} and a method that provides $C^n$ continuity was given by \cite{Dittmann2019,Dittmann2020}. Furthermore, $G^1$ mortar coupling, referred to as extended mortar coupling, was presented in \cite{Schuss2019} for Kirchhoff--Love shells, based on a coupling in least square sense. On the other hand, in \cite{Benvenuti2023} a mortar method to enforce $C^1$ coupling for the biharmonic equation was developed, where the Lagrange multiplier spaces are constructed similarly to \cite{Brivadis2015a} for $C^0$-coupling. An approach to reduce computational costs involved in finding Lagrange multipliers is called dual mortaring \cite{Wohlmuth2000}, where Lagrange multipliers are eliminated using a compact dual basis. This approach has been developed for Bezier elements \cite{Zou2018} and it has been applied for Kirchhoff--Love shells \cite{Miao2021} and a bi-orthogonal spline space has been presented for weak dual mortaring for patch coupling \cite{Wunderlich2019}. In \cite{Horger2019}, a hybrid method was provided and applied to Kirchhoff plates, which combines mortar methods and penalty methods. Lastly, a comparison of Nitsche, penalty and mortar methods is given by \cite{Apostolatos2019b}. For a more complete overview of mortar methods for isogeometric analysis, the reader is referred to \cite{Hesch2022}. In general, mortar methods have the advantage over Nitsche's methods that there are no parameters involved and that the implementation efforts are lower. However, the disadvantage is that a suitable spline space needs to be found for the Lagrange multipliers \cite{Brivadis2015a,Dornisch2021,Benvenuti2023}. Like Nitsche's and penalty methods, mortar methods can handle non-matching parameterizations and non-manifold interfaces, the latter by similar penalization of interfacing patches.

\subsection{Unstructured splines}\label{subsec:US}
Compared to weak coupling methods, unstructured spline constructions do not alter the system of equations to be solved. Instead, the computational basis is modified such that it satisfies continuity conditions across patch interfaces. Unstructured splines are typically constructed for in-plane (i.e. manifold) interfaces and not on out-of-plane (i.e. non-manifold) interfaces, since the notion of smoothness is uniquely defined only in the former setting. However, unstructured spline constructions for non-manifold interfaces are possible, e.g. as in \cite{Cirak2011,Ying2001,Moulaeifard2023} in the context of subdivision. Furthermore, unstructured spline constructions are typically constructed on conforming interfaces, i.e. interfaces with matching meshes, but, as long as the patch parameterizations are matching, this can be overcome by taking the knot vector union of the interface patches. However, the advantage of unstructured spline constructions is that as soon as the basis is constructed for a certain untrimmed geometry, there are no additional costs involved other than evaluation costs for changing shapes, which make unstructured spline bases suitable for shape optimization problems. In case of topology changes or large changes of the shape, however, the mesh topology of the unstructured spline space has to be changed as well. Unlike weak methods, which are typically based on the introduction of penalties (e.g. in terms of energy), unstructured spline constructions are typically provided as generic geometric methods that are applicable to any equation that requires $C^1$ coupling across multi-patch interfaces. With the advance of isogeometric analysis, the interest in parametrically $C^1$ and geometrically $G^1$ splines has grown. An overview of smooth multi-patch discretizations for isogeometric analysis can be found in \cite{Hughes2021}, and a small overview is provided below. We distinguish between enforcing parametric continuity, i.e., the type of continuity between mesh elements within a regular tensor-product spline patch, and general geometric continuity, cf. \cite{Groisser2015}. In the following, three types of constructions are classified, depending on their continuity on patch interfaces, around vertices and in the patch interior:
\begin{itemize}
    \item Patch coupling with geometric continuity on patch interfaces and parametric continuity inside patches.
    \item Patch coupling with parametric continuity everywhere.
    \item Patch coupling with parametric continuity almost everywhere.
\end{itemize}

Although other constructions outside of these categories exist, e.g., \cite{Majeed2017,Yuan2014}, our review is restricted to the aforementioned categories since the methods considered in \cref{sec:qualitative,sec:quantitative} fall into these categories.

\subsubsection*{Geometric continuity on patch interfaces and parametric continuity inside patches}
This first category of unstructured spline constructions assumes that a fixed $C^0$-matching multi-patch parameterization is given. On this multi-patch domain, a $C^1$-smooth isogeometric space is constructed. As shown in \cite{Groisser2015}, for any isogeometric function the $C^1$ condition over each interface is equivalent to a $G^1$ geometric continuity condition of the graph surface corresponding to the function. If the domain is planar and the patches are bilinear, then the $C^1$ constraints can be resolved and a $C^1$ spline space was constructed by \cite{Kapl2015} and applied to the isogeometric analysis of the biharmonic equation in \cite{Kapl2017}. It could be shown in \cite{Kapl2021a} and \cite{Groselj2020} that $C^1$ splines over bilinear quadrilaterals and mixed (bi)linear quadrilateral/triangle meshes possess optimal approximation properties. Furthermore, the work \cite{Kapl2021} studied the arbitrary $C^n$-smooth spline space for bi-linear multi-patch parameterizations, based on their previously published findings.\\

Considering general $C^0$-matching multi-patch domains, the work of \cite{Collin2016} introduces the class of analysis-suitable $G^1$ (AS-$G^1$) multi-patch parameterizations which includes bi-linear patches. This AS-$G^1$ condition is in general required to obtain optimal approximation properties. The condition implies that the gluing data for $G^1$ continuity is linear, which is explained in more detail in \cref{subsubsec:ASG1}.  While it could be shown in \cite{Kapl2018} that all planar multi-patch domains possess AS-$G^1$ reparameterizations, creating AS-$G^1$ surface domains is more difficult. Several strategies to achieve this were introduced in \cite{Farahat2023}, thus making $C^1$-smooth multi-patch parameterizations applicable to biharmonic equations and isogeometric Kirchhoff--Love shell models \cite{Farahat2023a}. In the work of \cite{Reichle2023} the construction of \cite{Collin2016} is used to develop a scaled-boundary model for smooth Kirchhoff--Love shells, similar to the approach of \cite{Arf2023} for Kirchhoff plates.\\

Alternatively to constructing an AS-$G^1$ parameterization, one can relax the smoothness conditions. This was done in \cite{Weinmuller2021}, where the construction of an approximate $C^1$ (Approx.~$C^1$) space is presented. The basis construction is explicit, possesses the same degree-of-freedom structure as an AS-$G^1$ space, but the $C^1$ condition is not satisfied exactly but only approximately. It defaults to the AS-$G^1$ construction when the AS-$G^1$ requirements are met. In \cite{Weinmuller2022} a comparison of the presented space with Nitsche's method was performed, yielding optimal convergence results without the need of a coupling terms. More details on the Approx.~$C^1$ method are provided in \cref{subsubsec:ApproxC1}.\\

\subsubsection*{Parametric continuity everywhere}
The starting point for this class of constructions is different from the previous. Here we create smooth splines in a parametric sense between neighboring mesh elements. Such parametric $C^1$ conditions are easy to resolve, but they lead to singularities at vertices of valencies other than four, so-called extraordinary vertices. This is due to the conflicting coupling conditions on partial derivatives around the EVs, which lead to all partial derivatives to vanish there. Inspired by the Degenerate Patch (D-Patch) approach from \cite{Reif1997}, the works of \cite{Nguyen2016,Toshniwal2017} provide $C^1$ smooth spline spaces for multi-patch geometries with parametric smoothness everywhere. On extraordinary vertices (EVs), which is a junction between 3 or 5 or more patches (i.e. valence $\nu>2,\nu\neq4$), the original D-Patch method shows a singularity of the basis in EVs combined with a reduction of degrees of freedom in this point. An improvement of the D-Patch method was presented in \cite{Nguyen2016}, by splitting elements around the EVs such that every element is associated to four degrees of freedom. However, this construction does not have non-negativity and is based on PHT splines, which have limited smoothness. A new design and analysis framework for multi-patch geometries was presented in \cite{Toshniwal2017}, based on D-Patches with T-splines for refinement and non-negative splines yielding optimal convergence properties. This was also demonstrated in \cite{Casquero2020} for isogeometric Kirchhoff--Love shells. In the work \cite{Wei2022}, it is motivated that this construction can also be used if only one element around the EV is isolated. More details on the D-Patch method are provided in \cref{subsubsec:dpatch}.\\

Alternatively, subdivision surface based constructions lead to unstructured splines that are parametrically continuous everywhere, cf. \cite{barendrecht2013isogeometric,riffnaller2016isogeometric,pan2016isogeometric,barendrecht2018efficient,zhang2018subdivision}. However, such approaches require an infinite number of polynomial pieces around each EV. Thus, we discard them for our comparison. Moreover, in general their approximation properties are severely reduced near EVs \cite{Takacs2023b}.\\

\subsubsection*{Parametric continuity almost everywhere}
As mentioned previously, imposing parametric continuity everywhere leads to singularities at all EVs. Thus, instead of constructing a space with full parametric continuity, spaces with parametric continuity almost everywhere except around the EVs can also be considered. This way, one ends up with regular, smooth rings around EVs which then need to be filled in some way. Such so-called \emph{hole-filling} techniques are commonplace in geometric modeling and can also be used to construct smooth spaces for isogeometric analysis, cf. \cite{Nguyen2016a,Karciauskas2015,Karciauskas2016,Karciauskas2021,Karciauskas2021a,Karciauskas2021b}. We focus here on the simplest possible way of resolving this issue, which is to enforce only $C^0$-smoothness near the EVs and $G^1$ at the EV, namely the Almost-$C^1$ construction proposed in \cite{Takacs2023}. Similar constructions, which enforce no additional smoothness near EVs were proposed for mixed quadrilateral/triangle meshes in \cite{Toshniwal2022} and for arbitrary degree multi-patch B-splines with enhanced smoothness (MPBES) in \cite{Buchegger2016}.\\

The Almost-$C^1$ construction we consider here yields piece-wise biquadratic splines which are $C^1$ in regular regions and which have reduced smoothness around extraordinary vertices,  independent of the valence or the location (i.e. interior or boundary EVs). In contrast to that, most commonly used hole-filling approaches yield exactly $C^1$-smooth spaces but introduce locally polynomials of higher degree, or require a higher degree to start with, such as the construction presented in \cite{Marsala2022a}, which converts Catmull--Clark subdivision surfaces to $G^1$-smooth piece-wise biquintic elements. While exact smoothness is of relevance for geometric modeling, it is not necessary from an analysis point of view.\\

%



\section{Qualitative comparison}\label{sec:qualitative}
In the qualitative comparison of this paper, we focus on the properties of different unstructured spline constructions and their implication on the application of these constructions in a workflow as in \cref{fig:workflow}. More precisely, we comment on the continuity of each construction and their nestedness properties and we aim to provide a set of requirements for the \emph{unstructured spline pre-processing} block in \cref{fig:IGA_Setup}. Since the qualitative comparison of the considered methods in this paper mostly covers properties of the methods and their implications, mathematical details about the construction or convergence properties are not provided. For more details, the reader is referred to \cite{Farahat2023} for the Analysis-Suitable $G^1$ (AS-$G^1$) method, which extends the 2D construction from \cite{Kapl2019}, to \cite{Weinmuller2021} for the Approximate $C^1$ (Approx.~$C^1$) method, to \cite{Toshniwal2017} for the Degenerate Patches (D-Patch) and to \cite{Takacs2023} for the Almost-$C^1$ method. However, for the qualitative comparison, some key terms are introduced as preliminaries.\\

\begin{figure}
    \centering
    \begin{subfigure}[t]{0.3\linewidth}
    \centering
    \resizebox{\linewidth}{!}{\begin{tikzpicture}
[
        ,boundary/.style={thick, draw=black},
        ,interior/.style={thin, draw=gray},
		,green_thin/.style={thin, draw=gray},
		,red_thin/.style={thin, draw=gray},
		,blue_thin/.style={thin, draw=gray},
		,orange_thin/.style={thin, draw=gray},
		,yellow_thin/.style={thin, draw=gray},
		,pink_thin/.style={thin, draw=gray},
		,purple_thin/.style={thin, draw=gray},
		,brown_thin/.style={thin, draw=gray},
		,lime_thin/.style={thin, draw=gray},
		,cyan_thin/.style={thin, draw=gray},
		,magenta_thin/.style={thin, draw=gray},
		,EV/.style={fill=gray},
		,bEV/.style={fill=black},
]
		\node [] (0) at (0, 3) {};
		\node [] (1) at (0, 3.5) {};
		\node [] (2) at (0, 3.95) {};
		\node [] (3) at (0, 2) {};
		\node [] (4) at (0, 1.5) {};
		\node [] (5) at (0, 1) {};
		\node [] (6) at (0, 0.5) {};
		\node [] (7) at (0, 4.5) {};
		\node [] (8) at (0.5, 2.5) {};
		\node [] (9) at (1, 2.5) {};
		\node [] (10) at (1.5, 2.5) {};
		\node [] (11) at (2, 2.5) {};
		\node [] (12) at (-0.5, 2.5) {};
		\node [] (13) at (-1, 2.5) {};
		\node [] (14) at (-1.45, 2.5) {};
		\node [] (15) at (-2, 2.5) {};
		\node [] (16) at (-0.675, 3.175) {};
		\node [] (17) at (-1, 3.5) {};
		\node [] (18) at (-1.425, 3.925) {};
		\node [] (19) at (0.675, 3.175) {};
		\node [] (20) at (1, 3.5) {};
		\node [] (21) at (1.425, 3.925) {};
		\node [] (22) at (-1.425, 1.075) {};
		\node [] (23) at (-1, 1.5) {};
		\node [] (24) at (-0.675, 1.8) {};
		\node [] (25) at (0.35, 2.15) {};
		\node [] (26) at (1, 1.5) {};
		\node [] (27) at (0.5, 1.5) {};
		\node [] (28) at (1, 2) {};
		\node [] (29) at (1.5, 2) {};
		\node [] (30) at (1.5, 1.15) {};
		\node [] (31) at (0.5, 1) {};
		\node [] (32) at (1, 1) {};
		\node [] (36) at (2, 2) {};
		\node [] (37) at (-0.325, 2.825) {};
		\node [] (38) at (-0.325, 2.175) {};
		\node [] (39) at (0.325, 2.825) {};
		\node [] (40) at (0, 0) {};
		\node [] (41) at (0.5, 0) {};
		\node [] (42) at (1, 0) {};
		\node [] (47) at (2, 1) {};
		\node [] (54) at (1.5, 0.5) {};
		\node [] (55) at (0.825, 0.5) {};
		\node [] (56) at (0.5, 0.5) {};
		\node [] (57) at (2, 1.5) {};
		\node [] (58) at (1.5, 1.5) {};

		\draw [style=boundary] (11.center) to (21.center);
		\draw [style=boundary] (21.center) to (7.center);
		\draw [style=boundary] (7.center) to (18.center);
		\draw [style=boundary] (18.center) to (15.center);
		\draw [style=boundary] (15.center) to (22.center);
		\draw [style=boundary] (22.center) to (6.center);
		\draw [style=boundary] (8.center) to (25.center);
		\draw [style=boundary] (25.center) to (3.center);
		\draw [style=interior] (19.center) to (20.center);
		\draw [style=boundary] (3.center) to (38.center);
		\draw [style=boundary] (38.center) to (12.center);
		\draw [style=boundary] (12.center) to (37.center);
		\draw [style=boundary] (37.center) to (0.center);
		\draw [style=boundary] (0.center) to (39.center);
		\draw [style=boundary] (39.center) to (8.center);
		\draw [style=interior] (24.center) to (38.center);
		\draw [style=interior] (12.center) to (13.center);
		\draw [style=interior] (37.center) to (16.center);
		\draw [style=interior] (0.center) to (1.center);
		\draw [style=interior] (39.center) to (19.center);
		\draw [style=interior] (16.center) to (17.center);
		\draw [style=interior] (20.center) to (21.center);
		\draw [style=interior] (1.center) to (2.center);
		\draw [style=interior] (2.center) to (7.center);
		\draw [style=interior] (17.center) to (18.center);
		\draw [style=interior] (14.center) to (15.center);
		\draw [style=interior] (23.center) to (22.center);
		\draw [style=interior] (23.center) to (24.center);
		\draw [style=interior] (13.center) to (14.center);
		\draw [style=boundary] (40.center) to (41.center);
		\draw [style=boundary] (41.center) to (42.center);
		\draw [style=boundary] (6.center) to (40.center);
		\draw [style={green_thin}] (32.center) to (26.center);
		\draw [style={green_thin}] (26.center) to (28.center);
		\draw [style={green_thin}] (28.center) to (9.center);
		\draw [style={green_thin}] (9.center) to (19.center);
		\draw [style={green_thin}] (19.center) to (1.center);
		\draw [style={green_thin}] (1.center) to (16.center);
		\draw [style={green_thin}] (16.center) to (13.center);
		\draw [style={green_thin}] (13.center) to (24.center);
		\draw [style={green_thin}] (24.center) to (4.center);
		\draw [style={green_thin}] (4.center) to (27.center);
		\draw [style={green_thin}] (27.center) to (26.center);
		\draw [style={red_thin}] (32.center) to (31.center);
		\draw [style={red_thin}] (31.center) to (5.center);
		\draw [style={red_thin}] (5.center) to (23.center);
		\draw [style={red_thin}] (23.center) to (14.center);
		\draw [style={red_thin}] (14.center) to (17.center);
		\draw [style={red_thin}] (17.center) to (2.center);
		\draw [style={red_thin}] (2.center) to (20.center);
		\draw [style={red_thin}] (20.center) to (10.center);
		\draw [style={red_thin}] (10.center) to (29.center);
		\draw [style={blue_thin}] (6.center) to (5.center);
		\draw [style={blue_thin}] (5.center) to (4.center);
		\draw [style={blue_thin}] (4.center) to (3.center);
		\draw [style={pink_thin}] (25.center) to (28.center);
		\draw [style={pink_thin}] (28.center) to (29.center);
		\draw [style={pink_thin}] (29.center) to (36.center);
		\draw [style={purple_thin}] (25.center) to (27.center);
		\draw [style={purple_thin}] (27.center) to (31.center);
		\draw [style={purple_thin}] (56.center) to (41.center);
		\draw [style={purple_thin}] (31.center) to (56.center);
		\draw [style={interior}] (11.center) to (10.center);
		\draw [style={interior}] (10.center) to (9.center);
		\draw [style={interior}] (9.center) to (8.center);
		\draw [style={boundary}](42.center) to (54.center);
		\draw [style={boundary}](54.center) to (47.center);
		\draw [style={orange_thin}] (32.center) to (54.center);
		\draw [style={lime_thin}] (30.center) to (47.center);
		\draw [style={brown_thin}] (55.center) to (42.center);
		\draw [style={cyan_thin}] (32.center) to (55.center);
		\draw [style={magenta_thin}] (32.center) to (30.center);
		\draw [style={yellow_thin}] (56.center) to (6.center);
		\draw [style={yellow_thin}] (55.center) to (56.center);
		\draw [style={green_thin}] (26.center) to (58.center);
		\draw [style={red_thin}] (58.center) to (30.center);
		\draw [style={red_thin}] (29.center) to (58.center);
		\draw [style={green_thin}] (58.center) to (57.center);
		\draw [style={boundary}](11.center) to (36.center);
		\draw [style={boundary}](36.center) to (57.center);
		\draw [style={boundary}](57.center) to (47.center);

		\fill[bEV] (6) circle [radius=0.05];
		\fill[EV] (55) circle [radius=0.05];
		\fill[EV] (30) circle [radius=0.05];
		\fill[EV] (32) circle [radius=0.05];
		\fill[bEV] (25) circle [radius=0.05];
\end{tikzpicture}}
    \caption{A simple mesh with boundary edges in black and interior edges in gray. The boundary extraordinary vertices (bEVs), i.e. the vertices on a boundary with valence $\nu\geq3$ are denoted by a black circle and the interior extraordinary vertices (iEVs), i.e. interior vertices with valence $\nu\geq4, \nu\neq4$ are denoted by gray circles.}
    \label{fig:Illustration_patches_from_mesh1}
    \end{subfigure}
    \hfill
    \begin{subfigure}[t]{0.3\linewidth}
    \centering
    \resizebox{\linewidth}{!}{\begin{tikzpicture}
[
        ,boundary/.style={thick, draw=black},
        ,interior/.style={thin, draw=gray},
		,green_thin/.style={thin, draw=green,-latex},
		,red_thin/.style={thin, draw=red,-latex},
		,blue_thin/.style={thin, draw=gray},
		,orange_thin/.style={thin, draw=orange,-latex},
		,yellow_thin/.style={thin, draw=gray},
		,pink_thin/.style={thin, draw=gray},
		,purple_thin/.style={thin, draw=gray},
		,brown_thin/.style={thin, draw=gray},
		,lime_thin/.style={thin, draw=gray},
		,cyan_thin/.style={thin, draw=cyan,-latex},
		,magenta_thin/.style={thin, draw=magenta,-latex},
		,EV/.style={fill=gray},
		,bEV/.style={fill=black},
]
		\node [] (0) at (0, 3) {};
		\node [] (1) at (0, 3.5) {};
		\node [] (2) at (0, 3.95) {};
		\node [] (3) at (0, 2) {};
		\node [] (4) at (0, 1.5) {};
		\node [] (5) at (0, 1) {};
		\node [] (6) at (0, 0.5) {};
		\node [] (7) at (0, 4.5) {};
		\node [] (8) at (0.5, 2.5) {};
		\node [] (9) at (1, 2.5) {};
		\node [] (10) at (1.5, 2.5) {};
		\node [] (11) at (2, 2.5) {};
		\node [] (12) at (-0.5, 2.5) {};
		\node [] (13) at (-1, 2.5) {};
		\node [] (14) at (-1.45, 2.5) {};
		\node [] (15) at (-2, 2.5) {};
		\node [] (16) at (-0.675, 3.175) {};
		\node [] (17) at (-1, 3.5) {};
		\node [] (18) at (-1.425, 3.925) {};
		\node [] (19) at (0.675, 3.175) {};
		\node [] (20) at (1, 3.5) {};
		\node [] (21) at (1.425, 3.925) {};
		\node [] (22) at (-1.425, 1.075) {};
		\node [] (23) at (-1, 1.5) {};
		\node [] (24) at (-0.675, 1.8) {};
		\node [] (25) at (0.35, 2.15) {};
		\node [] (26) at (1, 1.5) {};
		\node [] (27) at (0.5, 1.5) {};
		\node [] (28) at (1, 2) {};
		\node [] (29) at (1.5, 2) {};
		\node [] (30) at (1.5, 1.15) {};
		\node [] (31) at (0.5, 1) {};
		\node [] (32) at (1, 1) {};
		\node [] (36) at (2, 2) {};
		\node [] (37) at (-0.325, 2.825) {};
		\node [] (38) at (-0.325, 2.175) {};
		\node [] (39) at (0.325, 2.825) {};
		\node [] (40) at (0, 0) {};
		\node [] (41) at (0.5, 0) {};
		\node [] (42) at (1, 0) {};
		\node [] (47) at (2, 1) {};
		\node [] (54) at (1.5, 0.5) {};
		\node [] (55) at (0.825, 0.5) {};
		\node [] (56) at (0.5, 0.5) {};
		\node [] (57) at (2, 1.5) {};
		\node [] (58) at (1.5, 1.5) {};

		\draw [style=boundary] (11.center) to (21.center);
		\draw [style=boundary] (21.center) to (7.center);
		\draw [style=boundary] (7.center) to (18.center);
		\draw [style=boundary] (18.center) to (15.center);
		\draw [style=boundary] (15.center) to (22.center);
		\draw [style=boundary] (22.center) to (6.center);
		\draw [style=boundary] (8.center) to (25.center);
		\draw [style=boundary] (25.center) to (3.center);
		\draw [style=interior] (19.center) to (20.center);
		\draw [style=boundary] (3.center) to (38.center);
		\draw [style=boundary] (38.center) to (12.center);
		\draw [style=boundary] (12.center) to (37.center);
		\draw [style=boundary] (37.center) to (0.center);
		\draw [style=boundary] (0.center) to (39.center);
		\draw [style=boundary] (39.center) to (8.center);
		\draw [style=interior] (24.center) to (38.center);
		\draw [style=interior] (12.center) to (13.center);
		\draw [style=interior] (37.center) to (16.center);
		\draw [style=interior] (0.center) to (1.center);
		\draw [style=interior] (39.center) to (19.center);
		\draw [style=interior] (16.center) to (17.center);
		\draw [style=interior] (20.center) to (21.center);
		\draw [style=interior] (1.center) to (2.center);
		\draw [style=interior] (2.center) to (7.center);
		\draw [style=interior] (17.center) to (18.center);
		\draw [style=interior] (14.center) to (15.center);
		\draw [style=interior] (23.center) to (22.center);
		\draw [style=interior] (23.center) to (24.center);
		\draw [style=interior] (13.center) to (14.center);
		\draw [style=boundary] (40.center) to (41.center);
		\draw [style=boundary] (41.center) to (42.center);
		\draw [style=boundary] (6.center) to (40.center);
		\draw [style={green_thin}] (32.center) to (26.center);
		\draw [style={green_thin}] (26.center) to (28.center);
		\draw [style={green_thin}] (28.center) to (9.center);
		\draw [style={green_thin}] (9.center) to (19.center);
		\draw [style={green_thin}] (19.center) to (1.center);
		\draw [style={green_thin}] (1.center) to (16.center);
		\draw [style={green_thin}] (16.center) to (13.center);
		\draw [style={green_thin}] (13.center) to (24.center);
		\draw [style={green_thin}] (24.center) to (4.center);
		\draw [style={green_thin}] (4.center) to (27.center);
		\draw [style={green_thin}] (27.center) to (26.center);
		\draw [style={red_thin}] (32.center) to (31.center);
		\draw [style={red_thin}] (31.center) to (5.center);
		\draw [style={red_thin}] (5.center) to (23.center);
		\draw [style={red_thin}] (23.center) to (14.center);
		\draw [style={red_thin}] (14.center) to (17.center);
		\draw [style={red_thin}] (17.center) to (2.center);
		\draw [style={red_thin}] (2.center) to (20.center);
		\draw [style={red_thin}] (20.center) to (10.center);
		\draw [style={red_thin}] (10.center) to (29.center);
		\draw [style={blue_thin}] (6.center) to (5.center);
		\draw [style={blue_thin}] (5.center) to (4.center);
		\draw [style={blue_thin}] (4.center) to (3.center);
		\draw [style={pink_thin}] (25.center) to (28.center);
		\draw [style={pink_thin}] (28.center) to (29.center);
		\draw [style={pink_thin}] (29.center) to (36.center);
		\draw [style={purple_thin}] (25.center) to (27.center);
		\draw [style={purple_thin}] (27.center) to (31.center);
		\draw [style={purple_thin}] (56.center) to (41.center);
		\draw [style={purple_thin}] (31.center) to (56.center);
		\draw [style={interior}] (11.center) to (10.center);
		\draw [style={interior}] (10.center) to (9.center);
		\draw [style={interior}] (9.center) to (8.center);
		\draw [style={boundary}](42.center) to (54.center);
		\draw [style={boundary}](54.center) to (47.center);
		\draw [style={orange_thin}] (32.center) to (54.center);
		\draw [style={lime_thin}] (30.center) to (47.center);
		\draw [style={brown_thin}] (55.center) to (42.center);
		\draw [style={cyan_thin}] (32.center) to (55.center);
		\draw [style={magenta_thin}] (32.center) to (30.center);
		\draw [style={yellow_thin}] (56.center) to (6.center);
		\draw [style={yellow_thin}] (55.center) to (56.center);
		\draw [style={green_thin}] (26.center) to (58.center);
		\draw [style={red_thin}] (58.center) to (30.center);
		\draw [style={red_thin}] (29.center) to (58.center);
		\draw [style={green_thin}] (58.center) to (57.center);
		\draw [style={boundary}](11.center) to (36.center);
		\draw [style={boundary}](36.center) to (57.center);
		\draw [style={boundary}](57.center) to (47.center);

		\fill[bEV] (6) circle [radius=0.05];
		\fill[EV] (55) circle [radius=0.05];
		\fill[EV] (30) circle [radius=0.05];
		\fill[EV] (32) circle [radius=0.05];
		\fill[bEV] (25) circle [radius=0.05];
\end{tikzpicture}}
    \caption{Illustration of the interface tracing procedure. From each EV all outgoing edges are traced as illustrated until another EV or a boundary is hit. }
    \label{fig:Illustration_patches_from_mesh2}
    \end{subfigure}
    \hfill
    \begin{subfigure}[t]{0.3\linewidth}
    \centering
    \resizebox{\linewidth}{!}{\begin{tikzpicture}
[
        ,boundary/.style={thick, draw=black},
        ,interior/.style={thin, draw=gray},
		,green_thin/.style={thin, draw=green},
		,red_thin/.style={thin, draw=red},
		,blue_thin/.style={thin, draw=blue},
		,orange_thin/.style={thin, draw=orange},
		,yellow_thin/.style={thin, draw=yellow},
		,pink_thin/.style={thin, draw=pink},
		,purple_thin/.style={thin, draw=purple},
		,brown_thin/.style={thin, draw=brown},
		,lime_thin/.style={thin, draw=lime},
		,cyan_thin/.style={thin, draw=cyan},
		,magenta_thin/.style={thin, draw=magenta},
		,EV/.style={fill=gray},
		,bEV/.style={fill=black},
]
		\node [] (0) at (0, 3) {};
		\node [] (1) at (0, 3.5) {};
		\node [] (2) at (0, 3.95) {};
		\node [] (3) at (0, 2) {};
		\node [] (4) at (0, 1.5) {};
		\node [] (5) at (0, 1) {};
		\node [] (6) at (0, 0.5) {};
		\node [] (7) at (0, 4.5) {};
		\node [] (8) at (0.5, 2.5) {};
		\node [] (9) at (1, 2.5) {};
		\node [] (10) at (1.5, 2.5) {};
		\node [] (11) at (2, 2.5) {};
		\node [] (12) at (-0.5, 2.5) {};
		\node [] (13) at (-1, 2.5) {};
		\node [] (14) at (-1.45, 2.5) {};
		\node [] (15) at (-2, 2.5) {};
		\node [] (16) at (-0.675, 3.175) {};
		\node [] (17) at (-1, 3.5) {};
		\node [] (18) at (-1.425, 3.925) {};
		\node [] (19) at (0.675, 3.175) {};
		\node [] (20) at (1, 3.5) {};
		\node [] (21) at (1.425, 3.925) {};
		\node [] (22) at (-1.425, 1.075) {};
		\node [] (23) at (-1, 1.5) {};
		\node [] (24) at (-0.675, 1.8) {};
		\node [] (25) at (0.35, 2.15) {};
		\node [] (26) at (1, 1.5) {};
		\node [] (27) at (0.5, 1.5) {};
		\node [] (28) at (1, 2) {};
		\node [] (29) at (1.5, 2) {};
		\node [] (30) at (1.5, 1.15) {};
		\node [] (31) at (0.5, 1) {};
		\node [] (32) at (1, 1) {};
		\node [] (36) at (2, 2) {};
		\node [] (37) at (-0.325, 2.825) {};
		\node [] (38) at (-0.325, 2.175) {};
		\node [] (39) at (0.325, 2.825) {};
		\node [] (40) at (0, 0) {};
		\node [] (41) at (0.5, 0) {};
		\node [] (42) at (1, 0) {};
		\node [] (47) at (2, 1) {};
		\node [] (54) at (1.5, 0.5) {};
		\node [] (55) at (0.825, 0.5) {};
		\node [] (56) at (0.5, 0.5) {};
		\node [] (57) at (2, 1.5) {};
		\node [] (58) at (1.5, 1.5) {};

		\draw [style=boundary] (11.center) to (21.center);
		\draw [style=boundary] (21.center) to (7.center);
		\draw [style=boundary] (7.center) to (18.center);
		\draw [style=boundary] (18.center) to (15.center);
		\draw [style=boundary] (15.center) to (22.center);
		\draw [style=boundary] (22.center) to (6.center);
		\draw [style=boundary] (8.center) to (25.center);
		\draw [style=boundary] (25.center) to (3.center);
		\draw [style=interior] (19.center) to (20.center);
		\draw [style=boundary] (3.center) to (38.center);
		\draw [style=boundary] (38.center) to (12.center);
		\draw [style=boundary] (12.center) to (37.center);
		\draw [style=boundary] (37.center) to (0.center);
		\draw [style=boundary] (0.center) to (39.center);
		\draw [style=boundary] (39.center) to (8.center);
		\draw [style=interior] (24.center) to (38.center);
		\draw [style=interior] (12.center) to (13.center);
		\draw [style=interior] (37.center) to (16.center);
		\draw [style=interior] (0.center) to (1.center);
		\draw [style=interior] (39.center) to (19.center);
		\draw [style=interior] (16.center) to (17.center);
		\draw [style=interior] (20.center) to (21.center);
		\draw [style=interior] (1.center) to (2.center);
		\draw [style=interior] (2.center) to (7.center);
		\draw [style=interior] (17.center) to (18.center);
		\draw [style=interior] (14.center) to (15.center);
		\draw [style=interior] (23.center) to (22.center);
		\draw [style=interior] (23.center) to (24.center);
		\draw [style=interior] (13.center) to (14.center);
		\draw [style=boundary] (40.center) to (41.center);
		\draw [style=boundary] (41.center) to (42.center);
		\draw [style=boundary] (6.center) to (40.center);
		\draw [style={green_thin}] (32.center) to (26.center);
		\draw [style={green_thin}] (26.center) to (28.center);
		\draw [style={green_thin}] (28.center) to (9.center);
		\draw [style={green_thin}] (9.center) to (19.center);
		\draw [style={green_thin}] (19.center) to (1.center);
		\draw [style={green_thin}] (1.center) to (16.center);
		\draw [style={green_thin}] (16.center) to (13.center);
		\draw [style={green_thin}] (13.center) to (24.center);
		\draw [style={green_thin}] (24.center) to (4.center);
		\draw [style={green_thin}] (4.center) to (27.center);
		\draw [style={green_thin}] (27.center) to (26.center);
		\draw [style={red_thin}] (32.center) to (31.center);
		\draw [style={red_thin}] (31.center) to (5.center);
		\draw [style={red_thin}] (5.center) to (23.center);
		\draw [style={red_thin}] (23.center) to (14.center);
		\draw [style={red_thin}] (14.center) to (17.center);
		\draw [style={red_thin}] (17.center) to (2.center);
		\draw [style={red_thin}] (2.center) to (20.center);
		\draw [style={red_thin}] (20.center) to (10.center);
		\draw [style={red_thin}] (10.center) to (29.center);
		\draw [style={blue_thin}] (6.center) to (5.center);
		\draw [style={blue_thin}] (5.center) to (4.center);
		\draw [style={blue_thin}] (4.center) to (3.center);
		\draw [style={pink_thin}] (25.center) to (28.center);
		\draw [style={pink_thin}] (28.center) to (29.center);
		\draw [style={pink_thin}] (29.center) to (36.center);
		\draw [style={purple_thin}] (25.center) to (27.center);
		\draw [style={purple_thin}] (27.center) to (31.center);
		\draw [style={purple_thin}] (56.center) to (41.center);
		\draw [style={purple_thin}] (31.center) to (56.center);
		\draw [style={interior}] (11.center) to (10.center);
		\draw [style={interior}] (10.center) to (9.center);
		\draw [style={interior}] (9.center) to (8.center);
		\draw [style={boundary}](42.center) to (54.center);
		\draw [style={boundary}](54.center) to (47.center);
		\draw [style={orange_thin}] (32.center) to (54.center);
		\draw [style={lime_thin}] (30.center) to (47.center);
		\draw [style={brown_thin}] (55.center) to (42.center);
		\draw [style={cyan_thin}] (32.center) to (55.center);
		\draw [style={magenta_thin}] (32.center) to (30.center);
		\draw [style={yellow_thin}] (56.center) to (6.center);
		\draw [style={yellow_thin}] (55.center) to (56.center);
		\draw [style={green_thin}] (26.center) to (58.center);
		\draw [style={red_thin}] (58.center) to (30.center);
		\draw [style={red_thin}] (29.center) to (58.center);
		\draw [style={green_thin}] (58.center) to (57.center);
		\draw [style={boundary}](11.center) to (36.center);
		\draw [style={boundary}](36.center) to (57.center);
		\draw [style={boundary}](57.center) to (47.center);

		\fill[bEV] (6) circle [radius=0.05];
		\fill[EV] (55) circle [radius=0.05];
		\fill[EV] (30) circle [radius=0.05];
		\fill[EV] (32) circle [radius=0.05];
		\fill[bEV] (25) circle [radius=0.05];

		\fill[gray,opacity=0.1] (40.center)--(41.center)--(56.center)--(6.center)--cycle;
		\fill[gray,opacity=0.1] (41.center)--(42.center)--(55.center)--(56.center)--cycle;
		\fill[gray,opacity=0.1] (42.center)--(54.center)--(32.center)--(55.center)--cycle;
		\fill[gray,opacity=0.1] (54.center)--(47.center)--(30.center)--(32.center)--cycle;
		\fill[gray,opacity=0.1] (47.center)--(57.center)--(58.center)--(30.center)--cycle;
		\fill[gray,opacity=0.1] (58.center)--(57.center)--(36.center)--(29.center)--cycle;
		\fill[gray,opacity=0.1] (29.center)--(36.center)--(11.center)--(10.center)--cycle;

		\fill[gray,opacity=0.1] (6.center)--(56.center)--(31.center)--(5.center)--cycle;
		\fill[gray,opacity=0.1] (56.center)--(55.center)--(32.center)--(31.center)--cycle;

		\fill[gray,opacity=0.1] (5.center)--(31.center)--(27.center)--(4.center)--cycle;
		\fill[gray,opacity=0.1] (31.center)--(32.center)--(26.center)--(27.center)--cycle;
		\fill[gray,opacity=0.1] (32.center)--(30.center)--(58.center)--(26.center)--cycle;
		\fill[gray,opacity=0.1] (26.center)--(58.center)--(29.center)--(28.center)--cycle;
		\fill[gray,opacity=0.1] (28.center)--(29.center)--(10.center)--(9.center)--cycle;

		\fill[gray,opacity=0.1] (4.center)--(27.center)--(25.center)--(3.center)--cycle;
		\fill[gray,opacity=0.1] (27.center)--(26.center)--(28.center)--(25.center)--cycle;
		\fill[gray,opacity=0.1] (25.center)--(28.center)--(9.center)--(8.center)--cycle;

		\fill[gray,opacity=0.2] (8.center)--(9.center)--(19.center)--(1.center)--(16.center)--(13.center)--(24.center)--(4.center)--(3.center)--(38.center)--(12.center)--(37.center)--(0.center)--(39.center)--cycle;
		\fill[gray,opacity=0.3] (9.center)--(10.center)--(20.center)--(2.center)--(17.center)--(14.center)--(23.center)--(5.center)--(4.center)--(24.center)--(13.center)--(16.center)--(1.center)--(19.center)--cycle;
		\fill[gray,opacity=0.4] (10.center)--(11.center)--(21.center)--(7.center)--(18.center)--(15.center)--(22.center)--(6.center)--(5.center)--(23.center)--(14.center)--(17.center)--(2.center)--(20.center)--cycle;
\end{tikzpicture}}
    \caption{Result of interface tracing from all the EVs. Every patch is now bounded by a a set of boundary and traced interface curves. All patch corners are corners where a traced interface and/or a boundary edge form a corner. Along the hole, different patches are indicated with different shades of gray. In the part bottom-right of the hole, every face forms a patch, since all traced curves denoted by colors intersect with other traced curves.}
    \label{fig:Illustration_patches_from_mesh3}
    \end{subfigure}
    \caption{Procedure to find a multi-patch segmentation from a given mesh. The original mesh in (a) has 46 vertices, 81 edges and 45 faces and the final multi-patch (c) has 20 patches.}
    \label{fig:Illustration_patches_from_mesh}
\end{figure}
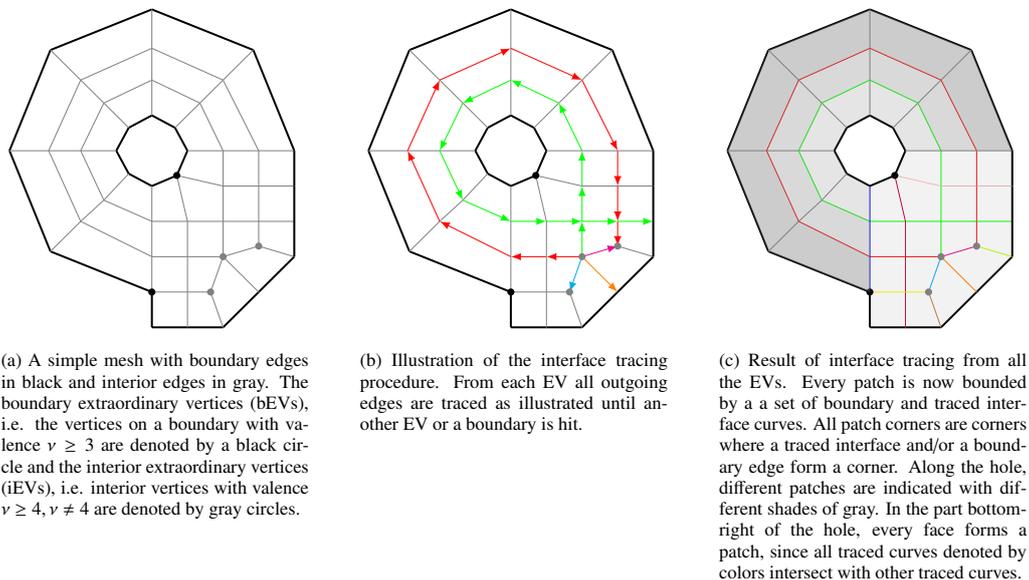

Firstly, a \emph{quadrilateral mesh} (\emph{quad mesh}) is a mesh of quadrilateral elements, representing a (planar) surface geometry. The quadrilaterals can be represented by tensor B-splines of any degree which can be mapped onto a parametric unit-square. Typically, when the tensor B-spline quadrilaterals have different sizes in different directions or even different refinement levels, assemblies of these \emph{patches} are typically referred to as \emph{multi-patches}. An example of a multi-patch is given in \cref{fig:segmentation}. The conversion of a quad-mesh with many elements to a multi-patch with a smaller number of patches derived from groups of elements can be done using the procedure described in \cref{fig:Illustration_patches_from_mesh}. Here, a half-edge mesh is traversed and elements are collected into groups corresponding to final patches. The vertices of the elements in one group (i.e. patch) form the control net of the bi-linear patch.\\

Secondly, for parametrically smooth constructions, different classes of vertices are considered. For so-called \emph{extraordinary vertices} (\emph{EVs}) , these constructions typically are different. An \emph{interior extraordinary vertex} (\emph{interior EV}) is a vertex on a quad mesh on which three or more than four patches meet. The number of patches coming together at a vertex is referred to as the \emph{valence}, denoted by $\nu$. Furthermore a \emph{boundary extraordinary vertex} (\emph{boundary EV}) is a vertex on the boundary of the quad mesh with valence $\nu\geq3$. For geometrically smooth constructions, the construction depends on the geometry around the vertex rather than the valence of the vertex. Hence, for these constructions the notion of EVs is irrelevant.\\

Lastly, a refinement of a spline space is called \emph{nested} if the refined spline space is fully contained in the unrefined space. As a consequence, the geometry is exact under element refinement, which is beneficial from an analysis point of view.

\subsection{Analysis-suitable $G^1$}\label{subsubsec:ASG1}

The analysis-suitable $G^1$ (AS-$G^1$) construction is a novel approach in isogeometric analysis that was introduced for planar geometries and surfaces in \cite{Collin2016}, but a construction which extends \cite{Kapl2019} for planar domains to surfaces is detailed in \cite{Farahat2023}. This construction ensures that basis functions at interfaces have $C^1$ continuity, while basis functions at vertices have $C^2$ continuity. The approach is based on the concept that $G^k$-smooth surfaces can produce $C^k$-smooth isogeometric functions \cite{Groisser2015}. When dealing with general $C^0$-matching multi-patch domains, the so-called AS-$G^1$ conditions must be satisfied to ensure optimal approximation. If these conditions are met, a $C^1$-smooth subspace of the isogeometric space can be constructed, which is sufficiently large. Such geometries are referred to as analysis-suitable geometries. However, it should be noted that the $C^1$-smooth multi-patch isogeometric space generally depends on the geometry, as discussed in \cite{Kapl2017a}. To overcome this issue, an Argyris-like space was proposed in \cite{Kapl2019}, which has a dimension that is independent of the geometry.\\

Given an interface between two patches, the $C^1$ continuity condition at the interface is defined by a linear combination of tangent vectors and transversal derivatives, which is referred to as gluing data \cite{Collin2016}. The $C^1$ smooth basis functions at the interface, or more generally at the edge, can be described by the first order Taylor expansion of the trace and the transversal derivative. It is shown in \cite{Collin2016}, that the ideal choice for the space-representation of the trace and transversal derivative is $\vb{\mathcal{S}} (\vb p, \vb r-1, \vb h)$\footnote{The notation $\vb{\mathcal{S}} (\vb p = (p,p), \vb r = (r,r), \vb h = (h,h))$ indicates a two-dimensional spline space with $p$ as the polynomial degree, $r$ as the regularity and $h$ as the mesh size in both directions.} and $\vb{\mathcal{S}} (\vb p-1, \vb r-2, \vb h)$, respectively. These basis functions have local support and are linearly independent, but they depend on the gluing data and, therefore, on the geometry reparameterization itself. To ensure that the basis functions form a $C^1$-smooth subspace of the isogeometric space, and to maintain the nestedness of the spline spaces, it is necessary to have gluing data as a linear function which fulfils all analysis-suitable geometries. For instance, all bi-linear patches meet this requirement. However, if a geometry is not analysis-suitable, it can be reparameterized using the technique presented in \cite{Kapl2018}.\\

For any vertices in the quad mesh, to describe the $C^1$ condition is not that straight-forward. In order to keep it general, the vertex basis functions is constructed by the $C^2$ interpolation using the $C^1$ basis functions from the corresponding edges. As a consequence, the vertex basis functions also have local support and are linearly independent.\\

Summarizing, the AS-$G^1$ construction can be constructed by three different, linearly independent sub-spaces: interior, edge and vertex space. They can be described as follows:
\begin{itemize}
    \item Interior space: basis functions that have zero values and derivatives on the patch edges and vertices.
    \item Interface space: basis functions that have vanishing function values up to the second derivatives at the vertices.
    \item Vertex space: $C^2$ interpolating functions at the vertex, i.e., basis functions that have non-vanishing $C^2$ data at the vertex.
\end{itemize}
The AS-$G^1$ construction with the interface and vertex constructions as described above are fully $C^1$ over the whole domain. In addition, the AS-$G^1$ construction can only be constructed when the degree of the basis is $p\geq3$ and the regularity is reduced as $r\leq p-2$.\\

\Cref{fig:continuity_ASG1} presents a local region around and EV with valence five with line styles indicating different continuity levels on patch or element boundaries (see the caption of \cref{fig:continuity0}). For the AS-$G^1$ construction, the continuity at the vertex is $C^2$ by construction. Furthermore, the continuity at the interior element interfaces is $C^{p-2}$ due to the restriction on keeping the isoparametric concept. Lastly, since the AS-$G^1$ construction  provides a $G^1$ surface, the patch interfaces are $C^1$ by construction \cite{Groisser2015}.\\

In sum, the core ideas behind the AS-$G^1$ construction are as follows:
\begin{itemize}
	\item \textbf{Degree, regularity, continuity}\\ The spline space is fully $C^1$, hence suitable to solve fourth-order problems. However, the computation of the space requires analysis-suitability of the parameterization as well as degree $p\geq3$ and regularity $r\leq p-2$ for the basis functions.
	\item \textbf{Limitations on construction}\\ The space can be constructed on fully unstructured quadrilateral meshes with both interior and boundary extraordinary vertices. The construction of the basis functions is independent of the location or valence of the EVs. However, the analysis-suitability condition imposes a requirement on the geometries on which the construction can be constructed. Furthermore, the geometry parameterization is not changed.
	\item \textbf{Nestedness}\\ The spline spaces are nested.
	\item \textbf{Refinement procedure}\\ Refinement procedure is standard (by knot insertion) since the parameterization does not change.
\end{itemize}

\subsection{Approximate $C^1$}\label{subsubsec:ApproxC1}
The Approximate $C^1$ construction \cite{Weinmuller2021} provides, as the name suggests, approximately $C^1$ continuity on interfaces and vertices, more precisely the construction provides $C^1$ continuity in the refinement limit. The Approx.~$C^1$ construction shares similarities with the AS-$G^1$ construction, but the main difference between the construction of the Approx.~$C^1$ and the AS-$G^1$ spaces is that it relaxes the AS-$G^1$ condition on the geometry, i.e., it allows geometries with non-linear gluing data. In fact, the exact gluing data are splines with higher polynomial degree and lower regularity or even piece-wise rational. As a consequence, trying to extend the construction for AS-$G^1$ parameterizations directly to non-AS-$G^1$ geometries yields complicated basis functions that are challenging to evaluate and integrate accurately. To overcome this issue and obtain a construction with more easily definable basis functions, the gluing data are approximated. However, this approximation means that the $C^1$ condition is no longer satisfied exactly but only approximately.\\

By utilizing the approximation of the gluing data, the Approximate $C^1$ construction incorporates the concept of different spline spaces found in the AS-$G^1$ construction. In this case, the interior, vertex, and interface basis functions fulfill the same conditions as in the AS-$G^1$ construction, but the degree and regularity differ between these spaces. Specifically, the sub-spaces for the AS-$G^1$ construction have $p\geq3$ and $r\leq p-2$, while the Approximate $C^1$ construction employs an interior space with $p\geq3$ and $r\leq p-1$, along with vertex and interface spaces that have locally reduced smoothness based on the approximation of the gluing data. Consequently, on the one hand the Approximate $C^1$ construction restores the potential for maximal smoothness of isogeometric functions in the refinement limit, but the nestedness of the basis is lost. On the other hand, the approximation of the gluing data in the Approximate $C^1$ construction does not require analysis-suitability for the optimal convergence rate, unlike the AS-$G^1$ construction. This feature makes the method applicable to more complex geometries. When the Approximate $C^1$ construction is applied to an analysis-suitable geometry with $p\geq3$ and $r\leq p-2$, and the gluing data approximation is exact, the construction becomes equivalent to the AS-$G^1$ construction.\\

\begin{figure}
 \centering
 \begin{subfigure}{0.45\linewidth}
  \centering
  \includegraphics[width=\linewidth]{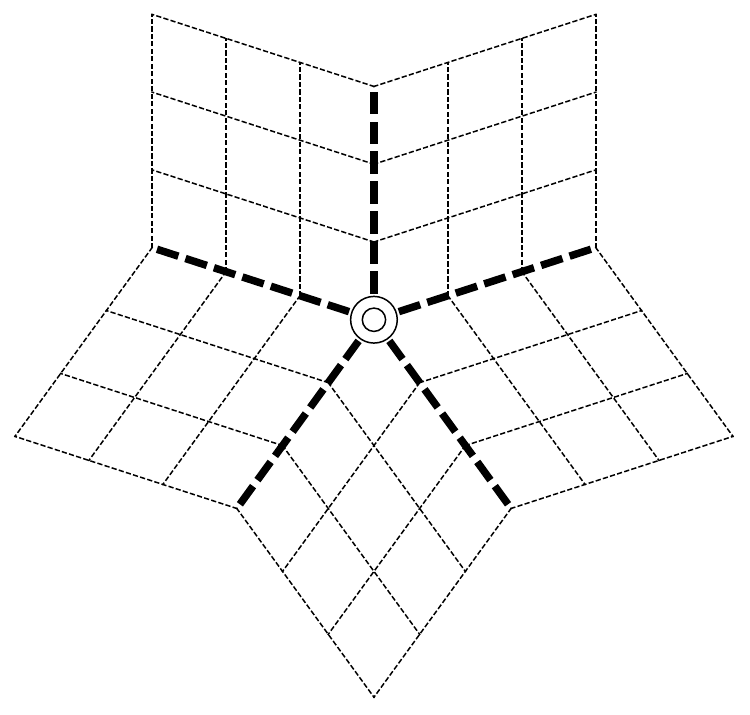}
  \caption{AS-$G^1$}
  \label{fig:continuity_ASG1}
\end{subfigure}
\hfill
 \begin{subfigure}{0.45\linewidth}
  \centering
  \includegraphics[width=\linewidth]{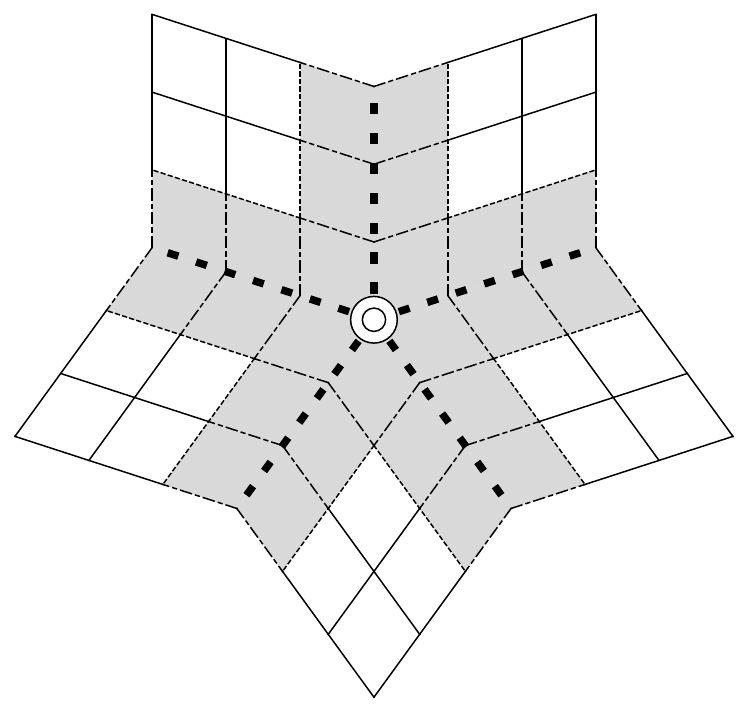}
  \caption{Approx. $C^1$}
  \label{fig:continuity_approxC1}
\end{subfigure}
\caption{Schematic representation of the continuity across element boundaries and patch interfaces for the (\subref{fig:continuity_ASG1}) AS-$G^1$ construction, (\subref{fig:continuity_approxC1}) Approx.~$C^1$ constructions. Thin lines indicate element boundaries and thick lines indicate patch interfaces. Solid lines represent $C^{p-1}$ continuity, dashed lines represent $C^{p-2}$ continuity, thick dashed lines represent $C^1$ interfaces and loosely dashed lines represent approximate $C^1$ interfaces. A double lined circle represent a $C^2$ continuous vertex, a filled circle represent a singular vertex and a white filled circled with a single line represents a $C^1$ continuous vertex. The gray shaded area for the Approx.~$C^1$ represents local reduced continuity.}
\label{fig:continuity0}
\end{figure}

\Cref{fig:continuity_approxC1} presents a local region around and EV with valence five with line styles indicating different continuity levels on patch or element boundaries (see the caption of \cref{fig:continuity0}). For the Approx.~$C^1$ construction on a fully smooth basis ($p\geq3$ and $r=p-1$), the interior basis recovers full smoothness on element boundaries, hence $C^{p-1}$ continuity. In the shaded region around the interfaces and the EV, the continuity is locally reduced by construction of the locally reduced continuous subspace and the approximation of the gluing data. Similar to the AS-$G^1$ construction, the continuity on the EV is $C^2$ by construction and the element boundaries are $C^1$ approximately.\\

In sum, the core ideas behind the Approx. $C^1$ construction are as follows:
\begin{itemize}
	\item \textbf{Degree, regularity, continuity}\\ The spline space is approximately $C^1$ and fully $C^1$ in the limit of refinement. This makes the spline space suitable to solve fourth-order problems. Contrary to the AS-G1 construction, the spline space approximates the gluing data, allowing maximal smoothness in the interior space ($r=p-1$) for degrees $p\geq3$. 
	\item \textbf{Limitations on construction}\\ The space can be constructed on fully unstructured quadrilateral meshes with both interior and boundary extraordinary vertices. The construction of the basis functions is independent of the location or valence of the vertices. Contrary to AS-$G^1$ the analysis-suitability condition is not needed. However, the construction requires a $G^1$ condition at the interfaces of surfaces.
	\item \textbf{Nestedness}\\ The spline spaces are not nested.
	\item \textbf{Refinement procedure}\\ Refinement procedure is standard since the parameterization does not change.
\end{itemize}

\subsection{D-patch}\label{subsubsec:dpatch}

The relative ease of imposing parametric smoothness for splines has led to the development of degenerate Bezier patches, or D-patches \cite{Reif1997}, which can be used to build $C^1$ smooth splines on unstructured quadrilateral meshes with no boundary extraordinary vertices. The constructions can be formulated for splines of any bi-degree \cite{Hughes2021}, and there are no restrictions on their smoothness in the locally-structured regions of the mesh. In the locally-unstructured regions of the mesh (i.e., in a neighbourhood of an extraordinary vertex), the splines are $C^1$ smooth and first-order degenerate. Note that this degeneracy means that the spline spaces are not necessarily $H^2$-conforming, but numerical evidence shows that they can still be used to solve fourth-order problems.\\

Specifically, imposition of strong $C^1$ smoothness around an extraordinary vertex requires that the splines vanish up to first order at the extraordinary vertex. This degeneracy trivially implies matching first derivatives at the extraordinary vertex (since all of them vanish) but does not imply $C^1$ smoothness of the resulting spline functions and the geometries built using them. As shown in \cite{Reif1997}, additional conditions can be imposed upon certain higher-order mixed derivatives to ensure this desired $C^1$ smoothness. Furthermore, the effect of these additional constraints can be localised to a neighbourhood of the extraordinary vertex by imposing them on a subdivided representation of the splines \cite{Nguyen2016}. This means that a patch-based representation of $C^1$ D-patch splines takes functions that are in $\vb{\mathcal{S}} (\vb p, \vb r, \vb h/2)$ on each patch, where almost all basis functions are in $\vb{\mathcal{S}} (\vb p, \vb r, \vb h)$, except a few basis functions supported in a neighbourhood of extraordinary points (the number of basis functions depends on the valence).\\

The D-patch construction allows for nested refinements of the spline spaces \cite{Reif1997}. If different orders of smoothness are being imposed in locally-structured and locally-unstructured regions of the mesh, then nested refinements produce spline spaces with a higher number of basis functions supported in the vicinity of extraordinary points (the number depends on the refinement-level), see \cite{Toshniwal2017} for instance. On the other hand, a patch-based approach allows for a simpler implementation by limiting the smoothness across patch interfaces to $C^1$; the smoothness in patch-interiors can still be arbitrarily chosen. However, special care should be taken when using D-patches with nested refinements -- the degeneracy of the splines near extraordinary vertices means that, with mesh refinements, the shape regularity of the mesh starts to worsen with refinements and the finite element matrices become very ill-conditioned.\\

In sum, the core ideas behind the D-patch spline construction are the following:
\begin{itemize}
    \item \textbf{Degree, regularity, continuity}\\ The spline space is fully $C^1$. In general, the degeneracy of derivatives means that the spaces are $H^2$-nonconforming, however numerical evidence supports their use in solving fourth-order problems. The construction can be formulated for splines of any degree and the smoothness away from extraordinary vertices can be chosen arbitrarily.
    \item \textbf{Limitations on construction}\\ The space can be constructed on unstructured quadrilateral meshes with no boundary extraordinary vertices.
    \item \textbf{Nestedness}\\ The spline spaces can be refined in a nested manner, however the resulting mesh have poor shape regularity and the corresponding finite element matrices may be very ill-conditioned.
    \item \textbf{Refinement procedure}\\ Refinement procedures can be derived from standard B-spline knot insertion.
\end{itemize}

\subsection{Almost $C^1$}\label{subsubsec:almostC1}
Almost-$C^1$ splines are defined on a general, conforming quadrilateral mesh. They are piece-wise biquadratic and possess mixed smoothness, i.e., they are $C^1$ in regular regions, while the smoothness near extraordinary vertices, i.e., vertices with valence different from four, is reduced. To be precise, they are $C^1$ smooth at all vertices (including extraordinary vertices) and across all edges except for the ones emanating from an extraordinary vertex. Moreover, while they are defined to be biquadratic on all regular elements, they are piece-wise biquadratic splines (with one inner knot in each direction) on all elements that are neighboring an extraordinary vertex. Details can be found in \cite{Takacs2023}. As a consequence, a patch based representation of Almost-$C^1$ splines takes functions that are in $\vb{\mathcal{S}} (\vb 2, \vb 1, \vb h/2)$ on each patch, where almost all basis functions are in $\vb{\mathcal{S}} (\vb 2, \vb 1, \vb h)$, except a few basis functions supported in a 1-ring neighbourhood of extraordinary points (the number depends on the valence). \\

A central feature of Almost-$C^1$ splines is the mixed smoothness imposition described above. In particular, this choice of mixed smoothness only depends on the current refinement level of the mesh. That is, standard $C^1$-smoothness is enforced across all edges at the current refinement level except the ones that are incident upon extraordinary vertices, where only $C^0$ smoothness is enforced. Additionally, these smoothness conditions are combined with $G^1$ smoothness imposition at each extraordinary vertex. This means that almost-$C^1$ splines do not yield nested spaces when refining. As a result, the refinement process essentially amounts to a projection of coarse Almost-$C^1$ splines onto the refined Almost-$C^1$ spline space. This projection can be chosen in many different ways and can have a significant impact on the limit surface description as well as isogeometric simulations using these spaces. In \cite{Takacs2023} a smoothing and refinement procedure is proposed that results in a $C^1$-smooth limit surface for sufficiently regular input data.\\

Let us briefly summarize the refinement procedure here. We assume that we are given a quad mesh and associate a control point with each face of the mesh. The initial smoothing step guarantees that all control points associated to the one ring around an extraordinary vertex are coplanar. Having given such an initial control point grid, we then refine the geometry using explicit subdivision rules as specified in \cite{Toshniwal2022,Takacs2023}. The rules are the same as for quadratic tensor-product B-splines in regular regions and maintain the coplanarity near extraordinary vertices.\\

In sum, the core ideas behind the Almost-$C^1$ spline construction are the following:
\begin{itemize}
	\item \textbf{Degree, regularity, continuity}\\ The spline space locally reproduces biquadratic polynomials and it is sufficiently smooth to be able to solve fourth order problems.
	\item \textbf{Limitations on construction}\\ The splines can be constructed on fully unstructured quadrilateral meshes, in particular, those that contain both interior and boundary extraordinary vertices.
	\item \textbf{Nestedness}\\ Since the spaces are not nested, the convergence behavior of Almost-$C^1$ splines depends on how the geometry parameterization is refined.
	\item \textbf{Refinement procedure}\\ An initial geometry and a refinement procedure can be constructed in such a way, that the limit geometry parameterization is normal continuous everywhere.
\end{itemize}

Thus, the concept introduced in \cite{Takacs2023} is quite flexible, since the initial smoothing procedure and the refinement procedure are not unique and can be tailored to the needs coming from geometric modeling, e.g., one may want to reproduce Doo-Sabin subdivision surfaces, thus having to modify the subdivision rule for refinement accordingly. The spline space that is introduced on each refinement level can be seen as a simple hole-filling construction, which is sufficient for numerical analysis.

\begin{figure}
 \centering
 \begin{subfigure}{0.45\linewidth}
  \centering
  \includegraphics[width=\linewidth]{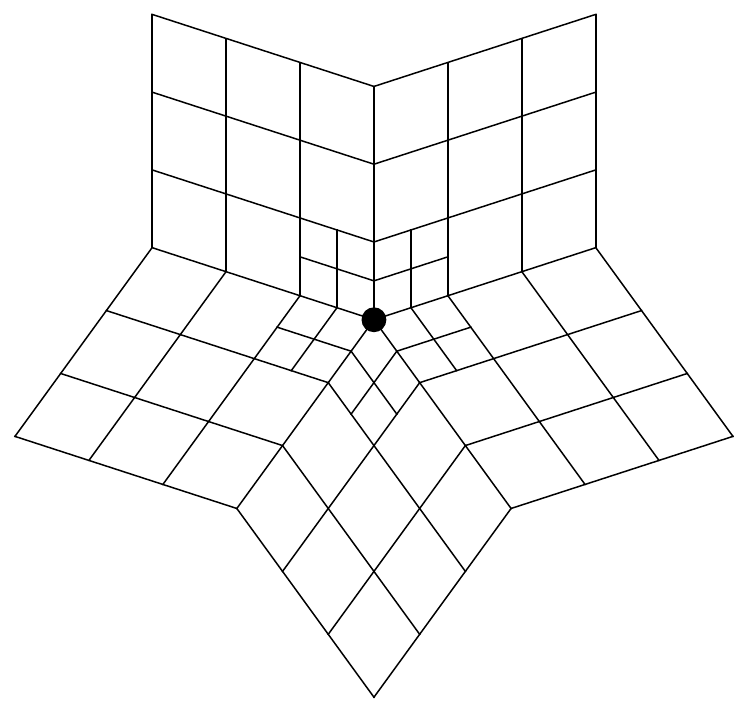}
 \caption{D-Patch}
 \label{fig:continuity_DPatch}
\end{subfigure}
\hfill
 \begin{subfigure}{0.45\linewidth}
  \centering
  \includegraphics[width=\linewidth]{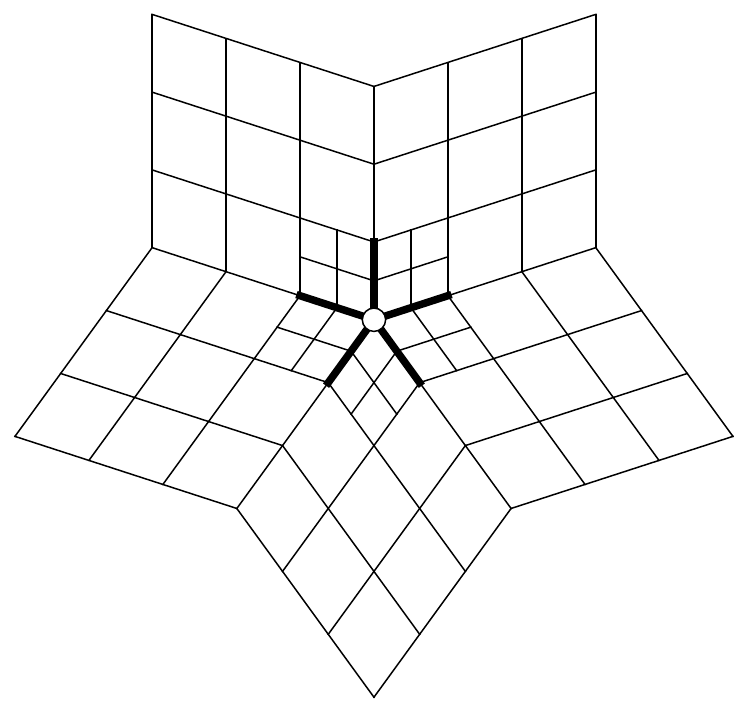}
 \caption{Almost-$C^1$}
 \label{fig:continuity_almostC1}
\end{subfigure}
\caption{Schematic representation of the continuity across element boundaries and patch interfaces for the (\subref{fig:continuity_DPatch}) D-Patch and (\subref{fig:continuity_almostC1}) Almost-$C^1$. Line styles are as in \cref{fig:continuity0}.}
\label{fig:continuity}
\end{figure}

\subsection{Conclusions}
In this section, a summary of the construction and the properties of the analysis-suitable $G^1$ (AS-$G^1$), the approximate $C^1$ (Approx.~$C^1$), the degenerate patches (D-patch) and the Almost-$C^1$ methods have been provided, referring to the relevant publications for the mathematical details. For each method, comments have been provided on the degree, regularity and continuity of the space, on the limitations of the construction in terms of the quadrilateral mesh, on nestedness for refinement and on the refinement procedure itself. In addition, \cref{fig:continuity0,fig:continuity} provides detailed information on the local continuity of the constructions around an extraordinary vertex.\\

The aim of the qualitative analysis of the methods in this paper is to provide a comparison of a set of properties and requirements of each method and their implications on their applicability. While the subsections presented before provide a brief description of the properties of the methods and the reason behind these properties and requirements, \cref{tab:qualitative} provides a side-by-side comparison of each method based on the subsections before. In particular, the table lists the \emph{(i,ii)} requirements on degree and regularity for the constructions, \emph{iii} geometrical or topological limitations if applicable, \emph{(iv,v)} the continuity of the constructed bases in the interior and on the interfaces and element boundaries and \emph{vi} nestedness of the constructed basis.\\

\begin{table}[htp]
    \centering
    \caption{Summary of the requirements for the construction and the properties for each of the considered bases. The construction requirements include the degree and regularity of the basis used for construction as well as geometrical or topological properties of the input geometry. The properties include the continuity on interfaces, vertices and in the interior of the unstructured spline construction, as well as the nestedness property.}
    \label{tab:qualitative}
    \begin{tabular}{p{0.2\linewidth}|p{0.15\linewidth}p{0.15\linewidth}p{0.15\linewidth}p{0.15\linewidth}}
        \toprule
        \textbf{Requirements} & \textbf{AS-$G^1$} & \textbf{Approx.~$C^1$} & \textbf{D-Patch} & \textbf{Almost-$C^1$} \\ \midrule
        (\emph{i}) Degree                  & $p\geq3$              & $p\geq3$              & $p\geq3$              & $p=2$\\
        (\emph{ii}) Regularity              & $r\leq p-2$           & $r\leq p-1$           & $r\leq p-1$           & $r=1$\\
        (\emph{v}) Geometrical / topological limitations & Analysis-suitability                  & $G^2$ continuity & BEVs: $\nu\leq3$, $C^1$ continuity & $C^1$ continuity \\
        \midrule
        \textbf{Properties} & \textbf{AS-$G^1$} & \textbf{Approx.~$C^1$} & \textbf{D-Patch} & \textbf{Almost-$C^1$} \\ 
        \midrule
        (\emph{iii}) Interface \& Vertex Continuity              & $C^1$                 & $C^1$ in the limit    & $C^1$                 & $C^1$ in the limit\\
        (\emph{iv}) Interior continuity              & $C^{p-2}$             & $C^{p-1}$             & $C^{p-1}$             & $C^1$\\
        (\emph{vi}) Nestedness & Yes & No & Yes & No \\
        \bottomrule
    \end{tabular}
\end{table}

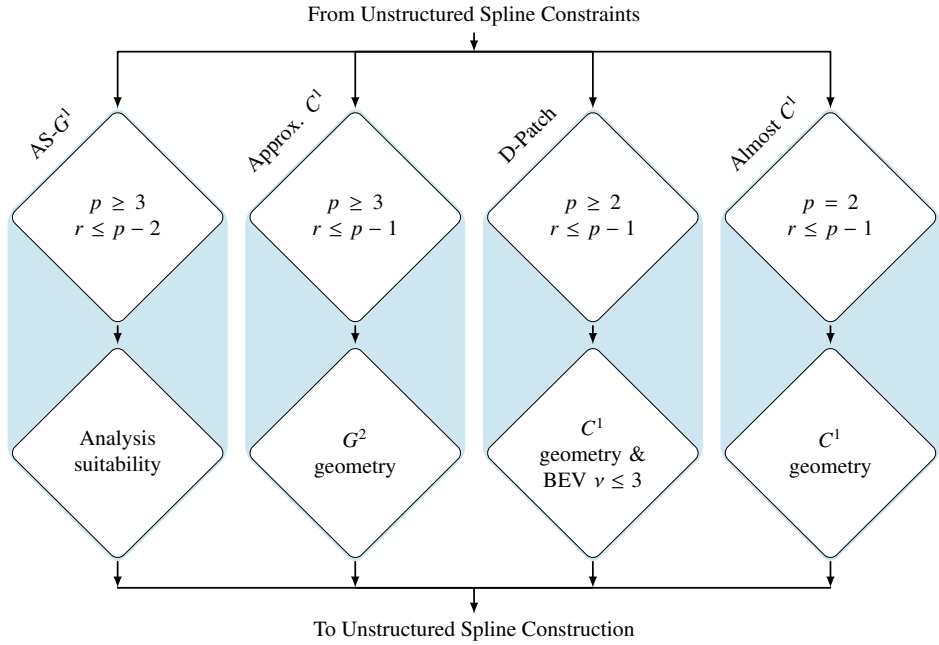
\begin{figure}
    \centering
    \resizebox{\linewidth}{!}{
%
%
%
%

\tikzstyle{arrow} = [draw, -latex,thick,line cap=round,line join=round]
\tikzstyle{line} = [draw,thick,line cap=round,line join=round]
\tikzstyle{block} = [rectangle, draw, fill=white,
    text width=10em, text centered, rounded corners, minimum height=3em]
\tikzstyle{circleblock} = [circle, draw, fill=white,
    text width=2.5em, text centered, rounded corners, minimum height=1.5em]
\tikzstyle{diamondblock} = [diamond, draw, fill=white,rounded corners,
    text width=5em, text centered, minimum height=10em, minimum width=10em]
\tikzstyle{dummyblock} = [rectangle, fill=none,
    text width=10em, text centered, rounded corners, minimum height=3em]
\tikzstyle{noblock} = [rectangle, draw=none, fill=none,
    text width=10em, rounded corners, minimum height=0em,node distance = 0.05\linewidth, text = blue]
\tikzstyle{textblock} = [text width=7em,text centered]
\tikzstyle{textblock2} = [text width=14em,text centered]
\tikzstyle{textblock4} = [text width=28em,text centered]
\tikzstyle{octagon} = [shape=regular polygon, regular polygon sides=8, draw, fill=red!20,draw=red,
    text width=12em, text centered, rounded corners, minimum height=3em,minimum width=3em,inner sep=-2em]

\begin{tikzpicture}[node distance = 0.3cm,scale=0.5]
	\node[diamondblock] (approxC1) {$p\geq 3$\\$r\leq p-1$};
	\node[diamondblock,left=of approxC1] (asG1)  {$p\geq 3$\\$r\leq p-2$};
	\node[diamondblock,right=of approxC1] (dpatch)  {$p\geq 2$\\$r\leq p-1$};
	\node[diamondblock,right=of dpatch] (almostC1)  {$p=2$\\$r\leq p-1$};
	\path (approxC1.north)--(dpatch.north) node[midway](mid){};

	\coordinate[above=of mid,yshift=0.5cm] (start);
	\coordinate[above=of start] (CADCAE);
	\path[arrow] (CADCAE.center) -- (start.center);
	\node[above] at (CADCAE)  {From Unstructured Spline Constraints};

	\node[diamondblock,below=of approxC1.south] (bapproxC1) {$G^2$ geometry};
	\node[diamondblock,below=of asG1] (basG1)  {Analysis suitability};
	\node[diamondblock,below=of dpatch] (bdpatch)  {$C^1$ geometry \& BEV $\nu\leq3$};
	\node[diamondblock,below=of almostC1] (balmostC1)  {$C^1$ geometry};

	\path[arrow] (start.center) -| (approxC1);
	\path[arrow] (start.center) -| (asG1);
	\path[arrow] (start.center) -| (dpatch);
	\path[arrow] (start.center) -| (almostC1);

	\path[arrow] (approxC1) -- (bapproxC1);
	\path[arrow] (asG1) -- (basG1);
	\path[arrow] (dpatch) -- (bdpatch);
	\path[arrow] (almostC1) -- (balmostC1);

	\node[below=of bapproxC1] (bbapproxC1) {};
	\node[below=of basG1] (bbasG1){};
	\node[below=of bdpatch] (bbdpatch){};
	\node[below=of balmostC1] (bbalmostC1){};
	\path[arrow] (bapproxC1)--(bbapproxC1.center);
	\path[arrow] (basG1)--(bbasG1.center);
	\path[arrow] (bdpatch)--(bbdpatch.center);
	\path[arrow] (balmostC1)--(bbalmostC1.center);

	\path (bbapproxC1.center)--(bbdpatch.center) node[midway](midsouth){};
	\path[line] (midsouth.center)--(bbapproxC1.center);
	\path[line] (midsouth.center)--(bbasG1.center);
	\path[line] (midsouth.center)--(bbdpatch.center);
	\path[line] (midsouth.center)--(bbalmostC1.center);

	\coordinate[below=of midsouth] (end);
	\path[arrow] (midsouth.center) -- (end.center);
	\node[below] at (end) {To Unstructured Spline Construction};

	\node[textblock,above=of dpatch.north west,rotate=45] {D-Patch};
	\node[textblock,above=of approxC1.north west,rotate=45] {Approx. $C^1$};
	\node[textblock,above=of asG1.north west,rotate=45] {AS-$G^1$};
	\node[textblock,above=of almostC1.north west,rotate=45] {Almost $C^1$};

\begin{pgfonlayer}{bg}    
	\draw[draw=col1!20,fill=col1!20][rounded corners]
([xshift=0,yshift=2]approxC1.north) --
([xshift=0,yshift=2]approxC1.east) --
([xshift=0,yshift=-2]bapproxC1.east) --
([xshift=0,yshift=-2]bapproxC1.south) --
([xshift=0,yshift=-2]bapproxC1.west) --
([xshift=0,yshift=2]approxC1.west)
--cycle;
	\draw[draw=col1!20,fill=col1!20][rounded corners]
([xshift=0,yshift=2]asG1.north) --
([xshift=0,yshift=2]asG1.east) --
([xshift=0,yshift=-2]basG1.east) --
([xshift=0,yshift=-2]basG1.south) --
([xshift=0,yshift=-2]basG1.west) --
([xshift=0,yshift=2]asG1.west)
--cycle;
	\draw[draw=col1!20,fill=col1!20][rounded corners]
([xshift=0,yshift=2]dpatch.north) --
([xshift=0,yshift=2]dpatch.east) --
([xshift=0,yshift=-2]bdpatch.east) --
([xshift=0,yshift=-2]bdpatch.south) --
([xshift=0,yshift=-2]bdpatch.west) --
([xshift=0,yshift=2]dpatch.west)
--cycle;
	\draw[draw=col1!20,fill=col1!20][rounded corners]
([xshift=0,yshift=2]almostC1.north) --
([xshift=0,yshift=2]almostC1.east) --
([xshift=0,yshift=-2]balmostC1.east) --
([xshift=0,yshift=-2]balmostC1.south) --
([xshift=0,yshift=-2]balmostC1.west) --
([xshift=0,yshift=2]almostC1.west)
--cycle;
\end{pgfonlayer}

\end{tikzpicture}}
    \caption{Inside the \emph{unstructured spline pre-processing} block from \cref{fig:IGA_Setup}. The unstructured spline requirements are depicted in diamond-shaped blocks for methods AS-$G^1$, Approx.~$C^1$, D-Patch and Almost $C^1$. The first row represents requirements on degree $p$ and regularity $r$. If not satisfied, the geometry can be projected onto a space that satisfies the requirement, or degree elevation or reduction steps can be performed together with refinement operations. The second row depicts the requirements on the geometry parameterization; these blocks can be satisfied by changing the geometry.}
    \label{fig:USconstraints}
\end{figure}

Following from \cref{tab:qualitative}, the requirements for construction of the unstructured spline bases are summarized in \cref{fig:USconstraints} as pre-processing conditions that have to be satisfied for each unstructured spline construction in the process depicted in \cref{fig:IGA_Setup}. The degree and regularity conditions (cf. \emph{i,ii} in \cref{tab:qualitative}) must be satisfied for each construction, e.g. by performing projections on suitable spline spaces or by knot insertion routines. Furthermore, the geometric or topological limitations (cf. \emph{iii} in \cref{tab:qualitative}) impose additional constraints that the geometry must satisfy. \\

\section{Quantitative comparison}\label{sec:quantitative}

In this section a quantitative comparison between the methods provided in \cref{sec:qualitative} is provided. In addition, variational coupling methods are compared if applicable. The quantitative comparison is composed of various benchmark problems, each providing a different conclusion with respect to the methods considered:
\begin{description}
    \item[Biharmonic problem on a planar domain (\cref{subsec:quantitative_biharmonic})] The first example entails solving the biharmonic problem on a planar domain. The goal of this example is to assess the convergence properties of all considered unstructured spline constructions, hence the problem will be solved on a simple analysis-suitable geometry without EVs on the boundary, such that every method from \cref{sec:qualitative} can be applied and compared to the manufactured solution.
    \item[Linear Kirchhoff--Love shell analysis on a surface (\cref{subsec:quantitative_linearShell})] The second example entails solving the Kirchhoff--Love shell equation on curved domains. The goal of this example is to demonstrate the performance of the unstructured spline construction for simple shell problems. Therefore, comparison will be made to single-patch results and penalty coupling from \cite{Herrema2019}.
    \item[Spectral analysis on a planar domain (\cref{subsec:quantitative_eigenvalue})] In the third example, spectral analysis of a plate equation is performed. The goal of this example is to assess the spectral properties of the unstructured spline methods compared to a variational approach and a single patch, since the spectral properties of highly continuous bases have been demonstrated to be superior over non-smooth bases \cite{Cottrell2007}.
    \item[Modal analysis of a complex geometry (\cref{subsec:quantitative_modalCar})] In the fourth example, a modal analysis is performed on a complex geometry extracted from a quad-mesh. The goal of this example is to demonstrate the applicability and performance of the unstructured spline methods on a large-scale, more complicated geometry.
    \item[Stress analysis in a curved shell (\cref{subsec:quantitative_stress})] Lastly, the fifth example involves the analysis of stress fields in shells. The goal of this example is to assess the performance of unstructured spline constructions and a penalty method when it comes to stress reconstruction in shells. For the Kirchhoff--Love shell, the stresses are obtained by taking gradients of the deformed geometry, hence of the solution. This means that for $C^1$ bases, stresses are $C^0$. This might be unfavourable in engineering applications where local stress fields are of importance, e.g.~fatigue analysis.
\end{description}

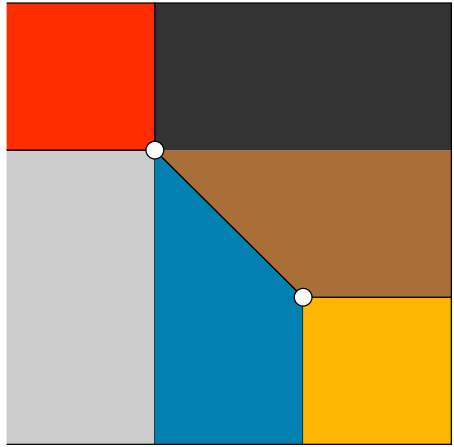
\begin{figure}
    \centering
    \resizebox{0.45\linewidth}{!}
    {
	    \begin{tikzpicture}[scale=5]
	    \filldraw[draw=black, fill=black!20] (0,0) -- (1/3,0) -- (1/3,2/3) -- (0,2/3);
	    \filldraw[draw=black, fill=col1] (1/3,0) -- (2/3,0) -- (2/3,1/3) -- (1/3,2/3);
	    \filldraw[draw=black, fill=col2] (2/3,0) -- (1,0) -- (1,1/3) -- (2/3,1/3);
	    \filldraw[draw=black, fill=col3] (0,2/3) -- (1/3,2/3) -- (1/3,1) -- (0,1);
	    \filldraw[draw=black, fill=col4] (1/3,2/3) -- (2/3,1/3) -- (1,1/3) -- (1,2/3);
	    \filldraw[draw=black, fill=black!80] (1/3,2/3) -- (1/3,1) -- (1,1) -- (1,2/3);
	    \filldraw[draw=black, fill=white] (1/3,2/3) circle[radius=0.02];
	    \filldraw[draw=black, fill=white] (2/3,1/3) circle[radius=0.02];
	    \end{tikzpicture}
	}
\caption{Multi-patch decomposition of a simple domain into six patches. The domain has two EVs in the interior (valence 3 and 5) and no boundary EVs.}
\label{fig:domain}
\end{figure}

In all examples except the complex geometry in \cref{subsec:quantitative_modalCar}, the domain decomposition from \cref{fig:domain} is used to decompose a simple domain into a domain with extraordinary vertices in the interior. Domains with EVs on the boundary are left out of scope, since the D-patch construction would change the outer boundaries of the domain, hence the comparison would involve a significantly different geometry. Since different methods have different constraints on the degree and regularity of the basis, different combinations of the degree $p$ and regularity $r$ are tested throughout the benchmark problems. In \cref{tab:comparison} the combinations of $p$ and $r$ and the methods that are compared for these bases are provided. For the biharmonic problem and the spectral analysis (\cref{subsec:quantitative_biharmonic,subsec:quantitative_eigenvalue}) Nitsche's method is used for comparison, see \cite{Weinmuller2022} for more details. When solving the Kirchhoff--Love shell equations, the penalty method is used for comparison, see \cite{Herrema2019} for more details. In all examples, Dirichlet boundary conditions are applied at the control points and clamped boundary conditions are applied weakly as in \cite{Herrema2019}. All results are obtained using the Geometry + Simulation modules \cite{Juttler2014,m2020} and will be published in a separate publication.\\

\begin{table}
    \centering
    \caption{Degree $p$ and regularity $r$ constraints for each considered method from \cref{sec:qualitative}, see \cref{tab:qualitative}.}
    \label{tab:comparison}
    \begin{tabular}{lccc}
        \toprule
                        & $p=2$, $r=1$  & $p=3$, $r=1$  & $p=3$, $r=2$  \\
        \midrule
        D-patch         & $\bigstar$    & $\bigstar$    & $\bigstar$    \\
        Almost-$C^1$    & $\bigstar$    &               &               \\
        Approx.~$C^1$   &               & $\bigstar$    & $\bigstar$    \\
        AS-$G^1$        &               & $\bigstar$    &               \\
        Nitsche/Penalty & $\bigstar$    & $\bigstar$    & $\bigstar$    \\
        \bottomrule
    \end{tabular}
\end{table}

As discussed in \cref{sec:qualitative}, the D-patch and Almost-$C^1$ constructions involve a pre-smoothing of the geometry. In case of mesh convergence results, refinements can be performed in different ways. On the one hand, the original geometry can be refined and a new construction with a new geometry approximation can be performed. On the other hand, the geometry resulting from the construction in the first refinement level can be refined in a nested way, such that the geometry does not change after the first mesh. In the quantitative comparison, all refinements are performed in a nested way, unless specified otherwise.

\subsection{Biharmonic equation on a planar domain}\label{subsec:quantitative_biharmonic}
The first benchmark entails the biharmonic equation on a planar domain. The purpose of this example is to assess the convergence properties of the unstructured spline methods described in \cref{sec:qualitative}. We mainly follow the structure of \cite{Weinmuller2022}. The biharmonic equation is solved on a unit square $\Omega = [0,1]^2$ with the patch segmentation from \cref{fig:domain}. The biharmonic equation is defined by
\begin{equation}\label{eq:biharmonic}
    \Delta^2\varphi = f.
\end{equation}
In the present example, convergence is analysed with respect to a manufactured solution
\begin{equation}
    \tilde{\varphi}(x_1,x_2) = (\cos (4 \pi x_1) - 1) (\cos (4 \pi x_2) - 1),
\end{equation}
such that the right-hand-side function becomes:
\begin{equation}
    f(x_1,x_2) = 256\pi^4(4\cos(4\pi x)\cos(4\pi y) - \cos(4\pi x) - \cos(4\pi y))
\end{equation}
Furthermore, on all boundaries of the domain, the manufactured solution and its derivatives are imposed as Dirichlet and Neumann boundary conditions, respectively:
\begin{equation}\label{eq:biharmonic_BCs}
\begin{aligned}
    \begin{rcases}
    \varphi &=\tilde{\varphi}(x_1,x_2)\\
    \partial_{\vb{n}}\varphi &= \partial_{\vb{n}}\tilde{\varphi}
    \end{rcases}
    \text{on }\Gamma,
\end{aligned}
\end{equation}
where $\Gamma = \partial \Omega$, $\vb{n}$ is the unit outward normal vector on $\Gamma$. The biharmonic equation from  \cref{eq:biharmonic} with boundary conditions \cref{eq:biharmonic_BCs} can be discretized by obtaining the weak formulation, see \cite{Weinmuller2022}, inserting \cref{eq:biharmonic_BCs} and by defining an approximation of the solution $\varphi$ as $\varphi_h$. Furthermore, a weak coupling can be established through Nitsche's method. For the mathematical details behind the discretization of the biharmonic equation and optionally adding Nitsche interface coupling terms, we refer to \cite{Weinmuller2021,Weinmuller2022}. For the D-Patch and Almost $C^1$ constructions, the geometry is smoothed upon construction. The geometry used for evaluation of the weak formulation is constructed by using an $L_2$-projection of the geometry from the coarsest space which is projected onto the smooth basis of each refinement level. For the D-Patch, the non-negative smoothness matrix for vertex smoothing is used. Although this matrix produces non-nested meshes, it provides the highest rates of convergence. Furthermore, the factor $\beta$ (cf. \cite[sec. 5.1]{Toshniwal2017}) is chosen as $\beta=0.4$ as used by \cite{Toshniwal2017}, or $\beta=1.2$, and halved in each refinement level.\\

To evaluate the unstructured spline constructions from \cref{sec:qualitative}, the numerical approximation $\varphi_h$ is compared to the manufactured solution $\tilde{\varphi}$ in the $L^2$-, $H^1$- and $H^2$-norms on the multi-patch segmentation from \cref{fig:domain}. The bi-linear segmentation is refined and degree elevated until the desired degree $p$ and regularity $r$ from \cref{tab:comparison} are obtained. In addition, a Nitsche coupling of the patches is employed for comparison. \\

The results for the comparison are presented in \cref{fig:biharmonic_results}. For degree $p=2$ and regularity $r=1$ the Almost-$C^1$, D-patch and Nitsche coupling methods are compared. As expected, the results show consistency between the Almost-$C^1$, D-patch and Nitsche's method with expected convergence. The results also show a slight dependency on the factor $\beta$ for the D-Patch. For degree $p=3$ and regularity $r=1$, the Approx.~$C^1$, AS-$G^1$, D-patch and Nitsche's method can be compared. The results of the Approx.~$C^1$ and AS-$G^1$ are exactly the same, since the original geometry is analysis-suitable and contains only bi-linear patches. Then applying the Approx.~$C^1$ to an analysis-suitable geometry with regularity $p-2$, the approximate gluing data becomes exact, hence the same as in the AS-$G^1$ construction. The D-patch in this case shows better convergence of the $L_2$, $H_1$ and $H_2$ errors for $\beta=1.2$ than for $\beta=0.4$. Though, for both choices of $\beta$, the convergence is sub-optimal, as was observed in the work by \cite{Casquero2020}. Furthermore, the $L_2$-norm increases at the last point of the D-Patch results, due to ill-conditioning of the system of equations. Lastly, for degree $p=3$ and regularity $r=2$ the Approx.~$C^1$, D-patch and Nitsche's method are compared. The observations are as for the $p=3$, $r=1$ case.\\

Overall, the results show expected convergence behaviour for all considered spline constructions compared to theoretical results and compared to a Nitsche coupling method. However, the D-Patch method does not converge for very fine meshes, due to ill-conditioning of the system matrix.

\input{Figures/Plots/Biharmonic.tex}

\subsection{Linear Kirchhoff--Love shell analysis on a surface}\label{subsec:quantitative_linearShell}
We solve the linear Kirchhoff--Love shell equations on two geometries to demonstrate the convergence behaviour of the methods on curved surfaces. To this end, two benchmark examples are considered. Firstly, a hyperbolic paraboloid surface is constructed with shape, inspired by \cite{Farahat2023}:
\begin{equation}
    \vb{r}(\xi_1,\xi_2) = \qty[\xi_1 \quad \xi_2 \quad \xi_1^2-\xi_2^2]
\end{equation}
The left-side of the hyperbolic paraboloid is clamped ($\vb{u}=\vb{0}$) and the other sides are free. Furthermore, a distributed load with magnitude $8000t$ is applied with $t$ the thickness, see \cref{fig:hyperboloid_setup}. Secondly, an elliptic paraboloid shaped-domain is modelled, with equation
\begin{equation}
    \vb{r}(\xi_1,\xi_2) = \qty[\xi_1 \quad \xi_2 \quad 1-2\qty(\xi_1^2+\xi_2^2)]
\end{equation}
For this shape, a point load with magnitude $10^8t$ is applied in the middle of the domain. The corners of the domain are only fixed in vertical $z$ direction to allow sliding in the $xy$-plane. One corner is fixed in all directions to create a well-posed problem. For both hyperbolic paraboloid (\cref{fig:hyperboloid_setup}) and elliptic paraboloid (\cref{fig:paraboloid_setup}) the multi-patch segmentation from \cref{fig:domain} is used. In both cases, the shells are modelled with a thickness of $t=0.01\:[\text{mm}]$ and with a Saint-Venant Kirchhoff material with Young's modulus $E=200\:[\text{GPa}]$ and Poisson's ratio $\nu=0.3\:[\text{-}]$. The refinement procedure as described in \cref{subsec:quantitative_biharmonic} is used for the D-Patch and Almost-$C^1$ constructions.\\

\begin{figure}
\centering
\begin{minipage}[t]{0.45\linewidth}
\centering
\resizebox{\linewidth}{!}{\begin{tikzpicture}[]
\begin{axis}
[
width = 0.8\linewidth,
height = 0.8\linewidth,
grid=none,
axis lines=left,
zmin = -0.5,
zmax = 0.5,
enlargelimits,
xlabel={$x$},
ylabel={$y$},
zlabel={$z$},
no markers
]
    \addplot3 [domain=-0.5:0.5, samples y=1] ({x},{0.5},{x^2-0.5^2});
\addplot3[colormap/blackwhite,line width=0pt,opacity=1,surf,domain=-0.5:0.5, y domain= -0.5:0.5,shader=interp,opacity=0.8,samples=20, samples y = 20]{x^2-y^2};
    \addplot3 [domain=-0.5:0.5, samples y=1] ({x},{-0.5},{x^2-0.5^2});
    \addplot3 [ultra thick,domain=-0.5:0.5, samples y=1] ({-0.5},{x},{0.5^2-x^2});
    \addplot3 [domain=-0.5:0.5, samples y=1] ({0.5},{x},{0.5^2-x^2});

\fill (axis cs: 0.5,0,0.25) node[above left] {$A$} circle [radius=10];
\foreach \x in {0,-0.05,...,-0.25}
{
    \foreach \y in {-0.25,-0.30,...,-0.50}
    {
        \edef\z{-\x^2+\y^2}
        \edef
        \temp{
                    \noexpand
                    \draw [latex-] (axis cs:\x,\y,\z) -- (axis cs:\x,\y,\z+0.2);
                }
        \temp
    }
}
\node[left] (q) at (axis cs: -0.25,-0.25,0) {$p_{x^3}$};
\end{axis}
\end{tikzpicture}}
\caption{Hyperbolic paraboloid shell geometry with coordinates $\vb{r}(\xi_1,\xi_2) =  \qty[\xi_1 \quad \xi_2 \quad \xi_1^2-\xi_2^2]$, $\xi_1,\xi_2\in[-1/2.1/2]$. The left-edge of the hyperbolic paraboloid is clamped, i.e.~the displacements and rotations are zero ($\vb{u}=\vb{0}$ and $\pdv{u_z}{x}=0$).}
\label{fig:hyperboloid_setup}
\end{minipage}
\hfill
\begin{minipage}[t]{0.45\linewidth}
\centering
\resizebox{\linewidth}{!}{\begin{tikzpicture}[]
\begin{axis}
[
width = 0.8\linewidth,
height = 0.8\linewidth,
grid=none,
axis lines=left,
zmin = 0.0,
zmax = 2.0,
enlargelimits,
xlabel={$x$},
ylabel={$y$},
zlabel={$z$},
no markers
]
    \addplot3 [domain=-0.5:0.5, samples y=1] ({x},{0.5},{1-2*(x^2+0.5^2)});
    \addplot3[colormap/blackwhite,line width=0pt,opacity=1,surf,domain=-0.5:0.5, y domain= -0.5:0.5,shader=interp,opacity=0.8,samples=20, samples y = 20]{1-2*(x^2+y^2)};
    \addplot3 [domain=-0.5:0.5, samples y=1] ({x},{-0.5},{1-2*(x^2+0.5^2)});
    \addplot3 [domain=-0.5:0.5, samples y=1] ({0.5},{x},{1-2*(0.5^2+x^2)});

\fill (axis cs: 0.0,0.0,1.0) node[above left] {$A$} circle [radius=10];
\draw[latex-] (axis cs: 0.0,0.0,1.0) -- (axis cs: 0.0,0.0,1.5) node[above]{$P$};
\end{axis}
\end{tikzpicture}}
\caption{Elliptic paraboloid shell geometry with coordinates $\vb{r}(\xi_1,\xi_2) =  \qty[\xi_1 \quad \xi_2 \quad 1-2\qty(\xi_1^2+\xi_2^2)]$, $\xi_1,\xi_2\in[-1/2.1/2]$. On the corners of the domain, the vertical displacements are set to zero $u_z=0$ and one corner is fixed in-plane as well. Furthermore, a point load with magnitude $P=10^8t$ is applied in the middle of the geometry.}
\label{fig:paraboloid_setup}
\end{minipage}
\end{figure}

The results of both analyses are given in \cref{fig:linearShell_hyperboloid,fig:linearShell_paraboloid}. Here, different unstructured spline constructions are tested on patch-bases with different degrees and regularities, as reported in \cref{tab:comparison}. For each combination of degree $p$ and regularity $r$, the energy norm $W_{\text{int}}^h=\frac{1}{2}\vb{u}_h^\top K_h \vb{u}_h$ is plotted against the number of degrees of freedom, with $\vb{u}_h$ the discrete displacement vector and $K_h$ the discrete linear stiffness matrix. From the results in \cref{fig:linearShell_hyperboloid,fig:linearShell_paraboloid}, a few observations can be made. Firstly, the Approx.~$C^1$ and AS-$G^1$ methods show slow convergence on the hyperbolic paraboloid geometry, while the convergence on the elliptic paraboloid geometry is similar to the single-patch convergence. The slow convergence for the hyperbolic paraboloid shell is also observed in \cite{Farahat2023}. Since the results of the same constructions on the elliptic paraboloid geometries do not show slower convergence, the slow convergence is hypothetically a result of the double curvature with different signs of the shell.
Secondly, the D-Patch and Approx.~$C^1$ show comparable convergence to the penalty method on both geometries, which is slightly slower than the convergence of the single-patch results. This is explained by the fact that the degrees of freedom are more optimally allocated for the single-patch parameterization. Lastly, the results obtained by the penalty method for different penalty parameters $\alpha$ show convergence with a rate similar to the D-Patch and Almost-$C^1$ methods for penalty parameters $\alpha\in\{1,10\}$. For $\alpha=100$ the penalty method is still converging to the same solution, but convergence starts after a few refinement steps.

\input{Figures/Plots/LinearShell_hyperboloid.tex}
\input{Figures/Plots/LinearShell_paraboloid.tex}

\subsection{Spectral analysis on a planar domain}\label{subsec:quantitative_eigenvalue}
In this example, the spectral properties of the unstructured spline constructions on a multi-patch domain are considered. From \cite{Cottrell2006} it is known that isogeometric analysis has the advantage over $C^0$ Finite Element Analysis with respect to spectra for eigenvalue problems. Smooth isogeometric discretization provide converging spectra with spline degree $p$, whereas the spectra obtained by $C^0$ FEA diverge with $p$ and typically have optical branches. Similarly, when patches with $C^0$ continuity are considered, optical branches are introduced and the accuracy of the spectral approximation decreases \cite{Nguyen2022a}. In this benchmark problem, we compare the basis constructions from \cref{tab:comparison} on \cref{fig:domain} on their spectral properties. For Nitsche's method, we use different values for the coupling parameter to assess its influence on the spectrum.\\

For the problem at hand, we consider a unit-square domain with parametric lay-out from \cref{fig:domain} for simplicity. We consider modal analysis using the plate equation. The stiffness operator of the free vibration plate equation is similar to the biharmonic equation from \cref{eq:biharmonic}, and the inertia is included on the right-hand-side:
\begin{equation}
    D\Delta^2 w = -\rho t \pdv[2]{w}{\tau}
\end{equation}
Assuming that $w(x,y,\tau)$ is harmonic, i.e.~$w(x,y,\tau)=\hat{w}(x,y)\exp{i\omega \tau}$ with $\omega$ a frequency, the equation simplifies to
\begin{equation}\label{eq:freevibration}
    D\Delta^2 \hat{w} = \omega^2 \pdv[2]{\hat{w}}{\tau}.
\end{equation}
Here, $D=Et^3/(12(1-\nu^2))$ is the flexural rigidity of the plate with $E=10^5\:[\text{Pa}]$ the Young's modulus of the plate, $t=10^{-2}\:[\text{m}]$ the thickness and $\nu=0.2\:[\text{-}]$ the Poisson's ratio. Furthermore, $\rho=10^5\:[\text{kg}]$ is the material density. \Cref{eq:freevibration} is a generalized eigenvalue problem with eigenpairs $(\omega_i,v_i)$ where $\omega_i$ is the $\text{i}^{\text{th}}$ eigenfrequency and $v_i$ the $\text{i}^{\text{th}}$ mode shape. The mode shape for a simply supported unit plate with $n\times m$ half-waves is given by
\begin{equation}
    v_{nm}(x,y) = \sin\qty(n\pi x)\sin\qty(m\pi y)
\end{equation}
with corresponding eigenfrequency
\begin{equation}
    \omega_{nm} = (n^2+m^2)\pi^2\sqrt{\frac{D}{\rho t}}.
\end{equation}
In addition, the numerical solution to \cref{eq:freevibration} is obtained by solving the following generalized eigenvalue problem
\begin{equation}
    D \int_\Omega\Delta w \Delta \varphi\dd{\Omega} = \omega^2\rho t \int_\Omega u\varphi\dd{\Omega}
\end{equation}
With $\varphi$ a test function, see \cref{subsec:quantitative_biharmonic}. In further representation of the solutions, we employ the index $i$ such that $\omega_i<\omega_{i+1}$ and we use the subscript $h$ for numerically obtained solutions.\\

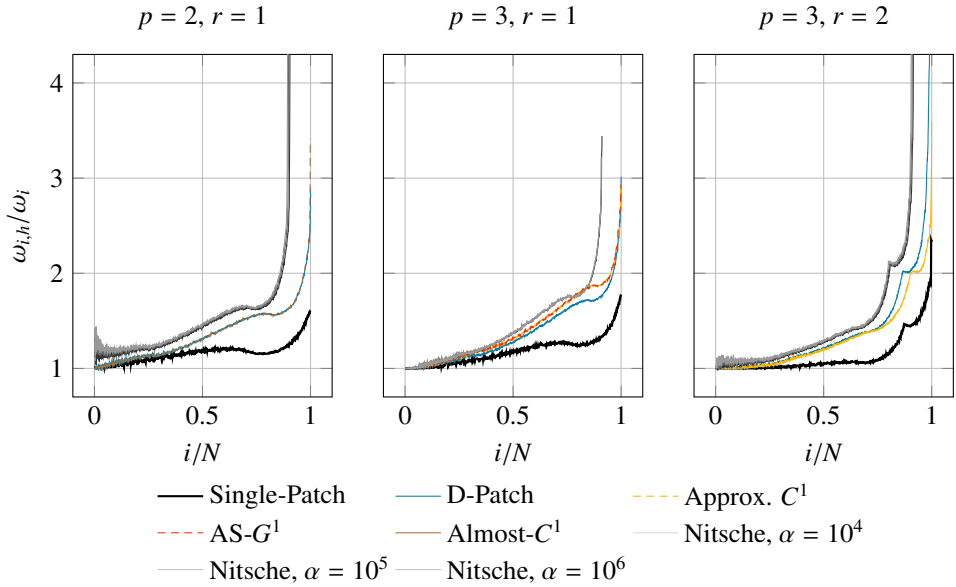
\begin{figure}
    \centering
    \begin{tikzpicture}
            \begin{groupplot}
            [
                group style={
                group size=3 by 1,
                ylabels at=edge left,
                y descriptions at=edge left,
                xlabels at=edge bottom,
                x descriptions at=edge bottom,
                vertical sep=0.05\textheight,
                horizontal sep=0.05\linewidth,
                },
                xlabel={$i/N$},
                ylabel={$\omega_{i,h}/\omega_i$},
                legend pos = north west,
                ymin = 1,
                ymax = 4,
                restrict y to domain = {0:5},
                width=0.37\linewidth,
                height=0.31\textheight,
                grid=major,
                legend style={fill=white, fill opacity=0.6, draw opacity=1,draw=none,text opacity=1},
                enlargelimits
            ]
            \nextgroupplot[title={$p=2$, $r=1$}]
            \addplot+[style=singlepatch,no markers,black,fill opacity=1.0] table[header=true,x index = {3},y index = {2}, col sep = comma]{Data/Biharmonic_eigenvalue/Single_p2_s1.csv};
            \addplot+[style=dpatch,no markers,fill opacity=1.0] table[header=true,x index = {3},y index = {2}, col sep = comma]{Data/Biharmonic_eigenvalue/DPatch_p2_s1.csv};
            \addplot+[style=almostC1,densely dashed,no markers,fill opacity=1.0] table[header=true,x index = {3},y index = {2}, col sep = comma]{Data/Biharmonic_eigenvalue/AlmostC1_p2_s1.csv};
            \addplot+[style=weak1,no markers,fill opacity=1.0] table[header=true,x index = {3},y index = {2}, col sep = comma]{Data/Biharmonic_eigenvalue/Nitsche_1e4_p2_s1.csv};
            \addplot+[style=weak2,no markers,fill opacity=1.0] table[header=true,x index = {3},y index = {2}, col sep = comma]{Data/Biharmonic_eigenvalue/Nitsche_1e5_p2_s1.csv};
            \addplot+[style=weak3,no markers,fill opacity=1.0] table[header=true,x index = {3},y index = {2}, col sep = comma]{Data/Biharmonic_eigenvalue/Nitsche_1e6_p2_s1.csv};

            \nextgroupplot[title={$p=3$, $r=1$}]
            \addplot+[style=singlepatch,no markers,black,fill opacity=1.0] table[header=true,x index = {3},y index = {2}, col sep = comma]{Data/Biharmonic_eigenvalue/Single_p3_s1.csv};
            \addplot+[style=dpatch,no markers,fill opacity=1.0] table[header=true,x index = {3},y index = {2}, col sep = comma]{Data/Biharmonic_eigenvalue/DPatch_p3_s1.csv};
            \addplot+[style=approxC1,no markers,fill opacity=1.0] table[header=true,x index = {3},y index = {2}, col sep = comma]{Data/Biharmonic_eigenvalue/ApproxC1_p3_s1.csv};
            \addplot+[style=exactC1,densely dashed,no markers,fill opacity=1.0] table[header=true,x index = {3},y index = {2}, col sep = comma]{Data/Biharmonic_eigenvalue/ASG1_p3_s1.csv};
            \addplot+[style=weak1,no markers,fill opacity=1.0] table[header=true,x index = {3},y index = {2}, col sep = comma]{Data/Biharmonic_eigenvalue/Nitsche_1e4_p3_s1.csv};
            \addplot+[style=weak2,no markers,fill opacity=1.0] table[header=true,x index = {3},y index = {2}, col sep = comma]{Data/Biharmonic_eigenvalue/Nitsche_1e5_p3_s1.csv};
            \addplot+[style=weak3,no markers,fill opacity=1.0] table[header=true,x index = {3},y index = {2}, col sep = comma]{Data/Biharmonic_eigenvalue/Nitsche_1e6_p3_s1.csv};

            \nextgroupplot[title={$p=3$, $r=2$},legend to name={CommonLegend},legend columns=3]
            \addlegendimage{style=singlepatch   ,mark=none}\addlegendentry[]{Single-Patch}
            \addlegendimage{style=dpatch        ,mark=none}\addlegendentry[]{D-Patch}
            \addlegendimage{style=approxC1,densely dashed      ,mark=none}\addlegendentry[]{Approx. $C^1$}
            \addlegendimage{style=exactC1,densely dashed       ,mark=none}\addlegendentry[]{AS-$G^1$}
            \addlegendimage{style=almostC1      ,mark=none}\addlegendentry[]{Almost-$C^1$}
            \addlegendimage{style=weak1         ,mark=none,opacity=0.2}\addlegendentry[]{Nitsche, $\alpha=10^4$}
            \addlegendimage{style=weak2         ,mark=none,opacity=0.4}\addlegendentry[]{Nitsche, $\alpha=10^5$}
            \addlegendimage{style=weak3         ,mark=none,opacity=0.6}\addlegendentry[]{Nitsche, $\alpha=10^6$}

            \addplot+[style=singlepatch,no markers,black,fill opacity=1.0] table[header=true,x index = {3},y index = {2}, col sep = comma]{Data/Biharmonic_eigenvalue/Single_p3_s2.csv};
            \addplot+[style=dpatch,no markers,fill opacity=1.0] table[header=true,x index = {3},y index = {2}, col sep = comma]{Data/Biharmonic_eigenvalue/DPatch_p3_s2.csv};
            \addplot+[style=approxC1,no markers,fill opacity=1.0] table[header=true,x index = {3},y index = {2}, col sep = comma]{Data/Biharmonic_eigenvalue/ApproxC1_p3_s2.csv};
            \addplot+[style=weak1,no markers,fill opacity=1.0] table[header=true,x index = {3},y index = {2}, col sep = comma]{Data/Biharmonic_eigenvalue/Nitsche_1e4_p3_s2.csv};
            \addplot+[style=weak2,no markers,fill opacity=1.0] table[header=true,x index = {3},y index = {2}, col sep = comma]{Data/Biharmonic_eigenvalue/Nitsche_1e5_p3_s2.csv};
            \addplot+[style=weak3,no markers,fill opacity=1.0] table[header=true,x index = {3},y index = {2}, col sep = comma]{Data/Biharmonic_eigenvalue/Nitsche_1e6_p3_s2.csv};
        \end{groupplot}

    \path (group c1r1.south west) -- node[below=25pt]{\ref*{CommonLegend}} (group c3r1.south east);

    \end{tikzpicture}
    \caption{Eigenvalue spectra for the biharmonic eigenvalue problem on the domain from \cref{fig:domain}. The horizontal axes depict the eigenvalue index $i$ over the total number of eigenvalues $N$. The vertical axes represent the numerical eigenvalue $\omega_{i,h}$ over the analytical eigenvalue $\omega_i$, both with index $i$. The results are plotted for different combinations of the degree $p$ and regularity $r$ of the basis. The results for a Nitsche method are given for different penalty parameters $\alpha$.}
    \label{fig:spectrum}
\end{figure}

\Cref{fig:spectrum} presents the spectra for different degrees, regularities and for different methods. Here, the vertical axis represents the ratio of the numerically obtained eigenfrequency over the analytical eigenfrequency with index $i$, thus $\omega_{h,i}/\omega_i$. The horizontal axis represents the fraction of the eigenfrequency index $i$ over the total number of eigenmodes. The total number of eigenmodes is equal to the number of degrees of freedom in the system. The results are presented for the degrees and regularities as in \cref{tab:comparison}.\\

Firstly, the $p=2$, $r=1$ plot shows that Nitsche's method oscillates for all considered values of the penalty parameter. Furthermore, in the part where it is not oscillating, the ratio $\omega_{i.h}/\omega_i$ is higher than for the D-patch and Almost-$C^1$ method. Additionally, the D-patch and Almost-$C^1$ methods show a significant difference with respect to the single patch result, which is due to the non-Cartesian multi-patch segmentation of \cref{fig:domain} and the fact that the analytical mode shapes are Cartesian. For the $p=3$, $r=1$ and $p=3$, $r=2$ bases similar conclusions can be drawn. Although for the $p=3$, $r=1$ case the Approx.~$C^1$ method seems worse than the D-patch method, the opposite is true for $p=3$, $r=2$. Hence, it can be concluded that no method outperforms another, but that all unstructured spline constructions perform better than Nitsche's method.\\

\subsection{Modal analysis of a complex geometry}\label{subsec:quantitative_modalCar}
\begin{figure}
\centering
\begin{subfigure}{\linewidth}
\centering
\includegraphics[height=0.25\textheight]{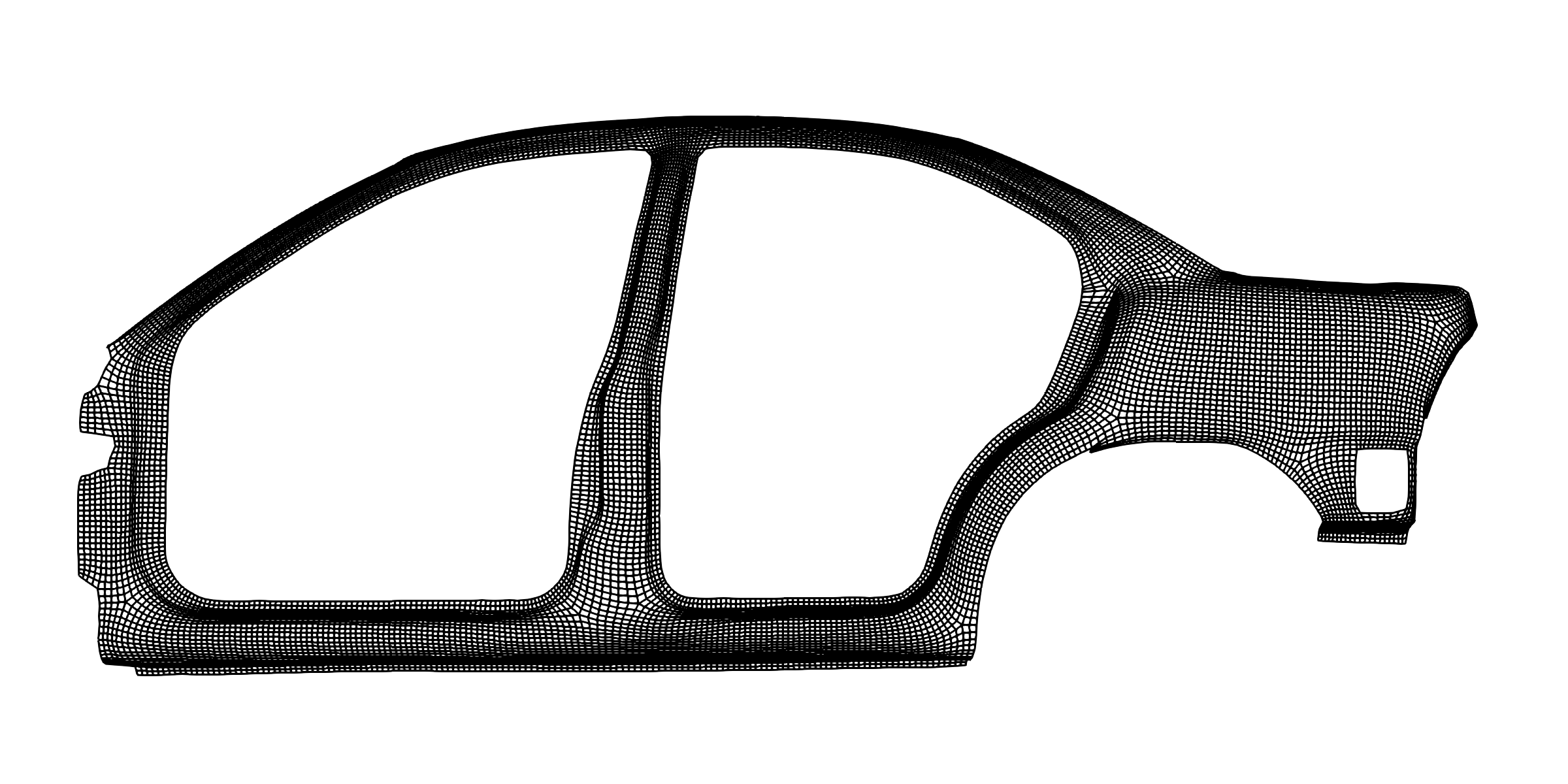}
\caption{Original quad mesh with 15895 vertices, 31086 edges and 62172 faces.}
\label{fig:car_mesh}
\end{subfigure}

\begin{subfigure}{\linewidth}
\centering
\includegraphics[height=0.25\textheight]{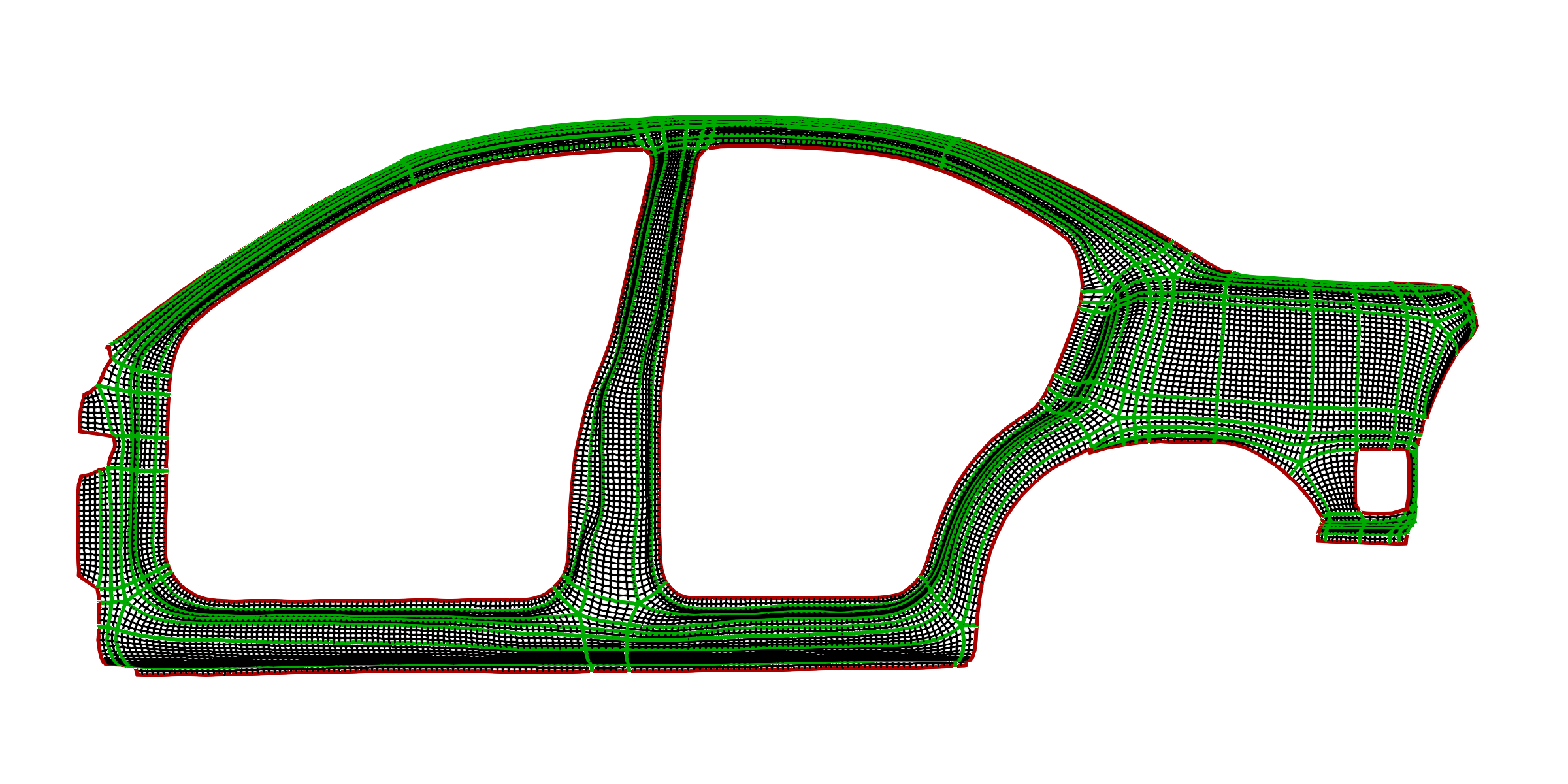}
\caption{Interface (green) and boundary (red) curves.}
\label{fig:car_patches}
\end{subfigure}

\begin{subfigure}{\linewidth}
\centering
\includegraphics[height=0.25\textheight]{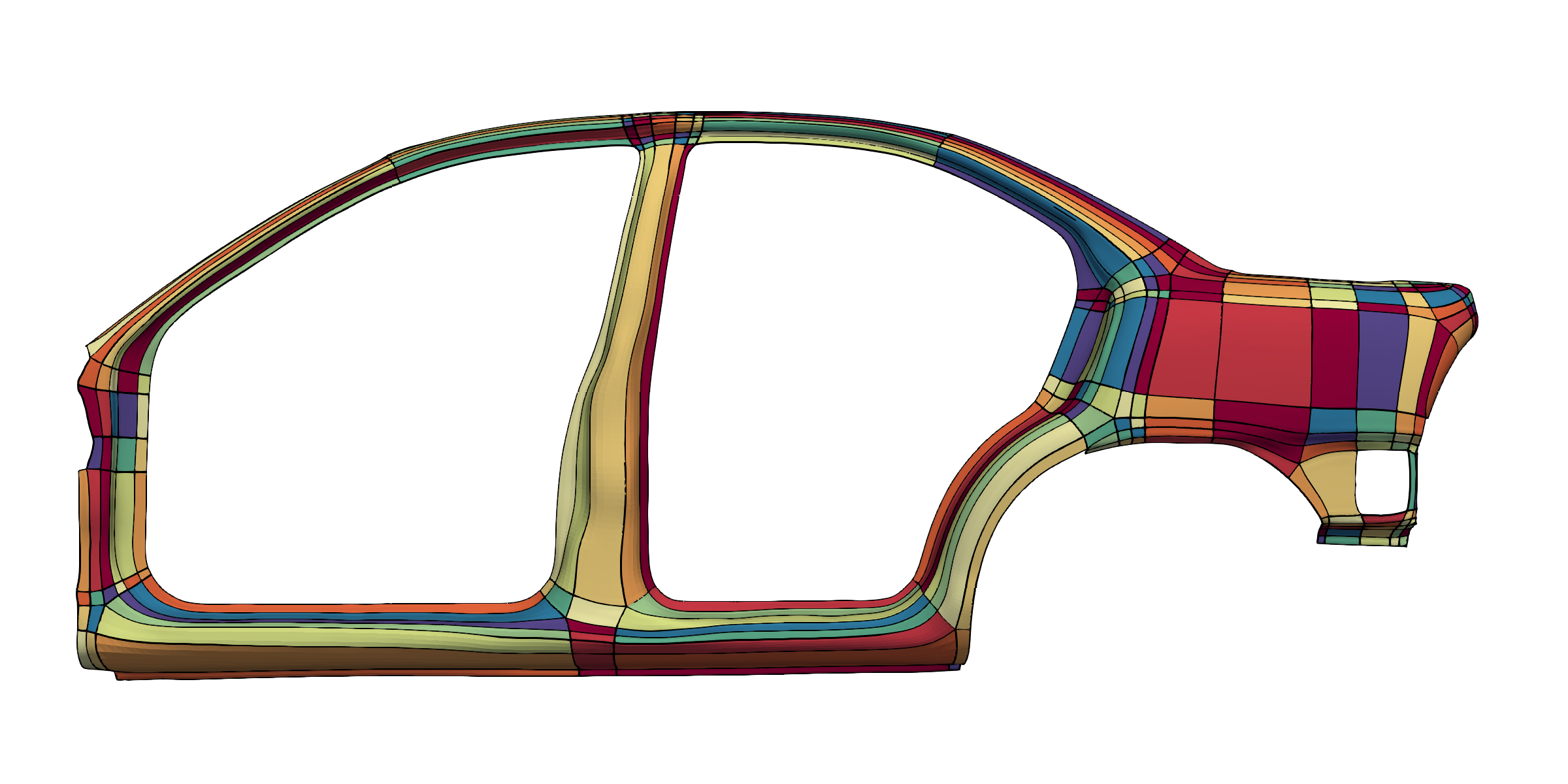}
\caption{Final multi-patch segmentation with 307 patches.}
\label{fig:car_almostC1}
\end{subfigure}
\caption{Geometry of the side panel of a car. The original mesh (a) is traced with the procedure from \cref{fig:Illustration_patches_from_mesh}, yielding a set of boundary and interface curves (b). From these curves, the multi-patch segmentation (c) for isogeometric analysis is constructed following \cref{fig:Illustration_patches_from_mesh3}.}
\label{fig:car_side}
\end{figure}

The next example for the quantitative analysis in this paper involves the modal analysis on a larger-scale complex geometry, depicted in \cref{fig:car_mesh}. The goal of this example is to show the usability of the considered constructions on an off-the-shelf industrial geometry. The geometry is represented as a mesh consisting of 15895 vertices, 31086 edges and 62172 faces. This geometry is converted to bi-linear patches using the procedure discussed in \cref{fig:Illustration_patches_from_mesh} in \cref{sec:qualitative}. The interface and boundary curves of the patches are given in \cref{fig:car_patches} and the final multi-patch object is given in \cref{fig:car_almostC1}. The latter has $3$ EVs of valence $3$, $10$ EVs of valence $5$ and $16$ bEVs. Moreover, the material parameters specified for a steel material. That is, the density of the material is $\rho=7850\cdot 10^{-6}\:[\text{tonnes}/\text{mm}^3]$, the shell thickness is $t=10\:[\text{mm}]$, the Young's modulus is $E=210\cdot 10^3\:[\text{MPa}]$ and the Poisson's ratio is $\nu=0.3\:[-]$. All the sides of the geometry are kept free, meaning that the modal analysis results will consist of six modes with zero eigenfrequencies: the rigid body modes. In the sequel, we list the results for deformation modes only.\\

After the creation of the linear multi-patch object, $h$-, $p$- and $k$-refinement steps can be performed to construct a multi-basis corresponding to the patch lay-out on which unstructured splines can be constructed. For the Almost-$C^1$ and D-Patch constructions, the bases are constructed by refining and elevating the initial linear basis up to the desired degree and regularity, after which the the Almost-$C^1$ and D-Patch basis and geometry are computed. An Almost-$C^1$ geometry is provided in \cref{fig:car_almostC1}.\\ 

The AS-$G^1$ construction requires an analysis-suitable geometry, which can be constructed following \cite{Farahat2023}, is based on the planar construction developed in \cite{Kapl2018}. However, the geometry from \cref{fig:car_patches} is only $C^0$-smooth due to the original linear mesh it is constructed from. An algorithm to automatically pre-process the geometry to obtain an analysis-suitable $G^1$ surface is not yet developed. The algorithm from \cite{Farahat2023} requires AS-$G^1$ gluing data, which cannot be prescribed directly on a $C^0$ surface. If the surface is not pre-processed to be AS-$G^1$, no suitable gluing data can be found and the basis construction is not applicable. Although the Approx. $C^1$ construction does not require an analysis-suitable re-parameterization, it does require $G^1$ smoothness at the interfaces. If this condition is not satisfied there exists no $C^1$ construction that can be approximated by this method. For both methods, the required pre-processing efforts are non-trivial or not demonstrated on industrial geometries, and therefore left out of the scope of this paper\footnote{In the case of a different starting point for this benchmark, such as a smooth mesh composed of higher-order quadrilateral elements, e.g. derived from a subdivision surface, instead of a bi-linear mesh, the pre-processing efforts required for the AS-$G^1$ and Approx.~$C^1$ will be different.}.\\

Furthermore, penalty methods have been used in the context of modal analysis on a 27 patch composite wind-turbine blade in \cite{Herrema2019}, where the variation of the element size of interface elements seems rather small. In the present paper, an attempt was made to apply the penalty method on the geometry in \cref{fig:car_almostC1}, but unidentifiable vibration modes were obtained, possibly because of the large variation of element lengths across the interfaces of the domain, challenging the determination of a suitable penalty parameter $\alpha$.\\

\begin{figure}
	\centering
	\begin{subfigure}{\linewidth}
		\centering
		\begin{minipage}{0.45\linewidth}
			\centering
				\includegraphics[width=\linewidth,trim=150 100 050 00, clip]{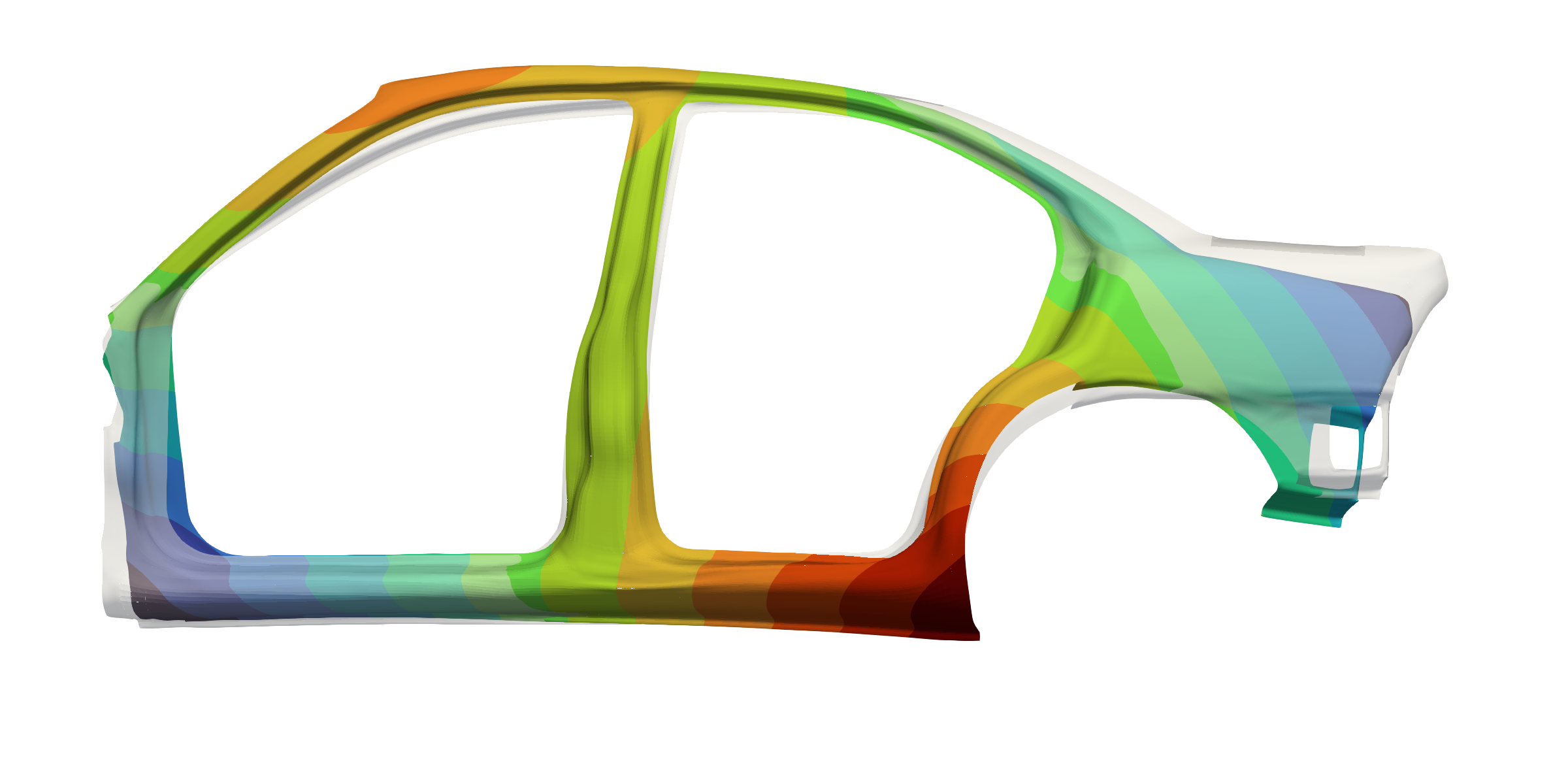}
		\end{minipage}
		\hfill
		\begin{minipage}{0.45\linewidth}
			\centering
				\includegraphics[width=\linewidth,trim=552 150 352 00, clip]{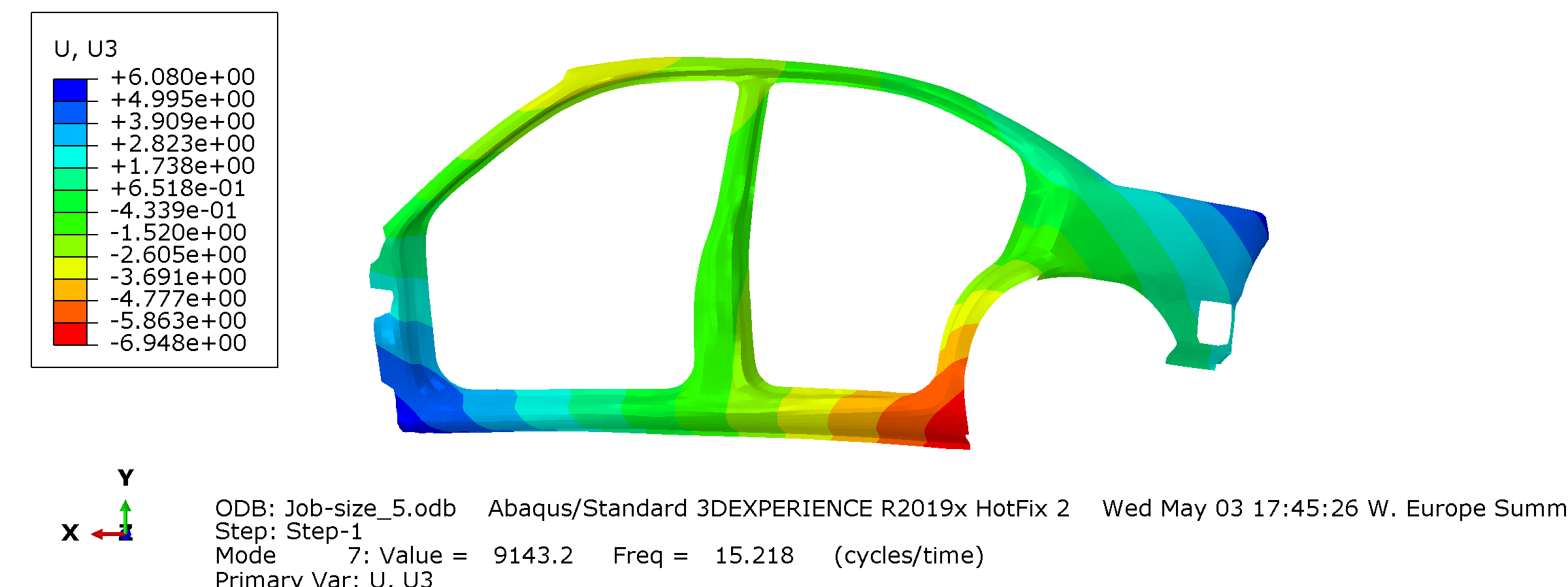}
		\end{minipage}
		\caption{Mode 1}
		\label{}
	\end{subfigure}
	
	\begin{subfigure}{\linewidth}
		\centering
		\begin{minipage}{0.45\linewidth}
			\centering
				\includegraphics[width=\linewidth,trim=150 100 050 00, clip]{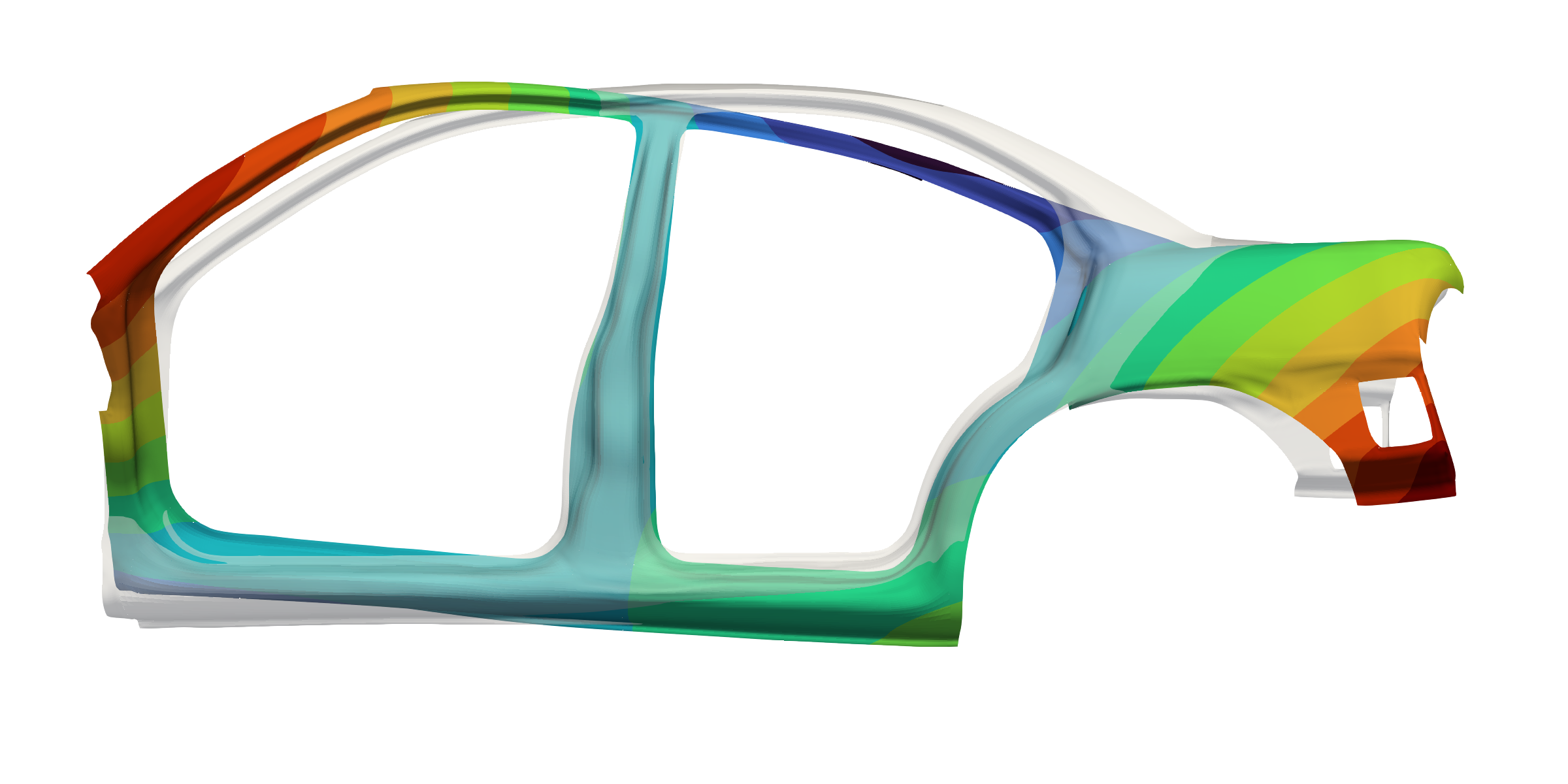}
		\end{minipage}
		\hfill
		\begin{minipage}{0.45\linewidth}
			\centering
				\includegraphics[width=\linewidth,trim=552 150 352 00, clip]{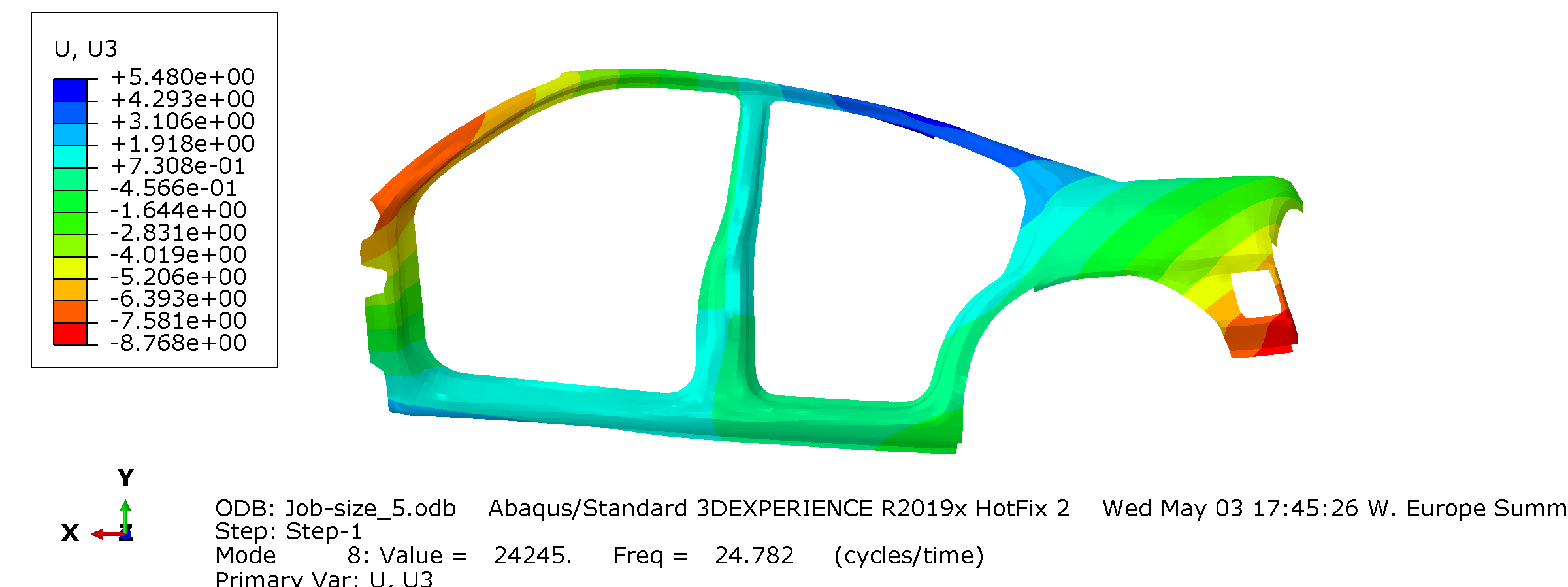}
		\end{minipage}
		\caption{Mode 2}
		\label{}
	\end{subfigure}
	
	\begin{subfigure}{\linewidth}
		\centering
		\begin{minipage}{0.45\linewidth}
			\centering
				\includegraphics[width=\linewidth,trim=150 100 050 00, clip]{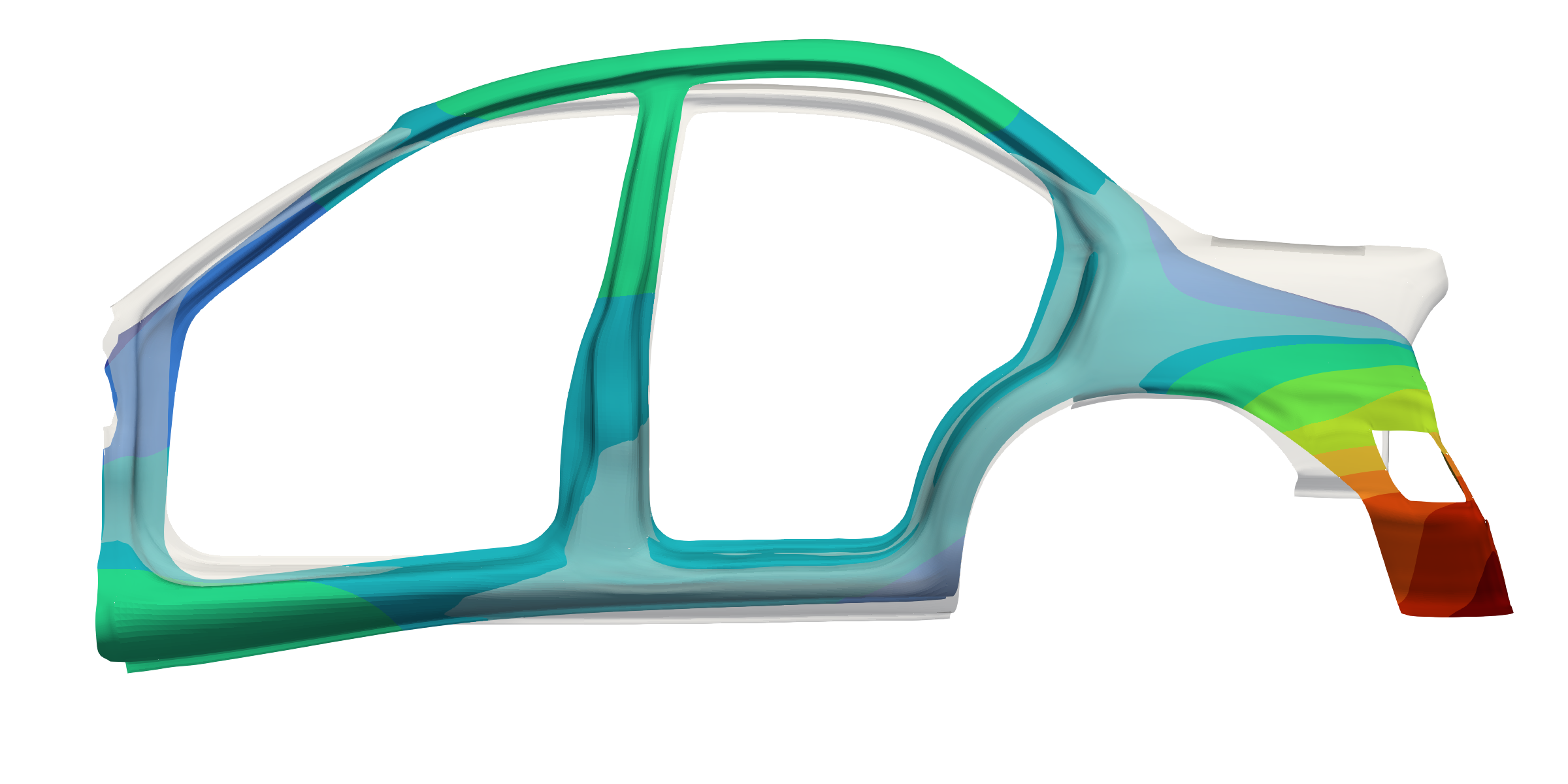}
		\end{minipage}
		\hfill
		\begin{minipage}{0.45\linewidth}
			\centering
				\includegraphics[width=\linewidth,trim=552 150 352 00, clip]{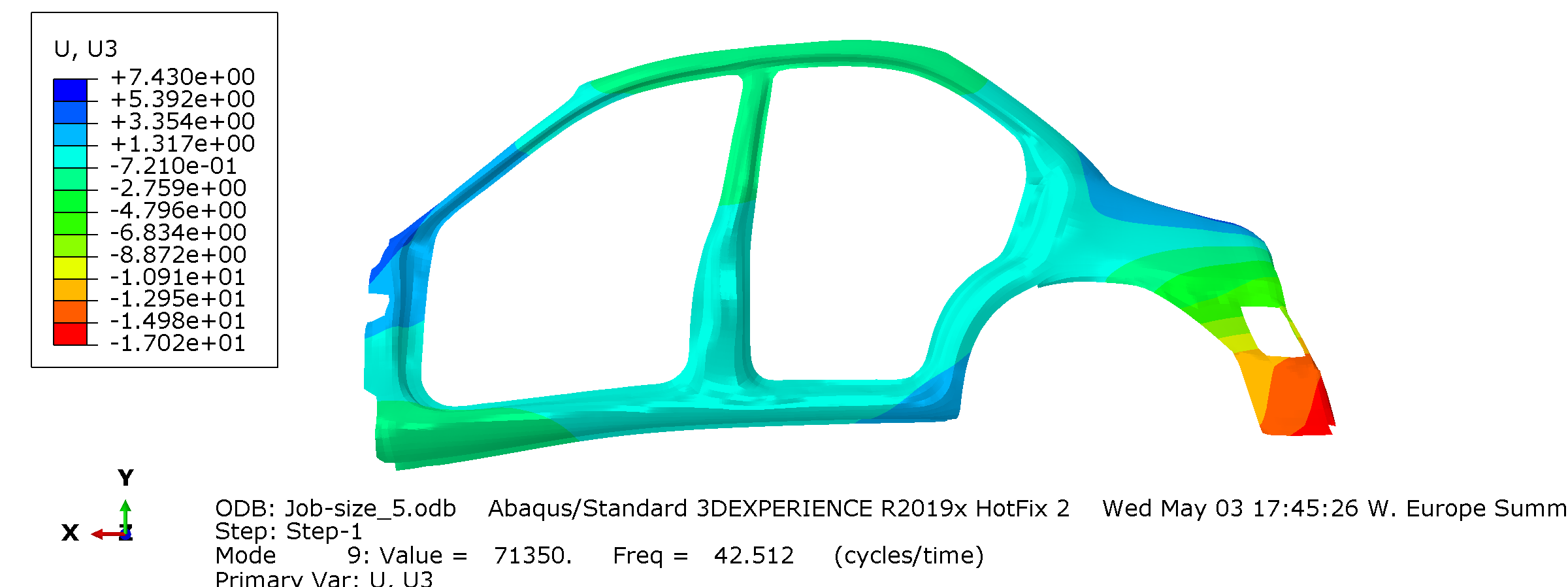}
		\end{minipage}
		\caption{Mode 3}
		\label{}
	\end{subfigure}
	
	\begin{subfigure}{\linewidth}
		\centering
		\begin{minipage}{0.45\linewidth}
			\centering
				\includegraphics[width=\linewidth,trim=150 100 050 00, clip]{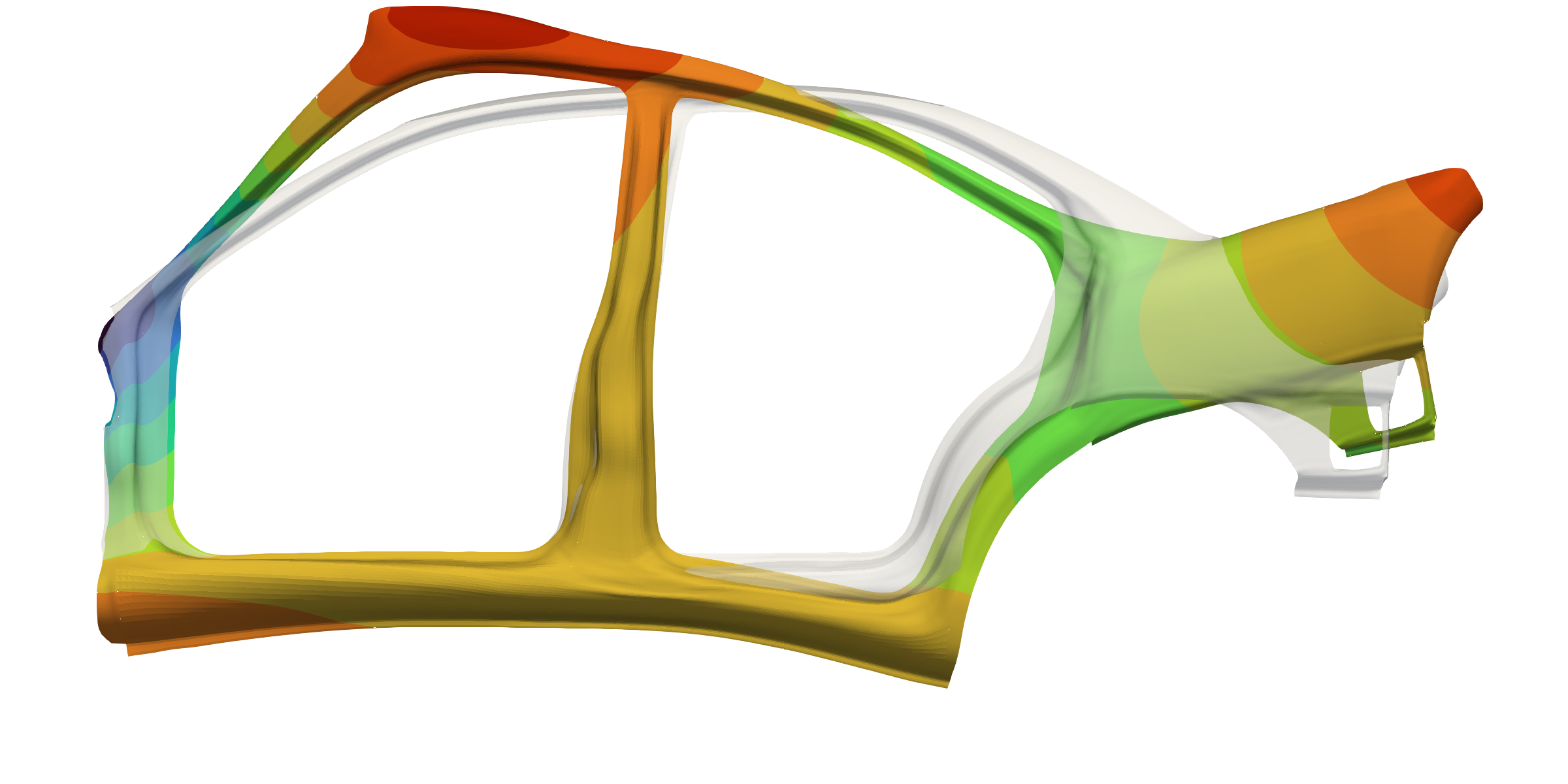}
		\end{minipage}
		\hfill
		\begin{minipage}{0.45\linewidth}
			\centering
				\includegraphics[width=\linewidth,trim=552 150 352 00, clip]{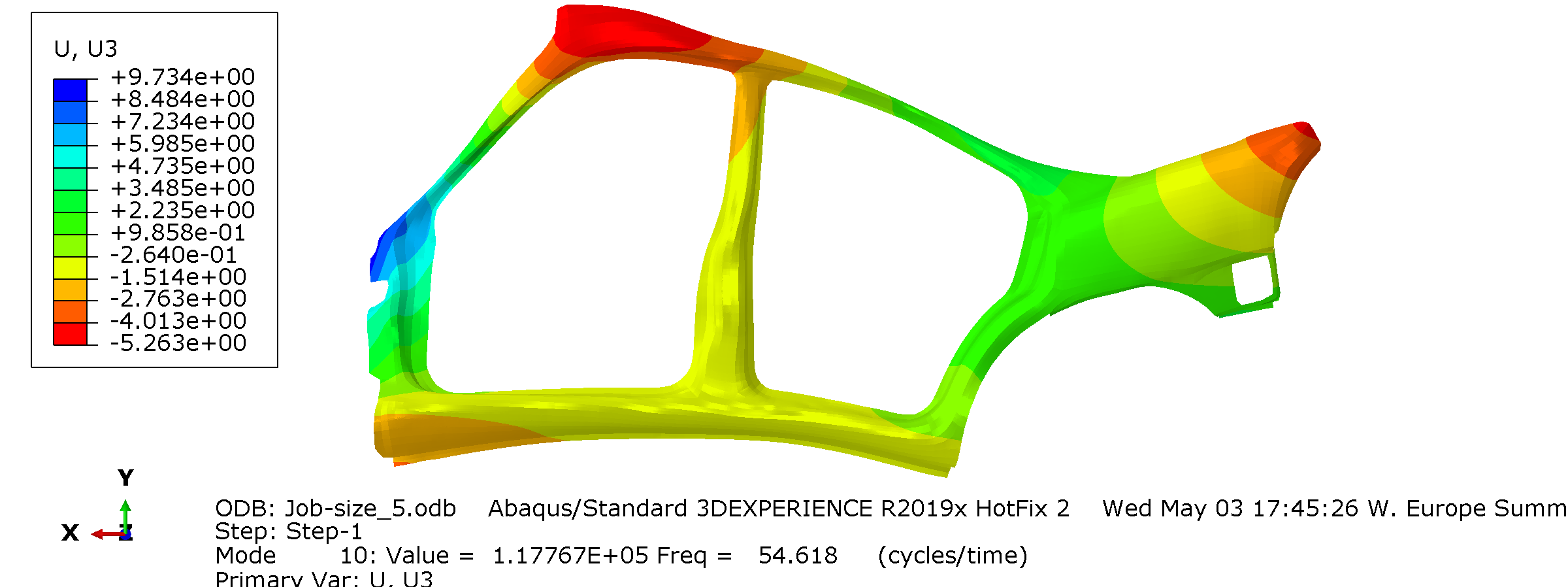}
		\end{minipage}
		\caption{Mode 4}
		\label{}
	\end{subfigure}
	
	
	
	\caption{Out-of-plane deformations of the first four vibration modes of the side of the car from \cref{fig:car_side}. The results on the left represent the results obtained by the D-Patch construction and the results on the right represent results obtained using ABAQUS (10mm). The mode shapes are all deformation modes warped by the deformation vector and plotted over the undeformed (transparent) geometry.}
	\label{fig:car_modes}
\end{figure}

\Cref{tab:ModalShell} presents the eigenfrequencies for the first four deformation modes of the car side panel for the D-Patch and the Almost-$C^1$ constructions with degree $p=2$ and regularity $r=1$ for the Almost-$C^1$ construction and with $(p,r)=(2,1)$, $(p,r)=(3,1)$ and $(p,r)=(3,2)$ for the D-Patch. \Cref{fig:car_modes} provides the corresponding mode shapes on the D-Patch geometry with $p=3$, $r=2$ and the mode shapes have been qualitatively matched to construct \cref{tab:ModalShell}. From these results, it can be observed that the Almost-$C^1$ and D-Patch methods provide eigenfrequencies in the same range and that the eigenfrequencies are mostly converging in the second digit. Moreover, the eigenfrequencies of the D-Patch and Almost-$C^1$ methods for coarse meshes and $p=2$, $r=1$ already provide reasonable estimates compared to higher degrees and refinements. On the other hand, the results obtained using an ABAQUS S4R element show convergence in the second digit, and slightly lower frequencies than the IGA results, possibly because the FEM uses a different geometry approximation. Overall, it can be concluded from this benchmark problem that the Almost-$C^1$ and D-Patch are more robust for industrial and large scale geometries, that are represented by at least $C^0$-conforming quadrilateral meshes, compared to the Approx.~$C^1$ and AS-$G^1$ methods due to the pre-processing efforts required by the latter. Furthermore, these methods are parameter-free, making them robust also with respect to penalty methods.

\begin{table}[]
  \centering
  \caption{Eigenfrequencies of the Almost-$C^1$ and D-Patch constructions for the car geometry in \cref{fig:car_side}. The results of an ABAQUS FEA simulation using the S4R element are provided as a reference. The mode-shapes are plotted in \cref{fig:car_modes}.}
  \label{tab:ModalShell}
  \begin{tabular}{lll|llll}
  \toprule
   Method & & \# DoFs & Mode 1 & Mode 2 & Mode 3 & Mode 4\\
   \midrule
    \multicolumn{2}{l}{\multirow[c]{4}{*}{Almost-$C^1$, $p=2$, $r=1$}}   & 13,731     & 15.740 & 25.567 & 43.829 & 56.654 \\
    \multicolumn{2}{l}{}                                                 & 49,758     & 15.762 & 25.564 & 43.429 & 56.778 \\
    \multicolumn{2}{l}{}                                                 & 189,654    & 15.776 & 25.552 & 43.269 & 56.785 \\
    \multicolumn{2}{l}{}                                                 & 740,814    & 15.774 & 25.531 & 43.177 & 56.746 \\
    \midrule
    \multicolumn{2}{l}{\multirow[c]{3}{*}{D-Patch, $p=2$, $r=1$}}        & 49,437     & 15.785 & 25.607 & 43.641 & 56.902\\
    \multicolumn{2}{l}{}                                                 & 189,333    & 15.780 & 25.561 & 43.323 & 56.807\\
    \multicolumn{2}{l}{}                                                 & 740,493    & 15.775 & 25.533 & 43.191 & 56.748 \\
    \midrule
    \multicolumn{2}{l}{\multirow[c]{2}{*}{D-Patch, $p=3$, $r=1$}}        & 136,839    & 15.749 & 25.593 & 43.348 & 56.786 \\
    \multicolumn{2}{l}{}                                                 & 630,459    & 15.760 & 25.581 & 43.231 & 56.801 \\
    \midrule
    \multicolumn{2}{l}{\multirow[c]{2}{*}{D-Patch, $p=3$, $r=2$}}        & 71,760     & 15.771 & 25.539 & 43.224 & 56.744 \\
    \multicolumn{2}{l}{}                                                 & 226,524    & 15.755 & 25.582 & 43.235 & 56.807 \\
    \midrule
    \multirow{3}*{ABAQUS S4R}& 10mm					                     & 126,966    & 15.303 & 24.881 & 42.629 & 54.887 \\
                             & 5mm 					                    & 440,076    & 15.224 & 24.780 & 42.516 & 54.627 \\
                             & 2.5mm 					                & 1,653,030 	& 15.119 & 24.640 & 42.338 & 54.277 \\

  \bottomrule
  \end{tabular}
\end{table}

\subsection{Stress analysis in a curved shell}\label{subsec:quantitative_stress}
An interesting application for smooth unstructured spline construction is for the use of thin shell analysis for engineering applications. Not only displacements (\cref{subsec:quantitative_linearShell}) or vibrations (\cref{subsec:quantitative_modalCar}) are of interest, but also stress evaluations, for example for fatigue analysis. In the last example, we demonstrate the performance of all methods in \cref{tab:comparison} on the evaluation of stresses in a curved Kirchhoff--Love shell. Since the Kirchhoff--Love shell formulation is displacement-based, the displacements are $C^1$ continuous across patch interfaces for $C^1$ constructions. The stresses, however, are based on the gradients of the displacements, hence their continuity theoretically is $C^0$ for a perfect $C^1$ coupling. In this example, we elaborate on the Von Mises membrane stress field resulting from the 6-patch elliptic paraboloid from \cref{fig:linearShell_paraboloid}. The stress fields are plotted for bases with degree and regularity from \cref{tab:comparison} and additionally for a basis with $p=4$, $r=2$. Note that the regularity $r$ of these bases is the regularity in the patch interior.\\

In \cref{fig:coupling_stresses}, the stress fields for the elliptic paraboloid example from \cref{fig:linearShell_paraboloid} are provided. From these results, it can immediately be seen that the stress field for a single patch parameterization with basis $p=2$, $r=1$ exposes the elements of the basis because of the $C^0$ continuity across elements. Similar effects are seen for the D-patch, Almost-$C^1$ and the penalty method. Increasing the degree of the basis while keeping the regularity the same results in a $p=3$, $r=1$ basis. The element continuity is still $C^0$ for the stresses, but the higher continuity of the basis within the element results in a slightly improved stress field, as can be seen from the single patch, the D-patch and penalty methods. The Approx.~$C^1$ and AS-$G^1$ methods in addition show a better stress field around the EVs compared to the D-patch with only small wiggles in the inner contour. Increasing the smoothness by going to $p=3$, $r=2$ shows that the Approx.~$C^1$ method predicts the stress field very well over the whole domain but with the wiggles in the inner contour, and that the D-patch suffers from the singularity at the EVs. Lastly, the $p=4$, $r=2$ plots show that the wiggles in the inner contour are eliminated for the Approx.~$C^1$ and the AS-$G^1$ methods and that the artifacts of the D-patch around the EV are still there but to a lesser extent. Finally, the results of the penalty method in \cref{fig:stresses_penalty} show it is able to provide an accurate representation of the stress fields. As seen from \cref{fig:linearShell_paraboloid}, penalty factors $\alpha=1$ and $\alpha=10$ provide good convergence in the bending energy norm. Indeed, the stress fields for the fixed $64\times64$ element meshes in \cref{fig:stresses_penalty} confirm that for these penalty factors the stress fields accurately represent the single patch stress fields, despite small artifacts around the EVs for $\alpha=1$. For a higher penalty factor of $\alpha=100$, the stress fields following from the penalty method are not guaranteed to be accurate, showing the downside of this method.\\

Overall, the stress analysis for multiple combinations and regularities shows that the Almost-$C^1$ method is generally unfavourable since it is only applicable for $p=2$, $r=1$ hence $C^0$ stress fields, suffering from a lack of continuity over the whole domain. This also makes the D-patch as applicable as the Approx.~$C^1$ method in terms of degree and regularity combinations. Comparing the D-patch with the Approx.~$C^1$ and the AS-$G^1$ methods, it is shown that the D-patch suffers from the singularity in the EVs when reconstructing stresses, whereas the other two methods are able to recover the stress fields without problems. Moreover, this example has also shown the advantage of smooth unstructured spline constructions for stress analyses, since their continuity across (almost) all of the domain is ensured, contrary to the penalty method. Lastly, this example shows the advantage of IGA in general over lower-order methods like FEA, since the higher-degree bases (e.g. $p=4$, $r=2$) provides smooth stress fields compared to lower-degree bases ($p=2$, $r=1$).\\

\input{Figures/Plots/StressParaboloid.tex}

\subsection{Conclusions}\label{subsec:conclusions}
In this section a quantitative comparison of the AS-$G^1$, the Approx.~$C^1$, the D-Patch and the Almost-$C^1$ constructions is provided. The methods have been assessed on different aspects: i) convergence of the biharmonic equation (\cref{subsec:quantitative_biharmonic}); ii) convergence of the linear Kirchhoff--Love shell (\cref{subsec:quantitative_linearShell}); iii) eigenvalue spectrum approximation (\cref{subsec:quantitative_eigenvalue}); iv) application to a large-scale complex geometry (\cref{subsec:quantitative_modalCar}) and; v) the reconstruction of stress fields (\cref{subsec:quantitative_stress}). From these analyses, the following conclusions can be drawn:
\begin{itemize}
 \item All methods converge in a theoretical setting to the same solution for the biharmonic equation (\cref{subsec:quantitative_biharmonic,subsec:quantitative_linearShell}). However, the convergence behaviour of the D-Patch method is sub-optimal and affected by conditioning issues for large meshes. Furthermore, the Approx.~$C^1$ and AS-$G^1$ methods give worse convergence compared to other methods for the hyperbolic paraboloid shell but good convergence rates for the elliptic paraboloid shell example.
 \item From a spectral analysis on the biharmonic equation \cref{subsec:quantitative_eigenvalue} it can be concluded that there is no best unstructured spline construction. Depending on the degree and regularity, small difference in the eigenvalue spectra are observed between the methods. Comparing with Nitsche's method, however, it is concluded that the unstructured spline constructions considered in this paper perform consistently better. This is also confirmed by the applied modal analysis on the car geometry \cref{subsec:quantitative_modalCar}, where penalty method fails to find accurate eigenfrequencies, possibly because of an unsuitable penalty parameter.
 \item From the applied modal analysis on a complex geometry, it can also be concluded that the Almost-$C^1$ and D-Patch constructions are more straight-forward to apply to a complex geometry extracted from a mesh. This is due to the fact that the Approx.~$C^1$ and AS-$G^1$ constructions require, respectively, a $G^2$ geometry and an analysis-suitable geometry, which are both not trivial to construct from an originally $C^0$-continuous mesh. Instead, the D-Patch and Almost-$C^1$ constructions require a $C^1$ geometry, which is easier to construct in general.
 \item From the stress fields presented in \cref{subsec:quantitative_stress} following from the analysis in \cref{subsec:quantitative_linearShell}, it can be concluded that the AS-$G^1$ and Approx.~$C^1$ methods provide excellent stress fields. The D-Patch also provides good stress fields, but inaccuracies are found around the EVs, possibly because of the singularity close to the EV. The Almost-$C^1$ method is considered inaccurate for stress analysis because of a lack of higher-degree generalizations. Lastly, comparison with penalty methods shows that the unstructured spline constructions generally provide a robust parameter-free approach for coupling, whereas the penalty method requires careful selection of the penalty parameter.
\end{itemize}
Overall, our finding suggest that the Almost-$C^1$ and D-Patch are generally easier to construct, but for certain problems they have limited accuracy. On the other hand, the AS-$G^1$ or Approx.~$C^1$ discretisations require more pre-processing efforts, but provide optimal convergence, hence accuracy. This, however, depends on the input geometry: generic quad-meshes might require more pre-processing efforts than $C^1$-matching parameterizations. Lastly, the results provided in this section have shown that strong coupling methods have certain advantages over weak methods, and therefore provide an interesting alternative.

\section{Conclusions and future work}\label{sec:conclusions}
In this paper, we provide a qualitative and quantitative comparison of unstructured spline constructions for smooth multi-patches in isogeometric analysis. The general advantage of unstructured spline constructions over trimming or variational coupling methods is that they are parameter-free, do not require specialized solvers and are typically once constructed in a shape optimization workflow. The goal of this paper is to compare the analysis-suitable $G^1$ (AS-$G^1$) the approximate $C^1$ (Approx.~$C^1$), the degenerate patches (D-Patch) and the Almost-$C^1$ constructions with respect to qualitative aspects (i.e. constraints for application) and quantitative aspects (i.e. numerical performance).\\

From the qualitative analysis, it followed that each method required a different set of constraints to be satisfied before the constructions could be applied, see \cref{fig:USconstraints,tab:qualitative}. Degree and regularity constraints can be satisfied by knot insertion routines or re-fitting, which are relatively straight-forward. The constraint on analysis-suitability for the AS-$G^1$ and the constraint on $G^2$ continuity for the Approx.~$C^1$ method require dedicated reparameterization routines, such as the one presented by \cite{Kapl2018}. The fact that D-Patches are restricted to geometries without boundary extraordinary vertices requires redefinition of the quadrilateral mesh. Lastly, the fact that the Almost-$C^1$ method is only defined for bi-quadratic bases ($p=2$) restricts the inter-element continuity to $C^1$ through the whole domain. Depending on the application and the availability of existing routines in software, different unstructured spline constructions are favourable, depending on the geometric flexibility or desired degree and regularity .\\

From the quantitative analysis, some conclusions can be drawn on the considered unstructured spline constructions and between unstructured spline constructions compared to variational methods such as Nitsche's method or a penalty method. From the analysis, it was in general observed that depending on the problem type, the different methods have their advantages and disadvantages. Firstly, simple biharmonic equations (see \cref{subsec:quantitative_biharmonic}) and linear shells (see \cref{subsec:quantitative_linearShell}) provided good results for all methods. However, the AS-$G^1$ and Approx.~$C^1$ methods showed slow convergence for the double-curved shell and the D-patch suffered from ill-conditioning for fine meshes. The Almost-$C^1$ provided good results in general, however it is only applicable on bi-quadratic splines. Secondly, all methods showed superiority over Nitsche's method for the computation of an eigenvalue spectrum for plate vibrations (see \cref{subsec:quantitative_eigenvalue}) and no significant differences between the unstructured spline constructions have been observed. Thirdly, the D-Patch and Almost-$C^1$ showed straight-forward applicability on the problem of a complex geometry (see \cref{subsec:quantitative_modalCar}), whereas the analysis-suitability requirement of the AS-$G^1$ method and the smoothness requirement of the Approx.~$C^1$ method are non-trivial to satisfy on off-the-shelf industrial geometries.\\

For the penalty method, no suitable penalty parameter was found, and probably optimal penalty parameters should be chosen per interface rather than globally. Lastly, the AS-$G^1$ and Approx.~$C^1$ methods provided superior results for stress reconstruction, where the D-Patch suffered around the EVs due to its singular parameterization and the Almost-$C^1$ method provided bad results due to a lack of higher degrees.\\

In conclusion, both comparisons give an overview of the applicability of the methods with respect to the requirements needed to construct them, on the notions of nestedness and in general on the performance of the methods. Overall, it can be concluded from both analyses that among the compared methods, there is no general best construction. More precisely, the quantitative analysis shows that different methods perform differently in different applications, given that they can be constructed. Furthermore, with the backgrounds and properties provided in the qualitative analysis section, we hope that the present paper provides valuable insights for application of the considered methods to multi-patch problems.\\

In addition, the comparisons in the present paper give directions for the improvements of the considered methods. For the AS-$G^1$ and Approx.~$C^1$ methods, restrictions on geometry and parameterization are a bottleneck in the industrial applications. Therefore, it is recommended to expand the applicability of these methods by developing dedicated geometric pre-processing routines. For the D-Patch construction, the limitation of the construction of the basis near $\nu>3$ boundary EVs calls for the development of routines to eliminate these EVs in quadrilateral multi-patches, as discussed in the qualitative comparison. Furthermore, the example of the bi-harmonic equation has shown that the D-Patch can suffer from ill-conditioned system, hence development of pre-conditioners for D-Patch constructions is advised. Lastly, although the Almost-$C^1$ resolves the downsides of the D-Patch construction, its restriction on the degree of the spline-space is a major disadvantage when plotting stress fields in shell analysis. Therefore, for the Almost-$C^1$ construction it is recommended to explore expansion to higher degrees.


\section*{Acknowledgements}
HMV is grateful to Delft University of Technology, specifically the faculty of Electrical Engineering, Mathematics and Computer Science (EEMCS) and the faculty of Mechanical, Maritime and Materials Engineering (3mE) for the financial support to conduct this research. AM acknowledges support from EU's Horizon 2020 research under the Marie Sklodowska-Curie grant 860843 and DT is grateful to ANSYS Inc. for their financial support. Furthermore, the authors thank Wei Jun Wong from Delft University of Technology for delivering the ABAQUS reference results, and the community of the Geometry + Simulation Modules (a.k.a. gismo) for their continuous effort on improving the code.
\bibliography{references}

\end{document}